\documentclass[10pt,a4paper]{article}
\usepackage{amsmath,amssymb,amsthm,graphics,graphicx,epic,eepic,multicol,ascmac}
\usepackage[mathscr]{euscript}
\usepackage[all]{xy}

\numberwithin{equation}{section}
\theoremstyle{theorem}
\newtheorem{thm}{Theorem}[section]
\newtheorem{prop}[thm]{Proposition}
\newtheorem{lem}[thm]{Lemma}
\newtheorem{rem}[thm]{Remark}

\newtheorem{ex}[thm]{Example}

\theoremstyle{definition}
\newtheorem{defn}[thm]{Definition}

\def\al{\alpha}

\def\ep{\epsilon}
\def\wht(#1){\widehat{\ #1\ }}

\newcommand{\cA}{{\mathcal A}}

\newcommand{\cF}{{\mathcal F}}

\newcommand{\cY}{{\mathcal Y}}

\newcommand{\frg}{\mathfrak g}
\newcommand{\frh}{\mathfrak h}

\newcommand{\frl}{\mathfrak l}

\newcommand{\frn}{\mathfrak n}
\newcommand{\fro}{\mathfrak o}
\newcommand{\frp}{\mathfrak p}

\newcommand{\frs}{\mathfrak s}

\newcommand{\bbC}{\mathbb C}

\newcommand{\ch}{\mathrm{ch}}

\newcommand{\lbr}{\begin{bmatrix}}
\newcommand{\rbr}{\end{bmatrix}}

\newcommand{\cd}{commutative diagram }

\def\ge{\frg}

\def\al{\alpha}

\def\beneme{\begin{enumerate}}
\def\beq{\begin{equation}}
\def\beqn{\begin{eqnarray}}
\def\beqnn{\begin{eqnarray*}}

\def\bfii0{{\bf i_0}}

\def\bbra#1,#2,#3{\left\{\begin{array}{c}\hspace{-5pt}
#1;#2\\ \hspace{-5pt}#3\end{array}\hspace{-5pt}\right\}}
\def\cd{\cdots}
\def\ci(#1,#2){c_{#1}^{(#2)}}
\def\Ci(#1,#2){C_{#1}^{(#2)}}
\def\mpp(#1,#2,#3){#1^{(#2)}_{#3}}
\def\bCi(#1,#2){\ovl C_{#1}^{(#2)}}
\def\ch(#1,#2){c_{#2,#1}^{-h_{#1}}}
\def\cc(#1,#2){c_{#2,#1}}

\def\del{\delta}
\def\Del{\Delta}

\def\di(#1,#2){D_{#1}^{(#2)}}
\def\dbi(#1,#2){\ovl D_{#1}^{(#2)}}

\def\eneme{\end{enumerate}}
\def\ep{\epsilon}
\def\eeq{\end{equation}}
\def\eeqn{\end{eqnarray}}
\def\eeqnn{\end{eqnarray*}}

\def\gau#1,#2{\left[\begin{array}{c}\hspace{-5pt}#1\\
\hspace{-5pt}#2\end{array}\hspace{-5pt}\right]}

\def\ji(#1,#2){j_{#1}^{(#2)}}

\def\lan{\langle}

\def\lm{\lambda}
\def\Lm{\Lambda}

\def\nd{\noindent}

\def\ovl{\overline}

\def\qq{\qquad}
\def\q{\quad}
\def\qed{\hfill\framebox[2mm]{}}

\def\ran{\rangle}

\def\TY(#1,#2,#3){#1^{(#2)}_{#3}}

\def\xxi(#1,#2,#3){\displaystyle {}^{#1}\Xi^{(#2)}_{#3}}
\def\xsi(#1,#2,#3){\displaystyle {}^{#1}\Sigma^{(#2)}_{#3}}
\def\xE(#1,#2,#3){\displaystyle {}^{#1}E_{#2}[#3]}
\def\xF(#1,#2){\displaystyle {}^{#1}F_{#2}}
\def\xx(#1,#2){\displaystyle {}^{#1}\Xi_{#2}}
\def\W1{W(\varpi_1)}

\def\m@th{\mathsurround=0pt}
\def\fsquare(#1,#2){
\hbox{\vrule$\hskip-0.4pt\vcenter to #1{\normalbaselines\m@th
\hrule\vfil\hbox to #1{\hfill$\scriptstyle #2$\hfill}\vfil\hrule}$\hskip-0.4pt
\vrule}}

\newcommand{\ba}{\begin{array}}
\newcommand{\ea}{\end{array}}

\newcommand{\eq}{\begin{eqnarray}}
\newcommand{\eneq}{\end{eqnarray}}

\title{\textbf{\large{Cluster algebras of finite type via Coxeter elements and Demazure Crystals of type B,C,D}}}
\author{\normalsize{YUKI KANAKUBO\thanks{Division of Mathematics, 
Sophia University, Kioicho 7-1, Chiyoda-ku, Tokyo 102-8554,
Japan: {j\_chi\_sen\_you\_ky@eagle.sophia.ac.jp}}}
}
\date{}

\begin{document}

\maketitle
\vspace{-10pt}

\begin{abstract}
For a classical group $G$ and a Coxeter element $c$ of the Weyl group, it is known that the coordinate ring $\mathbb{C}[G^{e,c^2}]$ of the double Bruhat cell $G^{e,c^2}:=B\cap B_-c^2B_-$ has a structure of cluster algebra of finite type, where $B$ and $B_-$ are opposite Borel subgroups. In this article, we consider the case $G$ is of type ${\rm B}_r$, ${\rm C}_r$ or ${\rm D}_r$ and describe all the cluster variables in $\mathbb{C}[G^{e,c^2}]$ as monomial realizations of certain Demazure crystals.
\end{abstract}

\nd

%Information for first author

\section{Introduction}

Fomin and Zelevinsky have invented cluster algebras for the study of total positivity and dual semicanonical bases (\cite{FZ2}). It is a commutative ring generated by so-called ``cluster 
variables''. It is known that the coordinate rings of many algebraic varieties related to semisimple algebraic groups carry cluster algebra structures. For instance, in \cite{A-F-Z,GY,GY2}, for simply connected, connected, complex simple algebraic group $G$ and 
its Weyl group elements $u,v\in W$, it is shown that $\mathbb{C}[G^{u,v}]$ is a cluster algebra, 
where $G^{u,v}:=BuB\cap B_-vB_-$ and $B$, $B_-$ are opposite
Borel subgroups. In \cite{GLS}, it is proved that the coordinate rings $\mathbb{C}[N(w)]$ and $\mathbb{C}[N^w]$ are cluster algebras by using the additive categorification via finite dimensional modules of the preprojective
algebras, where $N(w):=N\cap (w^{-1}N_-w), N^w:=N\cap (B_-w B_-)$ and $N$, $N_-$ are unipotent radicals. It is also proved that 
all the cluster variables are included in the dual semicanonical basis in the coordinate rings.

A cluster algebra is said to be of {\it finite type} if it has only finitely many cluster variables. In \cite{FZ3}, a complete classification of the cluster algebras of finite type are provided. More precisely, they are classified by the set of Cartan matrices up to coefficients. For a fixed Cartan matrix, all the cluster variables are parametrized by 
the set of ``almost positive roots'', which is, a union of all positive roots and negative simple roots corresponding to the Cartan matrix. By this classification, the {\it type} of each cluster algebra of finite type can be defined as the Cartan-Killing type of the corresponding Cartan matrix. Let $c\in W$ be a Coxeter element such that the length $l(c)$ satisfies 
$l(c^2)=2 l(c)=2$rank$(G)$. It is known that one can realize a cluster algebra of finite type on the coordinate ring $\bbC[G^{e,c^2}]$ and its type coincides with the Cartan-Killing type of $G$ \cite{A-F-Z}.

One purpose of our study is to reveal relation between cluster variables of the coordinate rings and Kashiwara's crystal bases (\cite{K0,K2}). The crystal bases were introduced for combinatorial study of the integrable modules over quantum groups and have many realizations, e.g., tableaux, paths, monomials, etc. In this article, we will treat the {\it monomial realization}, which is defined in \cite{K, Nj}.

In \cite{KaN}, we treated the {\it initial cluster variables} denoted by $\Delta(k;{\rm \bf{i}})$ of $\bbC[G^{u,e}]$ ($1\leq k\leq l(u)-r$, $u\in W$ and ${\rm \bf{i}}$ is a reduced word of $u$) in the case $G={\rm SL}_{r+1}(\mathbb{C})$. We found explicit formulas for $\{\Delta(k;{\rm \bf{i}})\}_{1\leq k\leq l(u)-r}$, which express them by Laurent polynomials with coefficients $1$. We also proved that the set of monomials appearing in $\Delta(k;{\rm \bf{i}})$ coincides with a monomial realization of certain Demazure crystal. In \cite{KaN2}, we considered the case $G$ is a classical algebraic group of type ${\rm B}_r$, ${\rm C}_r$  or ${\rm D}_r$, and gave explicit formulas for a part of the initial cluster variables $\{\Delta(k;{\rm \bf{i}})\}$ in $\bbC[G^{u,e}]$. Just as in the case $G={\rm SL}_{r+1}(\mathbb{C})$, the set of monomials appearing in $\Delta(k;{\rm \bf{i}})$ coincides with a monomial realization of certain Demazure crystal. In the both papers, we did not treat all the cluster variables but a part of the cluster variables. On the other hand, in \cite{KaN3}, we considered the case $G={\rm SL}_{r+1}(\mathbb{C})$ $(r\geq 3)$ and all the cluster variables in the coordinate ring $\bbC[G^{e,c^2}]$ ($c$ is a Coxeter element). As mentioned above, the algebra $\bbC[G^{e,c^2}]$ is a cluster algebra of finite type. We described each cluster variable $\varphi$ as a Laurent polynomial with coefficients $1$ and showed that the set of monomials appearing in $\varphi$ coincides with a monomial realization of the direct sum of certain Demazure crystals.

In this article, we consider the case $G$ is a classical algebraic group of type ${\rm B}_r$, ${\rm C}_r$ or ${\rm D}_r$ and the coordinate ring $\bbC[G^{e,c^2}]$, where $c=(s_r\cdots s_2s_1)$ is a Coxeter element. Our main result is that all the cluster variables in $\bbC[G^{e,c^2}]$ are described as Laurent polynomials with positive integers, and forgetting the coefficients, the set of monomials appearing in each cluster variable coincides with a monomial realization of the direct sum of certain Demazure crystals. 

For example, let us consider the case $G={\rm SO}_{5}(\mathbb{C})$ (type ${\rm B}_2$ algebraic group). Monomial realizations of the crystals  $B(\Lambda_1)$ and $B(\Lambda_2)$ of type ${\rm B}_2$ are
\[
Y_{1,1}\overset{1}{\longrightarrow} \frac{Y^2_{2,2}}{Y_{2,1}}
\overset{2}{\longrightarrow} \frac{Y_{2,2}}{Y_{3,2}}
\overset{2}{\longrightarrow} \frac{Y_{2,1}}{Y^2_{3,2}}
\overset{1}{\longrightarrow} \frac{1}{Y_{3,1}}\ \ \ {\rm and}\ \ \ 
Y_{1,2} \overset{2}{\longrightarrow} \frac{Y_{1,1}}{Y_{2,2}}
\overset{1}{\longrightarrow} \frac{Y_{2,2}}{Y_{2,1}}
\overset{2}{\longrightarrow} \frac{1}{Y_{3,2}},
\]
respectively. On the other hand, taking a Coxeter element $c=s_2s_1\in W$, specific initial cluster variables in 
$\mathbb{C}[G^{e,c^2}]$ are given by the {\it generalized minors} $\Delta_{\Lambda_i,s_1s_2\Lambda_i}$ $(i=1,2)$ (see \ref{bilingen}). Using the biregularly isomorphism $\ovl{x}^G_{\textbf{i}}:H\times (\mathbb{C}^{\times})^{4}\rightarrow G^{e,c^2}$ ($\textbf{i}:=(2,1,2,1)$) in Proposition
 \ref{gprime}, we have
\[ 
\Delta_{\Lambda_1,s_1s_2\Lambda_1}\circ \ovl{x}^G_{\textbf{i}}(a;\textbf{Y})=
a^{\Lm_1}(Y_{1,1}+\frac{Y^2_{2,2}}{Y_{2,1}}),\]
\[
\Delta_{\Lambda_2,s_1s_2\Lambda_2}\circ \ovl{x}^G_{\textbf{i}}(a;\textbf{Y})=
a^{\Lm_2}(Y_{1,2}+\frac{Y_{1,1}}{Y_{2,2}}+\frac{Y_{2,2}}{Y_{2,1}}),
\]
where we set $a\in H$ and $\textbf{Y}:=(Y_{1,2},Y_{1,1},Y_{2,2},Y_{2,1})\in(\mathbb{C}^{\times})^{4}$. Comparing with the above crystal graphs of $B(\Lambda_1)$ and $B(\Lambda_2)$,
we see that the set of monomials $\{Y_{1,1},\frac{Y^2_{2,2}}{Y_{2,1}}\}$ (resp. $\{Y_{1,2},\frac{Y_{1,1}}{Y_{2,2}},\frac{Y_{2,2}}{Y_{2,1}}\}$) appearing in $\Delta_{\Lambda_1,s_1s_2\Lambda_1}\circ \ovl{x}^G_{\textbf{i}}(a;\textbf{Y})$ 
(resp. $\Delta_{\Lambda_2,s_1s_2\Lambda_2}\circ \ovl{x}^G_{\textbf{i}}(a;\textbf{Y})$)
coincides with the monomial realization of the Demazure crystal $B(\Lambda_1)_{s_1}$ (resp. $B(\Lambda_2)_{s_2s_1}$) (see \ref{Demcrysub}). 

All other cluster variables in $\mathbb{C}[G^{e,c^2}]$ are
\[ \Delta_{\Lm_1,s_2s_1\Lm_1}\circ \ovl{x}^G_{\textbf{i}}=a^{\Lm_1}Y_{2,1},\quad (\Delta_{s_2\Lm_2,s_2s_1s_2\Lm_2}\cdot\Delta_{\Lm_2,\Lm_2})\circ \ovl{x}^G_{\textbf{i}}=a^{\Lm_1} Y_{2,2},\]
\[ \Delta_{\Lm_2,s_2\Lm_2}\circ \ovl{x}^G_{\textbf{i}}=a^{\Lm_2}(Y_{1,2}Y_{2,1}+
\frac{Y_{1,1}Y_{2,1}}{Y_{2,2}}),\]
\begin{eqnarray*}
& &(\Delta_{s_1s_2\Lm_2,s_2s_1s_2\Lm_2}\cdot \Delta_{\Lm_2,s_1s_2\Lm_2}\cdot\Delta_{\Lm_1,\Lm_1}-\Delta_{\Lm_2,\Lm_2}\cdot\Delta_{\Lm_2,s_2s_1s_2\Lm_2})\circ \ovl{x}^G_{\textbf{i}}\\
& &=a^{2\Lm_2}(Y^2_{1,2}Y_{2,1}+2\frac{Y_{1,1}Y_{1,2}Y_{2,1}}{Y_{2,2}}
+\frac{Y^2_{1,1}Y_{2,1}}{Y^2_{2,2}}+Y_{1,1}),
\end{eqnarray*}
and the sets of monomials appearing in these Laurent polynomials coincide with monomial realizations of the Demazure crystals $B(\Lambda_1)_e$, $B(\Lambda_2)_e$, $B(\Lambda_1+\Lambda_2)_{s_2}$ and $B(\Lambda_1+2\Lambda_2)_{s_2}\oplus B(\Lm_1)_e$, respectively. Note that the set of almost positive roots of type ${\rm B}_2$ is $\{-\al_1,-\al_2,\al_1,\al_2,\al_1+\al_2,\al_1+2\al_2\}$. Therefore, the number of the cluster variables in $\mathbb{C}[G^{e,c^2}]$ is $6$ from the result of \cite{FZ3}. 

The article is organized as follows. In section 2, we review the explicit forms of fundamental
representations of classical groups. Section 3 is devoted to recall properties of double Bruhat cells. In section 4, 
after a concise reminder on cluster algebras, we review isomorphisms between the coordinate rings of double Bruhat cells and cluster algebras $\mathcal{A}(\textbf{i})$. In section 5, we shortly review 
the monomial realizations of crystal bases. Section 6 presents our main results, which provide a relation between all the cluster variables in $\mathbb{C}[G^{e,c^2}]$ and monomial 
realizations of Demazure crystals, and we prove them in section 7.

\vspace{2mm}

\nd \textbf{Acknowledgement.} I would like to thank T. Nakashima for his helpful comments and discussions.

\section{Fundamental representations}\label{SectFund}

First, we review the fundamental representations of the complex simple Lie algebras $\ge$ of type ${\rm A}_r$,  ${\rm B}_r$,  ${\rm C}_r$,  and  ${\rm D}_r$ \cite{KN, N1} for calculations of generalized minors (see Subsection \ref{bilingen}). 
Let $I:=\{1,\cdots,r\}$, $A=(a_{ij})_{i,j\in I}$ 
be the Cartan matrix of $\ge$, and $(\frh,\{\al_i\}_{i\in I},\{h_i\}_{i\in I})$ the associated
root data 
satisfying $\al_j(h_i)=a_{ij}$, where 
$\al_i\in \frh^*$ is a simple root and 
$h_i\in \frh$ is a simple co-root.
Let $\{\Lm_i\}_{i\in I}$ be the set of the fundamental 
weights satisfying $\Lm_i(h_j)=\del_{i,j}$, $P=\bigoplus_{i\in I}\mathbb{Z}\Lm_i$ the weight lattice and $P^*=\bigoplus_{i\in I}\mathbb{Z}h_i$ the dual weight lattice. 

\subsection{Type ${\rm A}_r$}\label{SectFundA}

Let $\frg=\frs\frl(r+1,\mathbb{C})$ be the simple Lie algebra of type ${\rm A}_r$. The Cartan matrix $A=(a_{i,j})_{i,j\in I}$ of $\frg$ is as follows:
\[a_{i,j}=
\begin{cases}
2 & {\rm if}\ i=j, \\
-1 & {\rm if}\ |i-j|=1, \\
0 & {\rm otherwise.}  
\end{cases}
\]
For $\frg=\lan \frh,e_i,f_i(i\in I)\ran$, 
let us describe the vector representation 
$V(\Lm_1)$. Set ${\mathbf B}^{(r)}:=
\{v_i|\ i=1,2,\cd,r+1\}$ and define 
$V(\Lm_1):=\bigoplus_{v\in{\mathbf B}^{(r)}}\bbC v$. The weights of $v_i$ $(i=1,\cd,r+1)$ are given by ${\rm wt}(v_i)=\Lm_i-\Lm_{i-1}$, where $\Lm_0=\Lm_{r+1}=0$. We define the $\frg$-action on $V(\Lm_1)$ as follows: For $i\in I$ and $j$ $(1\leq j\leq r+1)$,
\[
h v_j=\lan h,{\rm wt}(v_j)\ran v_j\ \ (h\in P^*), \qquad
f_iv_i=v_{i+1},\ \ 
e_iv_{i+1}=v_i,
\]
and the other actions are trivial.

Let $\Lm_i$ be the $i$-th fundamental weight of type ${\rm A}_r$.
As is well-known that the fundamental representation 
$V(\Lm_i)$ $(1\leq i\leq r)$
is embedded in $\wedge^i V(\Lm_1)$
with multiplicity free.
The explicit form of the highest (resp. lowest) weight 
vector $u_{\Lm_i}$ (resp. $v_{\Lm_i}$)
of $V(\Lm_i)$ is realized in 
$\wedge^i V(\Lm_1)$ as follows:
\begin{equation}\label{A-h-l}
u_{\Lm_i}=v_1\wedge v_2\wedge\cdots\wedge v_i,\qq
v_{\Lm_i}=v_{i+1}\wedge v_{i+2}\wedge\cdots\wedge v_{r+1}.
\end{equation}

\subsection{Type ${\rm C}_r$}\label{SectFundC}

Let $\frg=\frs\frp(2r,\mathbb{C})$ be the simple Lie algebra of type ${\rm C}_r$. The Cartan matrix $A=(a_{i,j})_{i,j\in I}$ of $\frg$ is given by
\[a_{i,j}=
\begin{cases}
2 & {\rm if}\ i=j, \\
-1 & {\rm if}\ |i-j|=1\ {\rm and}\ (i,j)\neq (r-1,r), \\
-2 & {\rm if}\ (i,j)=(r-1,r), \\
0 & {\rm otherwise.}  
\end{cases}
\]
Note that $\al_i\ (i\neq r)$ are short roots and $\al_r$ is the long simple root. 

Define the total order on the set $J_{{\rm C}}:=\{i,\ovl i|1\leq i\leq r\}$ by 
\begin{equation}\label{C-order}
 1< 2<\cd< r-1< r
< \ovl r< \ovl{r-1}< \cd< \ovl
 2< \ovl 1.
\end{equation}
For $\frg=\lan \frh,e_i,f_i(i\in I)\ran$, 
let us describe the vector representation 
$V(\Lm_1)$. Set ${\mathbf B}^{(r)}:=
\{v_i,v_{\ovl i}|i=1,2,\cd,r\}$ and define 
$V(\Lm_1):=\bigoplus_{v\in{\mathbf B}^{(r)}}\bbC v$. The weights of $v_i$, $v_{\ovl{i}}$ $(i=1,\cd,r)$ are given by ${\rm wt}(v_i)=\Lm_i-\Lm_{i-1}$ and ${\rm wt}(v_{\ovl{i}})=\Lm_{i-1}-\Lm_{i}$,
where $\Lm_0=0$. We define the $\frg$-action on $V(\Lm_1)$ as follows:
\begin{eqnarray}
&& h v_j=\lan h,{\rm wt}(v_j)\ran v_j\ \ (h\in P^*,\ j\in J_{{\rm C}}), \label{C-wtv}\\
&&f_iv_i=v_{i+1},\ f_iv_{\ovl{i+1}}=v_{\ovl i},\q
e_iv_{i+1}=v_i,\ e_iv_{\ovl i}=v_{\ovl{i+1}}
\q(1\leq i<r),\label{C-f1}\\
&&f_r v_r=v_{\ovl r},\qq 
e_r v_{\ovl r}=v_r,\label{C-f2}
\end{eqnarray}
and the other actions are trivial.

Let $\Lm_i$ be the $i$-th fundamental weight of type ${\rm C}_r$.
As is well-known that the fundamental representation 
$V(\Lm_i)$ $(1\leq i\leq r)$
is embedded in $\wedge^i V(\Lm_1)$
with multiplicity free.
The explicit form of the highest (resp. lowest) weight 
vector $u_{\Lm_i}$ (resp. $v_{\Lm_i}$)
of $V(\Lm_i)$ is realized in 
$\wedge^i V(\Lm_1)$ as follows:
\begin{equation}\label{C-h-l}
u_{\Lm_i}=v_1\wedge v_2\wedge\cdots\wedge v_i,\quad
v_{\Lm_i}=v_{\ovl{1}}\wedge v_{\ovl{2}}\wedge\cdots \wedge v_{\ovl{i}}.
\end{equation}

\subsection{Type ${\rm B}_r$}\label{SectFundB}

Let $\frg=\frs\fro(2r+1,\mathbb{C})$ be the simple Lie algebra of type ${\rm B}_r$. The Cartan matrix $A=(a_{i,j})_{i,j\in I}$ of $\frg$ is given by
\[a_{i,j}=
\begin{cases}
2 & {\rm if}\ i=j, \\
-1 & {\rm if}\ |i-j|=1\ {\rm and}\ (i,j)\neq (r,r-1), \\
-2 & {\rm if}\ (i,j)=(r,r-1), \\
0 & {\rm otherwise.}  
\end{cases}
\]
Note that $\al_i\ (i\neq r)$ are long roots and $\al_r$ is the short simple root. 

Define the total order on the set $J_{{\rm B}}:=\{i,\ovl i|1\leq i\leq r\}\cup\{0\}$ by 
\begin{equation}\label{B-order}
 1< 2<\cd< r-1< r< 0<
 \ovl r< \ovl{r-1}< \cd< \ovl
 2< \ovl 1.
\end{equation}
For $\frg=\lan \frh,e_i,f_i(i\in I)\ran$, 
let us describe the vector representation 
$V(\Lm_1)$. Set ${\mathbf B}^{(r)}:=
\{v_i,v_{\ovl i}|i=1,2,\cd,r\}\cup\{v_0\}$ and define 
$V(\Lm_1):=\bigoplus_{v\in{\mathbf B}^{(r)}}\bbC v$. The weights of $v_i$, $v_{\ovl{i}}$ $(i=1,\cd,r)$ and $v_0$ are as follows:
\begin{equation}\label{B-wtv}
 {\rm wt}(v_i)=\Lm_i-\Lm_{i-1},\q {\rm wt}(v_{\ovl{i}})=\Lm_{i-1}-\Lm_{i}\q (1\leq i\leq r-1),
\end{equation}
\[  {\rm wt}(v_r)=2\Lm_r-\Lm_{r-1},\q {\rm wt}(v_{\ovl{r}})=\Lm_{r-1}-2\Lm_r,\q {\rm wt}(v_0)=0, \]
where $\Lm_0=0$. We define the $\frg$-action on $V(\Lm_1)$ as follows:
\begin{eqnarray}
&& h v_j=\lan h,{\rm wt}(v_j)\ran v_j\ \ (h\in P^*,\ j\in J_{{\rm B}}), \\
&&f_iv_i=v_{i+1},\ f_iv_{\ovl{i+1}}=v_{\ovl i},\q
e_iv_{i+1}=v_i,\ e_iv_{\ovl i}=v_{\ovl{i+1}}
\q(1\leq i<r),\ \ \label{B-f1}\\
&&f_r v_r=v_{0},\qq 
e_r v_{\ovl r}=v_0,\qq f_r v_0=2v_{\ovl r},\qq e_r v_0=2v_r,\label{B-f2}
\end{eqnarray}
and the other actions are trivial.

Let $\Lm_i$ $(1\leq i\leq r-1)$ be the $i$-th fundamental weight of type ${\rm B}_r$.
Similar to the ${\rm C}_r$ case, the fundamental representation 
$V(\Lm_i)$
is embedded in $\wedge^i V(\Lm_1)$
with multiplicity free. In 
$\wedge^i V(\Lm_1)$, the highest (resp. lowest) weight 
vector $u_{\Lm_i}$ (resp. $v_{\Lm_i}$)
of $V(\Lm_i)$ is realized as the same form as in (\ref{C-h-l}).

The fundamental representation $V(\Lm_r)$ is called the {\it spin representation}. It can be realized as follows: Set 
\[ {\mathbf B}_{{\rm sp}}^{(r)}:=
\{(\ep_1,\cdots,\ep_r)|\ \ep_i\in\{+,- \}\q (i=1,2,\cd,r)\}, \quad V_{{\rm sp}}^{(r)}:=\bigoplus_{v\in{\mathbf B}_{{\rm sp}}^{(r)}}\bbC v,\] 
and define the $\ge$-action on $V_{{\rm sp}}^{(r)}$ as follows:

\begin{equation}\label{Bsp-f0}
h_i(\ep_1,\cdots,\ep_r)=
\begin{cases}
\frac{\ep_i\cdot 1-\ep_{i+1}\cdot 1}{2}(\ep_1,\cdots,\ep_r) & {\rm if}\ i<r, \\
\ep_r(\ep_1,\cdots,\ep_r) & {\rm if}\ i=r,
\end{cases}
\end{equation}

\begin{equation}\label{Bsp-f1}
f_i(\ep_1,\cdots,\ep_r)=
\begin{cases}
(\ep_1,\cdots,\overset{i}{-},\overset{i+1}{+},\cdots,\ep_r) & {\rm if}\ \ep_i=+,\ \ep_{i+1}=-,\ i\neq r, \\
(\ep_1,\cdots,\ep_{r-1},\overset{r}{-}) & {\rm if}\ \ep_r=+,\ i=r, \\
0 & {\rm otherwise,}
\end{cases}
\end{equation}

\begin{equation}\label{Bsp-f2}
e_i(\ep_1,\cdots,\ep_r)=
\begin{cases}
(\ep_1,\cdots,\overset{i}{+},\overset{i+1}{-},\cdots,\ep_r) & {\rm if}\ \ep_i=-,\ \ep_{i+1}=+,\ i\neq r, \\
(\ep_1,\cdots,\ep_{r-1},\overset{r}{+}) & {\rm if}\ \ep_r=-,\ i=r, \\
0 & {\rm otherwise.}
\end{cases}
\end{equation}

Then the module $V_{{\rm sp}}^{(r)}$ is isomorphic to $V(\Lm_r)$ as a $\ge$-module.

\subsection{Type ${\rm D}_r$}\label{SectFundD}

Let $\frg=\frs\fro(2r,\mathbb{C})$ be the simple Lie algebra of type ${\rm D}_r$. The Cartan matrix $A=(a_{i,j})_{i,j\in I}$ of $\frg$ is given by
\[a_{i,j}=
\begin{cases}
2 & {\rm if}\ i=j, \\
-1 & {\rm if}\ |i-j|=1\ {\rm and}\ (i,j)\neq (r,r-1),\ (r-1,r),\ {\rm or}\ (i,j)=(r-2,r),\ (r,r-2), \\
0 & {\rm otherwise.}  
\end{cases}
\]

Define the partial order on the set $J_{{\rm D}}:=\{i,\ovl i|1\leq i\leq r\}$ by 
\begin{equation}\label{D-order}
 1< 2<\cd< r-1<\ ^{r}_{\ovl{r}}\ < \ovl{r-1}< \cd< \ovl
 2< \ovl 1.
\end{equation}
Note that there is no order between $r$ and $\ovl{r}$. For $\frg=\lan \frh,e_i,f_i(i\in I)\ran$, 
let us describe the vector representation 
$V(\Lm_1)$. Set ${\mathbf B}^{(r)}:=
\{v_i,v_{\ovl i}|i=1,2,\cd,r\}$ and define 
$V(\Lm_1):=\bigoplus_{v\in{\mathbf B}^{(r)}}\bbC v$. The weights of $v_i$, $v_{\ovl{i}}$ $(i\in I)$ are as follows:
\begin{equation}\label{D-wtv}
 {\rm wt}(v_i)=\Lm_i-\Lm_{i-1},\q {\rm wt}(v_{\ovl{i}})=\Lm_{i-1}-\Lm_{i}\q (1\leq i\leq r-2,\ i=r),
\end{equation}
\[  {\rm wt}(v_{r-1})=\Lm_r+\Lm_{r-1}-\Lm_{r-2},\q {\rm wt}(v_{\ovl{r-1}})=\Lm_{r-2}-\Lm_{r-1}-\Lm_r, \]
where $\Lm_0=0$. We define the $\frg$-action on $V(\Lm_1)$ as follows:
\begin{eqnarray}
&& h v_j=\lan h,{\rm wt}(v_j)\ran v_j\ \ (h\in P^*,\ j\in J_{{\rm D}}), \\
&&f_iv_i=v_{i+1},\ f_iv_{\ovl{i+1}}=v_{\ovl i},\q
e_iv_{i+1}=v_i,\ e_iv_{\ovl i}=v_{\ovl{i+1}}
\q(1\leq i<r),\ \ \label{D-f1}\\
&&f_r v_r=v_{\ovl{r-1}},\ f_r v_{r-1}=v_{\ovl r},\qq 
e_r v_{\ovl r}=v_{r-1},\ e_r v_{\ovl{r-1}}=v_{r},\label{D-f2}
\end{eqnarray}
and the other actions are trivial.

Let $\Lm_i$ $(1\leq i\leq r-2)$ be the $i$-th fundamental weight of type ${\rm D}_r$. Similar to the ${\rm B}_r$ and ${\rm C}_r$ cases, the fundamental representation 
$V(\Lm_i)$
is embedded in $\wedge^i V(\Lm_1)$
with multiplicity free. In 
$\wedge^i V(\Lm_1)$, the highest (resp. lowest) weight 
vector $u_{\Lm_i}$ (resp. $v_{\Lm_i}$)
of $V(\Lm_i)$ is realized as the formula (\ref{C-h-l}).

The fundamental representations $V(\Lm_{r-1})$ and $V(\Lm_r)$ are also called the {\it spin representations}. They can be realized as follows: Set 
\[ {\mathbf B}_{{\rm sp}}^{(+,r)}\ ({\rm resp.}\ {\mathbf B}_{{\rm sp}}^{(-,r)}):=
\{(\ep_1,\cdots,\ep_r)|\ \ep_i\in\{+,- \},\q \ep_1\cdots\ep_r=+\ ({\rm resp.}\ -)\}, \]
\[ V_{{\rm sp}}^{(+,r)}\ ({\rm resp.}\ V_{{\rm sp}}^{(-,r)}):=\bigoplus_{v\in{\mathbf B}_{{\rm sp}}^{(+,r)}({\rm resp.}\ {\mathbf B}_{{\rm sp}}^{(-,r)})}\bbC v,\] 
and define the $\ge$-action on $V_{{\rm sp}}^{(\pm,r)}$ as follows:

\begin{equation}\label{Dsp-f0}
h_i(\ep_1,\cdots,\ep_r)=
\begin{cases}
\frac{\ep_i\cdot 1-\ep_{i+1}\cdot 1}{2}(\ep_1,\cdots,\ep_r) & {\rm if}\ i<r, \\
\frac{\ep_{r-1}\cdot 1+\ep_{r}\cdot 1}{2}(\ep_1,\cdots,\ep_r) & {\rm if}\ i=r,
\end{cases}
\end{equation}

\begin{equation}\label{Dsp-f1}
f_i(\ep_1,\cdots,\ep_r)=
\begin{cases}
(\ep_1,\cdots,\overset{i}{-},\overset{i+1}{+},\cdots,\ep_r) & {\rm if}\ \ep_i=+,\ \ep_{i+1}=-,\ i\neq r, \\
(\ep_1,\cdots,\overset{r-1}{-},\overset{r}{-}) & {\rm if}\ \ep_{r-1}=+,\ \ep_r=+,\ i=r, \\
0 & {\rm otherwise,}
\end{cases}
\end{equation}

\begin{equation}\label{Dsp-f2}
e_i(\ep_1,\cdots,\ep_r)=
\begin{cases}

(\ep_1,\cdots,\overset{i}{+},\overset{i+1}{-},\cdots,\ep_r) & {\rm if}\ \ep_i=-,\ \ep_{i+1}=+,\ i\neq r, \\
(\ep_1,\cdots,\overset{r-1}{+},\overset{r}{+}) & {\rm if}\ \ep_{r-1}=-,\ \ep_r=-,\ i=r, \\
0 & {\rm otherwise.}
\end{cases}
\end{equation}

Then the module $V_{{\rm sp}}^{(+,r)}$ (resp. $V_{{\rm sp}}^{(-,r)}$) is isomorphic to $V(\Lm_r)$ (resp. $V(\Lm_{r-1})$) as a $\ge$-module.

\section{Factorization theorem}\label{DBCs}

In this section, we shall introduce double Bruhat cells $G^{u,v}$ and their properties\cite{B-Z, F-Z}. For $l\in \mathbb{Z}_{>0}$, we set $[1,l]:=\{1,2,\cdots,l\}$.

\subsection{Double Bruhat cells}\label{factpro}

Let $G$ be a classical algebraic group over $\mathbb{C}$, $B$ and $B_-$ be two opposite Borel subgroups in $G$, $N\subset B$ and $N_-\subset B_-$ be their unipotent radicals, 
$H:=B\cap B_-$ a maximal torus. We set $\frg:={\rm Lie}(G)$ with the triangular decomposition $\frg=\frn_-\oplus \frh \oplus \frn$. Let $e_i$, $f_i$ $(i\in[1,r])$ be the generators of $\frn$, $\frn_-$. For $i\in[1,r]$ and $t \in \mathbb{C}$, we set
\begin{equation}\label{xiyidef} 
x_i(t):={\rm exp}(te_i),\ \ \ y_{i}(t):={\rm exp}(tf_i).
\end{equation}
Let $W:=\lan s_i |i=1,\cdots,r \ran$ be the Weyl group of $\frg$, where
$\{s_i\}$ are the simple reflections. We identify the Weyl group $W$ with ${\rm Norm}_G(H)/H$. An element 
\begin{equation}\label{smpl}
\ovl{s_i}:=x_i(-1)y_i(1)x_i(-1)
\end{equation}
is in ${\rm Norm}_G(H)$, which is a representative of $s_i\in W={\rm Norm}_G(H)/H$ \cite{N1}. For $u\in W$, let $u=s_{i_1}\cdots s_{i_n}$ be its reduced expression. Then we write $\ovl{u}=\ovl{s_{i_1}}\cdots \ovl{s_{i_n}}$, call $l(u):=n$ the length of $u$. We have two kinds of Bruhat decompositions of $G$ as follows:
\[ G=\displaystyle\coprod_{u \in W}B\ovl{u}B=\displaystyle\coprod_{u \in W}B_-\ovl{u}B_- .\]
Then, for $u$, $v\in W$, 
we define the {\it double Bruhat cell} $G^{u,v}$ as follows:
\[ G^{u,v}:=B\ovl{u}B \cap B_-\ovl{v}B_-. \]

\begin{defn}\label{redworddef}
Let $v=s_{j_n}\cdots s_{j_1}$ be a reduced expression of $v\in W$ $(j_n,\cdots,j_1\in [1,r])$. Then the finite sequence $\textbf{i}:=(j_n,\cdots,j_1)$ is called a {\it reduced word} for $v$.
\end{defn}

For example, the sequence $(3,2,1,3,2,1)$ is a reduced word of the element $s_3s_2s_1s_3s_2s_1$ of the Weyl group of type ${\rm B}_3$ or ${\rm C}_3$. In this paper, we mainly treat double Bruhat cells of the form $G^{e,v}:=B \cap B_-\ovl{v} B_-$.

\subsection{Factorization theorem}\label{fuctorisec}

In this subsection, we shall introduce isomorphisms between double Bruhat cells $G^{e,v}$ and $H\times (\mathbb{C}^{\times})^{l(v)}$. For $i \in [1,r]$ and $t\in \mathbb{C}^{\times}$, we set $\alpha_i^{\vee}(t):=t^{h_i}$. 

For a sequence $\textbf{i}=(i_1, \cdots ,i_n)$ 
($i_1,\cdots,i_n\in[1,r]$), 
we define a map $x^G_{\textbf{i}}:H\times \mathbb{C}^n \rightarrow G$ as 
\begin{equation}\label{xgdef}
x^G_{\textbf{i}}(a; t_1, \cdots, t_n):=a\cdot x_{i_1}(t_1)\cdots x_{i_n}(t_n).
\end{equation}

\begin{thm}\label{fp}${\cite{B-Z,F-Z}}$ For $v\in W$ and its reduced word ${\rm \bf{i}}$, the map $x^G_{{\rm \bf{i}}}$ is a biregular isomorphism from $H\times (\mathbb{C}^{\times})^{l(v)}$ to a Zariski open subset of $G^{e,v}$. 
\end{thm}

For $\textbf{i}=(i_1, \cdots ,i_n)$
($i_1,\cdots,i_n\in[1,r]$), we define a map
$\ovl{x}^G_{\textbf{i}}:H\times(\mathbb{C}^{\times})^{n}\rightarrow
G^{e,v}$ as
\[ \ovl{x}^G_{\textbf{i}}(a;t_1,\cdots,t_n)
=a x_{i_1}(t_1)\alpha_{i_1}^{\vee}(t_1)x_{i_2}(t_2)\alpha_{i_2}^{\vee}(t_2)\cdots x_{i_n}(t_n)\alpha_{i_n}^{\vee}(t_n) , \]
where $a\in H$ and $(t_1,\cdots,t_n)\in (\mathbb{C}^{\times})^{n}$. 

Now, let $G$ be a classical algebraic group of type ${\rm B}_r$, ${\rm C}_r$ or ${\rm D}_r$, and $c\in W$ be a Coxeter element such that a reduced word ${\rm \bf{i}}$ of $c^2$ can be written as
\begin{equation}\label{redwords2}
{\rm \bf{i}}=
(r,r-1,\cdots,2,1)^2.
\end{equation}

\begin{rem}\label{importantrem}
In the rest of the paper, we use double indexed variables $Y_{s,j}$ $(s\in \mathbb{Z}$, $j\in [1,r])$. If we see the variables $Y_{s,0}$, $Y_{s,j}$ 
$(r+1\leq j)$ then 
we understand $Y_{s,0}=Y_{s,j}=1$. For example, if $l=1$ then $Y_{s,l-1}=1$.
\end{rem}

\begin{prop}\label{gprime}
In the above setting, the map $\ovl{x}^G_{{\rm \bf{i}}}$ is a biregular isomorphism between $H\times(\mathbb{C}^{\times})^{2r}$ and a Zariski open subset of $G^{e,c^2}$.
\end{prop}
\nd
[{\sl Proof.}]

Let us prove this proposition in the case $G$ is type ${\rm C}_r$. In the case $G$ is type ${\rm B}_r$ or ${\rm D}_r$, we can prove it in the same way. 

In this proof, we use the notation
\[ \textbf{Y}:=(Y_{1,r},Y_{1,r-1}\cdots,Y_{1,1},Y_{2,r},Y_{2,r-1},\cdots,Y_{2,1}), \]
for variables instead of $(t_1,\cdots,t_{2r})\in (\mathbb{C}^{\times})^{2r}$.

We define a map
$\phi:H\times(\mathbb{C}^{\times})^{2r}\rightarrow
H\times(\mathbb{C}^{\times})^{2r}$,
\begin{equation*}
\phi(a;\textbf{Y})=(\Phi_H(a;\textbf{Y});\Phi_{1,r}(\textbf{Y}),\cdots,\Phi_{1,1}(\textbf{Y}),
\Phi_{2,r}(\textbf{Y}),\cdots,\Phi_{2,2}(\textbf{Y}),\Phi_{2,1}(\textbf{Y})),
\end{equation*} 
as
\begin{equation}\label{mbasea} 
 \Phi_H(a;\textbf{Y}):=a\cdot \prod^{r}_{i=1}\prod^{2}_{j=1}\al_i^{\vee}(Y_{j,i}), 
\end{equation}
and for $l\in\{1,2,\cdots,r \}$, 

\begin{equation}\label{mbase0} 
\Phi_{1,l}(\textbf{Y}):=
\begin{cases}
\frac{Y_{1,l-1}Y_{2,l-1}Y_{2,l+1}}{Y_{1,l}Y_{2,l}^{2}} & {\rm if}\ l<r, \\
\frac{Y^2_{1,r-1}Y^2_{2,r-1}}{Y_{1,r}Y_{2,r}^{2}} & {\rm if}\ l=r,
\end{cases}
\end{equation}
\begin{equation}\label{mbase01} 
\Phi_{2,l}(\textbf{Y}):=
\begin{cases}
\frac{Y_{2,l-1}}{Y_{2,l}} & {\rm if}\ l<r, \\
\frac{Y^2_{2,r-1}}{Y_{2,r}} & {\rm if}\ l=r.
\end{cases}
\end{equation}

Note that $\phi$ is a biregular isomorphism since we can construct the inverse map $\psi:H\times(\mathbb{C}^{\times})^{2r}\rightarrow
H\times(\mathbb{C}^{\times})^{2r}$, 
\[ \psi(a;\textbf{Y})=(\Psi_H(a;\textbf{Y});\Psi_{1,r}(\textbf{Y}),\cdots,\Psi_{1,1}(\textbf{Y}),\Psi_{2,r}(\textbf{Y}),\cdots,\Psi_{2,1}(\textbf{Y}))
\]
of $\phi$ as follows:
\[
\Psi_{2,l}(\textbf{Y}):=
\begin{cases}
(Y_{2,1}\cdots Y_{2,l-1}Y_{2,l})^{-1} & {\rm if}\ l<r, \\
(Y_{2,1}\cdots Y_{2,r-2}Y_{2,r-1})^{-2}Y^{-1}_{2,r}  & {\rm if}\ l=r,
\end{cases}
\]
and $\Psi_{1,l}(\textbf{Y})$ is defined inductively as
\[
\Psi_{1,l}(\textbf{Y}):=
\begin{cases}
\frac{\Psi_{1,l-1}(\textbf{Y})\Psi_{2,l-1}(\textbf{Y})\Psi_{2,l+1}(\textbf{Y})}{\Psi^2_{2,l}(\textbf{Y})Y_{1,l}} & {\rm if}\ l<r, \\
\frac{\Psi^2_{1,r-1}(\textbf{Y})\Psi^2_{2,r-1}(\textbf{Y})}{\Psi^2_{2,r}(\textbf{Y})Y_{1,r}}  & {\rm if}\ l=r,
\end{cases}
\]
and
\[ \Psi_{H}(a;\textbf{Y}):=
a\cdot(\prod^{r}_{i=1}\prod^{2}_{j=1}\al_i^{\vee}(\Psi_{j,i}(\textbf{Y})))^{-1}.
\]
Then, the map $\psi$ is the inverse map of $\phi$.

Let us prove
\[ \ovl{x}^G_{\textbf{i}}(a;\textbf{Y})=(x^G_{\textbf{i}}\circ\phi)(a;\textbf{Y}), \]
which implies that $\ovl{x}^G_{\textbf{i}}:H\times(\mathbb{C}^{\times})^{2r}\rightarrow G^{e,c^2}$ is a biregular isomorphism by Theorem \ref{fp}. First, it is known that for $1\leq i,\ j\leq r$ and $s,\ t\in \mathbb{C}^{\times}$,
\begin{equation}\label{base2}
\al_i^{\vee}(s)x_{j}(t)=
x_{j}(s^{a_{i,j}}t)\al_i^{\vee}(s) ,
\end{equation}
where $(a_{i,j})_{i,j\in I}$ is the Cartan matrix of type ${\rm C}_r$. On the other hand, it follows from the definition (\ref{xgdef}) of $x^G_{\textbf{i}}$ and $(\ref{mbasea})$ that
\begin{multline*}(x^G_{\textbf{i}}\circ\phi)(a;\textbf{Y})
=a\cdot \left(\prod^{r}_{i=1}\prod^{2}_{s=1}\al_i^{\vee}(Y_{s,i})\right)
 \times x_{r}(\Phi_{1,{r}}(\textbf{Y}))
\cdots x_{1}(\Phi_{1,{1}}(\textbf{Y}))\\
\times x_{r}(\Phi_{2,r}(\textbf{Y}))
\cdots x_{2}(\Phi_{2,2}(\textbf{Y}))x_{1}(\Phi_{2,1}(\textbf{Y})).
\end{multline*}

For each $l$ $(1\leq l< r)$, we can move 
\[
\prod^{l}_{i=1}\al_{i}^{\vee}(Y_{1,i})
\prod^{r}_{i=1}\al_{i}^{\vee}(Y_{2,i})
\]
to the right of
$x_{l}(\Phi_{1,l}(\textbf{Y}))$ by using the relations (\ref{base2}):
\begin{eqnarray*}
& &\left(\prod^{l}_{i=1}\al_{i}^{\vee}(Y_{1,i})
\prod^{r}_{i=1}\al_{i}^{\vee}(Y_{2,i})\right)
x_{l}(\Phi_{1,l}(\textbf{Y}))\\
&=&x_{l}(\Phi_{1,l}(\textbf{Y})\frac{Y_{1,l}^2Y_{2,l}^2}{Y_{1,l-1}Y_{2,l-1}Y_{2,l+1}})
\prod^{l}_{i=1}\al_{i}^{\vee}(Y_{1,i})\prod^{r}_{i=1}\al_{i}^{\vee}(Y_{2,i})\\
&=&x_{l}(Y_{1,l})
\prod^{l}_{i=1}\al_{i}^{\vee}(Y_{1,i})
\prod^{r}_{i=1}\al_{i}^{\vee}(Y_{2,i}).
\end{eqnarray*}
In the same way, we see that 
\[
\left(\prod^{r}_{i=1}\al_{i}^{\vee}(Y_{1,i})
\prod^{r}_{i=1}\al_{i}^{\vee}(Y_{2,i})\right)
x_{r}(\Phi_{1,r}(\textbf{Y}))=x_{r}(Y_{1,r})\left(\prod^{r}_{i=1}\al_{i}^{\vee}(Y_{1,i})
\prod^{r}_{i=1}\al_{i}^{\vee}(Y_{2,i})\right).
\]
Similarly, for $1\leq l\leq r$, we can also move $\prod^{l}_{i=1}\al_{i}^{\vee}(Y_{2,i})$ to the right of $x_{l}(\Phi_{2,l}(\textbf{Y}))$:
\[
\left(\prod^{l}_{i=1}\al_{i}^{\vee}(Y_{2,i})\right)x_{l}(\Phi_{2,l}(\textbf{Y}))
=x_{l}(Y_{2,l})\left(\prod^{l}_{i=1}\al_{i}^{\vee}(Y_{2,i})\right).
\]

Thus, we get
\begin{multline*} 
(x^G_{\textbf{i}}\circ\phi)(a;\textbf{Y})=
a\cdot x_{r}(Y_{1,r})\al^{\vee}_{r}(Y_{1,r})\cdots x_{1}(Y_{1,1})\al^{\vee}_{1}(Y_{1,1})\\
 x_{r}(Y_{2,r})\al^{\vee}_{r}(Y_{2,r})\cdots
x_{2}(Y_{2,2})\al^{\vee}_{2}(Y_{2,2}) x_{1}(Y_{2,1})\al^{\vee}_{1}(Y_{2,1})=\ovl{x}^G_{\textbf{i}}(a;\textbf{Y}).
\end{multline*}
\qed

\section{Cluster algebras and generalized minors}\label{CluSect}
Following \cite{A-F-Z,F-Z,FZ2,M-M-A}, we review the definitions of cluster algebras and their generators called cluster variables. It is known that the coordinate rings of double Bruhat cells have cluster algebra structures, and generalized minors are their initial cluster variables \cite{GY,GY2}. We will refer to a relation between cluster variables on double Bruhat cells and crystal bases in Sect.\ref{gmc}.

We set $[1,l]:=\{1,2,\cdots,l\}$ and $[-1,-l]:=\{-1,-2,\cdots,-l\}$ for $l\in \mathbb{Z}_{>0}$. For $n,m\in \mathbb{Z}_{>0}$, let $x_1, \cdots ,x_n,x_{n+1}, \cdots
,x_{n+m}$ be commuting variables and $\cF:=\mathbb{C}(x_{1}, \cdots ,x_{n},x_{n+1},\cdots,x_{n+m})$ 
be the field of rational functions.

\subsection{Cluster algebras of geometric type}

In this subsection, we recall the definitions of cluster algebras. Let $\tilde{B}=(b_{ij})_{1\leq i\leq n+m,\ 1\leq j \leq n}$ be an $(n+m)\times
n$ integer matrix. The {\it principal part} $B$ of $\tilde{B}$ is obtained from $\tilde{B}$ by deleting the last $m$ rows. For $\tilde{B}$ and $k\in [1,n]$, the new $(n+m)\times n$ integer matrix $\mu_k(\tilde{B})=(b'_{ij})$ is defined by
\[b_{ij}':=
\begin{cases}
	-b_{ij} & {\rm if}\ i=k\ {\rm or}\ j=k, \\
	b_{ij}+\frac{|b_{ik}|b_{kj}+b_{ik}|b_{kj}|}{2} & {\rm otherwise}.
\end{cases}
\]
One calls $\mu_k(\tilde{B})$ the {\it matrix mutation} in direction $k$ of $\tilde{B}$. If there exists a positive 
integer diagonal matrix $D$ such that $DB$ is skew symmetric, we say $B$ is {\it skew symmetrizable}. Then we also say $\tilde{B}$ is skew symmetrizable. It is easily verified that if $\tilde{B}$ is skew symmetrizable then $\mu_k(\tilde{B})$ is also skew symmetrizable${\cite[Proposition3.6]{M-M-A}}$. We can also verify that $\mu_k\mu_k(\tilde{B})=\tilde{B}$. Define $\textbf{x}:=(x_1,\cdots,x_{n+m})$ and we call the pair $(\textbf{x}, \tilde{B})$ {\it initial seed} and $x_1,\cdots,x_{n+m}$ {\it initial cluster variables}. For $k\in[1,n]$, a new cluster variable $x_k'$ is defined by the following {\it exchange relation}.
\begin{equation}\label{exrel} x_k x_k' = 
\prod_{1\leq i \leq n+m,\ b_{ik}>0} x_i^{b_{ik}}
+\prod_{1\leq i \leq n+m,\ b_{ik}<0} x_i^{-b_{ik}}. \end{equation}
Let $\mu_k(\textbf{x})$ be the set of variables obtained from $\textbf{x}$ by replacing $x_k$ by $x'_k$. Ones call the pair $(\mu_k(\textbf{x}), \mu_k(\tilde{B}))$ the {\it mutation} in direction $k$ of the seed $(\textbf{x}, \tilde{B})$ and denote it by $\mu_k((\textbf{x}, \tilde{B}))$.

Now, we can repeat this process of mutation and obtain a set of seeds inductively. Hence, each seed consists of an $(n+m)$-tuple of variables and a matrix. Ones call this $(n+m)$-tuple and matrix {\it cluster} and {\it exchange matrix} respectively. Variables in cluster is called {\it cluster variables}. In particular, the variables $x_{n+1},\cdots,x_{n+m}$ are called {\it frozen cluster variables}.

\begin{defn}${\cite{F-Z, M-M-A}}$\label{clusterdef}
Let $\tilde{B}$ be an integer matrix whose principal part is skew symmetrizable and $\Sigma=(\textbf{x},\tilde{B})$ a seed. We set ${\mathbb A}:={\mathbb Z}[x_{n+1}^{\pm1}, \cdots ,x_{n+m}^{\pm1}]$. The cluster algebra (of geometric type) $\cA=\cA(\Sigma)$ over $\mathbb A$ associated with
seed $\Sigma$ is defined as the ${\mathbb A}$-subalgebra of $\cF$ generated by all cluster variables in all seeds which can be obtained from $\Sigma$ by sequences of mutations.
\end{defn}

\subsection{Cluster algebra $\cA({\rm \bf{i}})$}\label{cAi}

Let $G$ be a classical algebraic group, $\ge:={\rm Lie}(G)$ and $A=(a_{i,j})$ be its Cartan matrix. In Definition \ref{redworddef}, we defined a reduced word ${\rm \bf{i}}=(j_{n},\cdots,j_2,j_{1})$ for an element $v$ of Weyl group $W$. In this subsection, we define a cluster algebra $\cA({\rm \bf{i}})$, which is obtained from ${\rm \bf{i}}$. It satisfies that $\cA({\rm \bf{i}})\otimes \mathbb{C}$ is isomorphic to the coordinate ring $\mathbb{C}[G^{e,v}]$ of the double Bruhat cell \cite{A-F-Z}. Let $j_k$ $(k\in[1,n])$ be the $k$-th index of ${\rm \bf{i}}$ from the right. Let us add $r$ additional entries $j_{-r},\cdots,j_{-1}$ at the beginning of ${\rm \bf{i}}$ by setting $j_{-t}=-t$ $(t\in [1,r])$.

For $l\in[1,n]$, we denote by $l^-$ the largest index $k\in[1,n]$ such that $k<l$ and 
$j_l=j_k$. For $l\in[-1,-r]$, let $l^-$ be the largest index $k\in[1,n]$ such that 
$| j_l| =| j_k| $. If such $k$ does not exist, we set $l^-:=0$. 
For example, if $[-1,-3]\cup{\mathbf{i}}=(-3,-2,-1,3,2,1,3,2,1)$ then, $(-1)^-=4$,
 $(-2)^-=5$, $(-3)^-=6$, $4^-=1$, $5^-=2$, $6^-=3$, and $3^-=2^-=1^-=0$. 
We define a set $e({{\bf i}})$ as 
\[ e({\mathbf{i}}):= \{k\in[1,n]\ |\ k^->0\}. 
\]
Following \cite{A-F-Z}, we define a directed graph $\Gamma_{{\rm \bf{i}}}$ as follows. The vertices of $\Gamma_{{\rm \bf{i}}}$ are the variables $x_k$ ($k\in[-1,-r]\cup[1,n]$). For two vertices $x_k$ $(k\in [-1,-r]\cup[1,n])$ and $x_l$ $(l\in e({\rm \bf{i}}))$ with either $l<k$ or $k\in[-1,-r]$, there exists an arrow $x_k\rightarrow x_l$ (resp. $x_l\rightarrow x_k$) if and only if $l=k^-$ (resp. $l^-<k^-<l$ and $a_{|j_k|,|j_l|}<0$). For two
vertices $x_k$ $(k\in[1,n]\setminus e({\mathbf{i}}))$ and
 $x_l$ $(l\in e({\mathbf{i}}))$ with
$k<l$, there exists an arrow 
$x_l\rightarrow x_k$ (resp.
$x_k\rightarrow x_l$) if and only if $k=l^-$ 
(resp. $k^-<l^-<k$ and
$a_{|j_l|,|j_k|}<0$). Next, let us define a matrix $\tilde{B}=\tilde{B}({\rm \bf{i}})$. 

\begin{defn}
Let $\tilde{B}({\rm \bf{i}})$ be an integer matrix with rows labelled by all the indices in $[-1,-r]\cup [1,n]$ and columns labelled by all the indices in $e({\rm \bf{i}})$. For $k\in[-1,-r]\cup [1,n]$ and $l\in e({\rm \bf{i}})$, an entry $b_{kl}$ of $\tilde{B}({\rm \bf{i}})$ is determined as follows: If there exists an arrow $x_k\rightarrow x_l$ (resp. $x_l\rightarrow x_k$) in $\Gamma_{{\rm \bf{i}}}$, then
\[
b_{kl}:=\begin{cases}
		1\ ({\rm resp.}\ -1)& {\rm if}\ |j_k|=|j_l|, \\
		-a_{|j_k||j_l|}\ ({\rm resp.}\ a_{|j_k||j_l|})& {\rm if}\ |j_k|\neq|j_l|.
	\end{cases}
\]
If there exist no arrows between $k$ and $l$, we set $b_{kl}=0$. The principal part $B({\rm \bf{i}})$ of $\tilde{B}({\rm \bf{i}})$ is the submatrix $(b_{i,j})_{i,j\in e({\rm \bf{i}})}$.
\end{defn}

\begin{prop}\label{propss}${\cite{A-F-Z}}$
$\tilde{B}({\rm \bf{i}})$ is skew symmetrizable. 
\end{prop}

\begin{defn}${\cite{A-F-Z}}$
We set $\textbf{x}=(x_i)_{i\in[-1,-r]\cup[1,n]}$ and define the cluster algebra $\cA({\rm \bf{i}})$ over $\mathbb{Z}[x^{\pm1}_i |\ i\in[-1,-r]\cup[1,n]\setminus e(\textbf{i})]$ as $\cA({\rm \bf{i}}):=\mathcal{A}((\textbf{x},\tilde{B}(\textbf{i})))$.
\end{defn}
In this definition, we use the notation $x_i$ ($i\in[-1,-r]\cup[1,n]\setminus e(\textbf{i}))$ for frozen cluster variables instead of $x_{n+1},\cdots,x_{n+m}$ in Definition \ref{clusterdef}.
 
\subsection{Generalized minors}\label{bilingen}

Set $\cA({\rm \bf{i}})_{\mathbb{C}}:=\cA({\rm \bf{i}})\otimes \mathbb{C}$. It is known that the coordinate ring $\mathbb{C}[G^{e,v}]$ of the double Bruhat cell is isomorphic to $\cA({\rm \bf{i}})_{\mathbb{C}}$ (Theorem \ref{clmainthm}). To describe this isomorphism explicitly, we need generalized minors.  

We set $G_0:=N_-HN$, and let $x=[x]_-[x]_0[x]_+$ with $[x]_-\in N_-$, $[x]_0\in H$, $[x]_+\in N$ be the corresponding decomposition. 

\begin{defn}
For $i\in[1,r]$ and $w,\ w'\in W$, the {\it generalized minor} $\Delta_{w'\Lambda_i,w\Lambda_i}$ is a regular function on $G$ whose restriction to the open set $\ovl{w'}G_0\ovl{w}^{-1}$ is given by $\Delta_{w'\Lambda_i,w\Lambda_i}(x)=([\ovl{w'}^{-1}x \ovl{w} ]_0)^{\Lambda_i}$. Here, $\Lambda_i$ is the $i$-th  fundamental weight and $\ovl w$ is the one we defined in (\ref{smpl}).
\end{defn}
The generalized minor $\Delta_{w'\Lambda_i,w\Lambda_i}$ depends on $w'\Lambda_i,$ $w\Lambda_i$ and does not depend on $w',\ w$. By definition, for $a\in H$, $x\in G$, $w\in W$, $i,j\in I$ and $t\in\mathbb{C}$,
\begin{equation}\label{genbasic}
\Delta_{\Lambda_i,w\Lambda_i}(ax)=a^{\Lambda_i}\Delta_{\Lambda_i,w\Lambda_i}(x),\ 
\Delta_{\Lambda_i,\Lambda_i}(x x_{j}(t))=\Delta_{\Lambda_i,\Lambda_i}(x),
\end{equation}
where $x_{j}(t)\in N$ is the one in (\ref{xiyidef}).

Let $\omega:\ge\to\ge$ be the anti-involution 
\[
\omega(e_i)=f_i,\q
\omega(f_i)=e_i,\q \omega(h)=h,
\] and extend it to $G$ by setting
$\omega(x_i(c))=y_{i}(c)$, $\omega(y_{i}(c))=x_i(c)$ and $\omega(t)=t$
$(t\in H)$. Here, $x_i(t)$ and $y_i(t)$ were defined in (\ref{xiyidef}). One can calculate the generalized minors as follows. There exists a $\ge$ (or $G$)-invariant bilinear form on the
finite-dimensional  irreducible
$\ge$-module $V(\lm)$ such that 
\[
 \lan au,v\ran=\lan u,\omega(a)v\ran,
\q\q(u,v\in V(\lm),\,\, a\in \ge\ (\text{or }G)).
\]
For $g\in G$, 
we have the following simple fact:
\[
 \Del_{\Lm_i,\Lm_i}(g)=\lan gu_{\Lm_i},u_{\Lm_i}\ran,
\]
where $u_{\Lm_i}$ is a properly normalized highest weight vector in
$V(\Lm_i)$. Hence, for $w,\ w'\in W$, we have
\begin{equation}\label{minor-bilin}
 \Del_{w'\Lm_i,w\Lm_i}(g)=
\Del_{\Lm_i}(\ovl{w'}^{-1}g\ovl w)=
\lan g\ovl w\cdot u_{\Lm_i},\ovl{w'}\cdot u_{\Lm_i}\ran.
\end{equation}
Note that $\omega(\ovl s_i^{\pm})=\ovl s_i^{\mp}$.

\subsection{Cluster algebras on double Bruhat cells}

For a reduced expression $v=s_{j_n}s_{j_{n-1}}\cdots s_{j_1}\in W$ and $k\in [1,n]$, we set
\begin{equation}\label{inc}
v_{> k}=v_{> k}({\rm \bf{i}}):=s_{j_{1}}s_{j_{2}}\cdots s_{j_{n-k}}.
\end{equation}
%For $k \in [-1,-r]$, we set $v_{> k}:=v^{-1}$ and $j_k:=k$. 
For $k \in [1,n]$, we define $\Delta(k;{\rm \bf{i}})(x):=\Delta_{\Lambda_{j_k},v_{>n-k+1}\Lambda_{j_k}}(x)$, and for $k \in [-1,-r]$, $\Delta(k;{\rm \bf{i}})(x):=\Delta_{\Lambda_{|k|},v^{-1}\Lambda_{|k|}}(x)$. Finally, we set $F({\rm \bf{i}}):=\{ \Delta(k;{\rm \bf{i}})(x)|k \in [-1,-r]\cup[1,n] \}$. It is known that the set $F({\rm \bf{i}})$ is an algebraically independent generating set for the field of rational functions $\mathbb{C}(G^{e,v})$ \cite[Theorem 1.12]{F-Z}. Then, we have the following.
\begin{thm}\label{clmainthm}${\cite{A-F-Z, GLS, GY, GY2}}$
The isomorphism of fields $\varphi :\cF \rightarrow \mathbb{C}(G^{e,v})$ defined by $\varphi (x_k)=\Delta(k;{\rm \bf{i}})\ (k \in [-1,-r]\cup [1,n] )$ restricts to an isomorphism of algebras $\cA({\rm \bf{i}})_{\mathbb{C}}\rightarrow \mathbb{C}[G^{e,v}]$.
\end{thm}

\begin{ex}\label{clmainthmex}
Let $G$ be a classical algebraic group of type ${\rm B}_r$, ${\rm C}_r$ or ${\rm D}_r$, $v=c^2$ be the square of a Coxeter element such that whose reduced word ${\rm \bf{i}}$ is written as in $(\ref{redwords2})$. Then for $k\in[1,r]$, we have $j_{r+k}=j_k=k$ and the isomorphism is given by
\begin{eqnarray*}
 x_{-k}&\mapsto& \Delta_{\Lambda_{k},c^{-2}\Lambda_{k}}=
\Delta_{\Lambda_{k},(s_{1}s_{2}\cdots s_{r})^2\Lambda_{k}}
=\Delta_{\Lambda_{k},(s_{1}s_{2}\cdots s_{r})(s_{1}s_{2}\cdots s_{k})\Lambda_{k}}
=\Delta_{\Lambda_{k},c^2_{>r-k}\Lambda_{k}}
,\\
 x_{r+k}&\mapsto & \Delta_{\Lambda_{k},c^2_{>r-k+1}\Lambda_{k}}=
\Delta_{\Lambda_{k},(s_{1}s_{2}\cdots s_{r})(s_{1}s_{2}\cdots s_{k-1})\Lambda_{k}}\\
& &=\Delta_{\Lambda_{k},(s_{1}s_{2}\cdots s_{k})\Lambda_{k}}
=\Delta_{\Lambda_{k},c^2_{>2r-k}\Lambda_{k}},\\
 x_k &\mapsto& \Delta_{\Lambda_{k},c^2_{>2r-k+1}\Lambda_{k}}=
\Delta_{\Lambda_{k},(s_{1}s_{2}\cdots s_{k-1})\Lambda_{k}}
=\Delta_{\Lambda_{k},\Lambda_{k}}.
\end{eqnarray*}
 
\end{ex}

\subsection{Finite type}

Let $\mathcal{S}$ be the set of seeds of a cluster algebra $\mathcal{A}$. If $\mathcal{S}$ is finite, then $\mathcal{A}$ is said to be {\it finite type}. In this subsection, we shall review the cluster algebras of finite type \cite{FZ3}.

Let $B=(b_{ij})$ be an integer square matrix. The {\it Cartan counter part} of $B$ is a generalized Cartan matrix $A=A(B)=(a_{i,j})$ defined as follows:
\[ a_{i,j}=
\begin{cases}
2 & {\rm if}\ i=j,\\
-|b_{i,j}| & {\rm if}\ i\neq j.
\end{cases}
\]

\begin{thm}\label{finthm1}$\cite{FZ3}$
For a cluster algebra $\cA$ with the set $\mathcal{S}$ of seeds, the following statements are equivalent:
\begin{enumerate}
\item The cluster algebra $\cA$ is of finite type.
\item There exists a seed $\Sigma=(\textbf{y},\tilde{B})$ such that $\cA= \cA(\Sigma)$ and $A(B)$ is a Cartan matrix of finite type, where $B$ is the principal part of $\tilde{B}$.
\item Let $(\textbf{y}',\tilde{B}')$ be an arbitrary seed in $\mathcal{S}$ and $(b_{i,j})$ be the principal part of $\tilde{B}'$. Then $|b_{i,j}b_{j,i}|\leq3$. %%2-finite
\end{enumerate}
\end{thm}

By this theorem, we can define the {\it type} of each cluster algebra of finite type mirroring the Cartan-Killing classification.

Let $\Phi$ be the root system associated with a Cartan matrix, with the set of simple roots $\Pi=\{\al_i|\ i\in I\}$ and the set of positive roots $\Phi_{>0}$. The set of {\it almost positive roots} is defined as $\Phi_{\geq -1}:=\Phi_{>0}\cup -\Pi$.

\begin{thm}\label{finthm2}$\cite{FZ3}$
\begin{enumerate}
\item For a cluster algebra $\cA$ of finite type, the number of the cluster variables included in $\cA$ is equal to $|\Phi_{\geq -1}|$. Here, $\Phi$ is the root system associated with the Cartan matrix of the same type as $\cA$.
\item Let $c\in W$ be a Coxeter element of a classical algebraic group $G$ whose length $l(c)$ satisfies 
$l(c^2)=2 l(c)=2${\rm rank}$(G)$. Then the coordinate ring $\mathbb{C}[G^{e,c^2}]$ has a structure of cluster algebra of finite type under the isomorphism in Theorem \ref{clmainthm}, and its type is the Cartan-Killing type of $G$.
\end{enumerate}
\end{thm}

Next, we define the following graph, which is a similar notion to the {\it weighted graph} introduced in \cite{FZ3}.

\begin{defn}\label{wei-dia}
Let $\Sigma=(\textbf{y},\tilde{B})$ be a seed with $\textbf{y}=(y_i)_{i\in[-1,-r]\cup[1,n]}$ and 
an $(r+n)\times |e({\rm \bf{i}})|$-skew symmetrizable matrix $\tilde{B}=(b_{i,j})$ which satisfies
$b_{i,j}\in\{-2,-1,0,1,2\}$ (the rows of $\tilde{B}$ are labelled by $[-1,-r]\cup[1,n]$ as above), where $y_i=x_i$ for $i\in[-1,-r]\cup[1,n]\setminus e(\textbf{i})$. We suppose that if $i,j\in e({\rm \bf{i}})$ then $|b_{i,j}b_{j,i}|\leq 3$.
We define $\Gamma(\Sigma)$ as the labelled directed graph whose vertices are $y_{-r},\cdots,y_{-1},y_1,\cdots,y_{n}$, and whose arrows and its labels are determined as follows: For $i,j\in e({\rm \bf{i}})$, there exists the arrow  $\overset{2}{y_i\rightarrow y_j}$ (resp. $\overset{-2}{y_j\rightarrow y_i}$) if and only if $b_{i,j}=2$ and $b_{j,i}=-1$
(resp. $b_{i,j}=-2$ and $b_{j,i}=1$). Further, there exists the arrow $y_i\rightarrow y_j$ if and only if $b_{i,j}=1$ and $b_{j,i}=-1$. For $i\in [-1,-r]\cup[1,n]\setminus e({\rm \bf{i}})$ and $j\in e({\rm \bf{i}})$, there exists the arrow  $\overset{2}{y_i\rightarrow y_j}$ (resp. $\overset{-2}{y_j\rightarrow y_i}$) if and only if $b_{i,j}=2$
(resp. $b_{i,j}=-2$). Further, there exists the arrow $y_i\rightarrow y_j$ (resp. $y_j\rightarrow y_i$) if and only if $b_{i,j}=1$ (resp. $b_{i,j}=-1$). We call the graph $\Gamma(\Sigma)$ {\it mutation diagram} of $\Sigma$. We understand the arrows $y_i\rightarrow y_j$ $(i,j\in[-1,-r]\cup[1,n])$ have the labels $1$ and do not denote it. 
\end{defn}

Note that the graph $\Gamma((\textbf{x},\tilde{B}({\rm \bf{i}})))$ is obtained from $\Gamma_{{\rm \bf{i}}}$ by labelling
arrows properly. %%% \Gamma_{{\rm \bf{i}}} の部分を、$\Gamma((\textbf{x},\tilde{B}({\rm \bf{i}})))$と書き直す。

\begin{lem}$\cite{FZ3,M-M-A}$\label{mutgamlem}
Let $\Sigma=(\textbf{y},\tilde{B})$ be a seed as in Definition \ref{wei-dia}. For $k\in e({\rm \bf{i}})$, the graph $\Gamma(\mu_k(\Sigma))$ has vertices $y_{-r},\cdots,y_{-1},y_1,\cdots,y_k',\cdots,y_{n}$, and edges or subgraphs of $\Gamma(\Sigma)$ are transformed to those of $\Gamma(\mu_k(\Sigma))$ by $\mu_k$ as follows:
\begin{enumerate}
\item[$(1)$]If $y_i\rightarrow y_k$ (resp. $y_k\rightarrow y_i$) in $\Gamma(\Sigma)$ then $y_k'\rightarrow y_i$ (resp. $y_i\rightarrow y_k'$) in $\Gamma(\mu_k(\Sigma))$. If $\overset{\pm2}{y_i\rightarrow y_k}$ in $\Gamma(\Sigma)$ then $\overset{\mp2}{y_k'\rightarrow y_i}$ in $\Gamma(\mu_k(\Sigma))$.
\item[$(2)$] For $i,j\in [-1,-r]\cup[1,n]$, we suppose that either $i$ or $j$ (or both) belong to $e({\rm \bf{i}})$. Then 
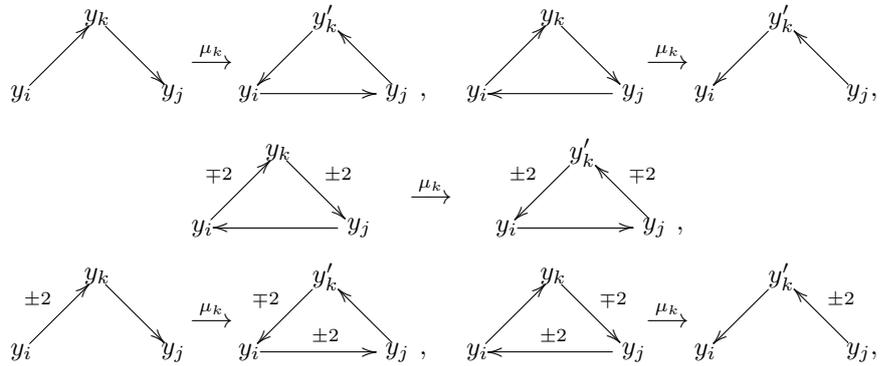
\[
\begin{xy}
(20,90)*{y_j}="j",
(10,100) *{y_k}="k",
(0,90)*{y_i}="i",
(25,95) *{\overset{\mu_k}{\longrightarrow}}="12",
(50,90)*{\ y_j\ ,}="jj",
(40,100) *{y'_k}="kk",
(30,90)*{y_i}="ii",
(80,90)*{\ y_j}="jjj",
(70,100) *{y_k}="kkk",
(60,90)*{y_i}="iii",
(85,95) *{\overset{\mu_k}{\longrightarrow}}="1212",
(110,90)*{\ y_j,}="jjjj",
(100,100) *{y'_k}="kkkk",
(90,90)*{y_i}="iiii",
\ar@{->} "i";"k"
\ar@{->} "k";"j"
\ar@{->} "jj";"kk"
\ar@{->} "kk";"ii"
\ar@{->} "ii";"jj"
\ar@{->} "iii";"kkk"
\ar@{->} "kkk";"jjj"
\ar@{->} "jjjj";"kkkk"
\ar@{->} "kkkk";"iiii"
\ar@{->} "jjj";"iii"
\end{xy} 
\]
\[
\begin{xy}
(20,90)*{\ y_j}="j",
(10,100) *{y_k}="k",
(0,90)*{y_i}="i",
(30,95) *{\overset{\mu_k}{\longrightarrow}}="12",
(60,90)*{\ y_j\ ,}="jj",
(50,100) *{y'_k}="kk",
(40,90)*{y_i}="ii",
\ar@{->} "j";"i"
\ar@{->} "i";"k"^{\mp2}
\ar@{->} "k";"j"^{\pm2}
\ar@{->} "jj";"kk"_{\mp2}
\ar@{->} "kk";"ii"_{\pm2}
\ar@{->} "ii";"jj"
\end{xy} 
\]
\[
\begin{xy}
(20,90)*{y_j}="j",
(10,100) *{y_k}="k",
(0,90)*{y_i}="i",
(25,95) *{\overset{\mu_k}{\longrightarrow}}="12",
(50,90)*{\ y_j\ ,}="jj",
(40,100) *{y'_k}="kk",
(30,90)*{y_i}="ii",
(80,90)*{\ y_j}="jjj",
(70,100) *{y_k}="kkk",
(60,90)*{y_i}="iii",
(85,95) *{\overset{\mu_k}{\longrightarrow}}="1212",
(110,90)*{\ y_j,}="jjjj",
(100,100) *{y'_k}="kkkk",
(90,90)*{y_i}="iiii",
\ar@{->} "i";"k"^{\pm2}
\ar@{->} "k";"j"
\ar@{->} "jj";"kk"
\ar@{->} "kk";"ii"_{\mp2}
\ar@{->} "ii";"jj"^{\pm2}
\ar@{->} "iii";"kkk"
\ar@{->} "kkk";"jjj"^{\mp2}
\ar@{->} "jjjj";"kkkk"_{\pm2}
\ar@{->} "kkkk";"iiii"
\ar@{->} "jjj";"iii"_{\pm2}
\end{xy} 
\]
\[
\begin{xy}
(20,90)*{y_j}="j",
(10,100) *{y_k}="k",
(0,90)*{y_i}="i",
(25,95) *{\overset{\mu_k}{\longrightarrow}}="12",
(50,90)*{\ y_j\ ,}="jj",
(40,100) *{y'_k}="kk",
(30,90)*{y_i}="ii",
(80,90)*{\ y_j}="jjj",
(70,100) *{y_k}="kkk",
(60,90)*{y_i}="iii",
(85,95) *{\overset{\mu_k}{\longrightarrow}}="1212",
(110,90)*{\ y_j,}="jjjj",
(100,100) *{y'_k}="kkkk",
(90,90)*{y_i}="iiii",
\ar@{->} "i";"k"^{\pm2}
\ar@{->} "k";"j"
\ar@{->} "j";"i"_{\mp2}
\ar@{->} "jj";"kk"
\ar@{->}
 "kk";"ii"_{\mp2}
\ar@{->} "iii";"kkk"
\ar@{->} "kkk";"jjj"^{\pm2}
\ar@{->} "jjjj";"kkkk"_{\mp2}
\ar@{->}
 "iiii";"jjjj"^{\pm2}
\ar@{->} "kkkk";"iiii"
\end{xy}
\]

where the arrow $y_i\overset{-2\quad \ }{\longrightarrow y_j}$ with $i\notin e({\rm \bf{i}})$, $j\in e({\rm \bf{i}})$ implies the arrow $y_i\longrightarrow y_j$, and the arrow $y_i\overset{2\quad \ }{\longrightarrow y_j}$ with $i\in e({\rm \bf{i}})$, $j\notin e({\rm \bf{i}})$ implies the arrow $y_i\longrightarrow y_j$.
\item [$(3)$] If two vertices $y_i$ and $y_j$ are not connected via a two-arrow oriented
path going through $y_k$ in $\Gamma(\Sigma)$, the arrows between $y_i$ and $y_j$ and their labels remain unchanged by $\mu_k$.
\end{enumerate}
\end{lem}

We will use this lemma in Sect \ref{mainthmpr}.

\begin{rem}
In the above lemma, we do not mention to several subgraphs. For example, the subgraph
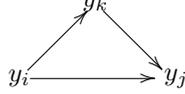
\[
\begin{xy}
(20,90)*{\ y_j}="j",
(10,100) *{y_k}="k",
(0,90)*{y_i}="i",
\ar@{->} "i";"j"
\ar@{->} "i";"k"
\ar@{->} "k";"j"
\end{xy}
\]
is not mentioned. But we will not treat these graphs in this article. 
\end{rem}

\begin{ex}\label{gammaex}
Let us consider the case $G={\rm SO}_{5}(\mathbb{C})$ $(type\ {\rm B}_2)$ and ${\rm \bf{i}}=(2,1,2,1)$. The graph $\Gamma((\textbf{x},\tilde{B}({\rm \bf{i}})))$ is described as
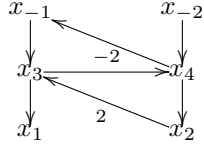
\[
\begin{xy}
(70,98) *{x_{-2}}="8",
(70,90)*{x_4}="4",
(70,82)*{x_2}="-2",
(50,98) *{x_{-1}}="6",
(50,90)*{x_3}="2",
(50,82)*{x_1}="-1",
\ar@{->} "8";"4"
\ar@{->} "6";"2"
\ar@{->} "4";"-2"
\ar@{->} "2";"-1"
\ar@{->} "2";"4"^{-2}
\ar@{->} "4";"6"
\ar@{->} "-2";"2"^2
\end{xy} \]
In general, let us consider the case $G={\rm SO}_{2r+1}(\mathbb{C})$ $(type\ {\rm B}_r)$ and ${\rm \bf{i}}$ is the sequence in $(\ref{redwords2})$. 
For $k$ $(1\leq k\leq r-2)$, vertices and arrows around the vertex $x_{r+k}$ in the graph $\Gamma((\textbf{x},\tilde{B}({\rm \bf{i}})))$ are
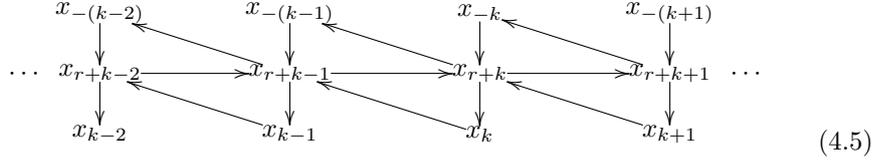
\begin{equation}\label{gammaex-1}
\begin{xy}
(95,90)*{\cdots}="emp1",
(85,98)*{x_{-(k+1)}}="7",
(85,90) *{x_{r+k+1}}="3",
(85,82)*{x_{k+1}}="-4",
(60,98)*{x_{-k}}="r+k",
(60,90) *{x_{r+k}}="k",
(60,82)*{x_k}="-j_k",
(35,98)*{x_{-(k-1)}}="8",
(35,90) *{x_{r+k-1}}="4",
(35,82)*{x_{k-1}}="-2",
(10,98)*{x_{-(k-2)}}="6",
(10,90) *{x_{r+k-2}}="2",
(10,82)*{x_{k-2}}="-1",
(0,90)*{\cdots}="emp",
\ar@{->} "7";"3"
\ar@{->} "r+k";"k"
\ar@{->} "8";"4"
\ar@{->} "6";"2"
\ar@{->} "3";"-4"
\ar@{->} "k";"-j_k"
\ar@{->} "k";"8"
\ar@{->} "4";"-2"
\ar@{->} "2";"-1"
\ar@{->} "k";"3"
\ar@{->} "4";"k"
\ar@{->} "-j_k";"4"
\ar@{->} "2";"4"
\ar@{->} "3";"r+k"
\ar@{->} "4";"6"
\ar@{->} "-4";"k"
\ar@{->} "-2";"2"
\end{xy} \end{equation}
Vertices and arrows around the vertex $x_{2r}$ in the graph $\Gamma((\textbf{x},\tilde{B}({\rm \bf{i}})))$ are
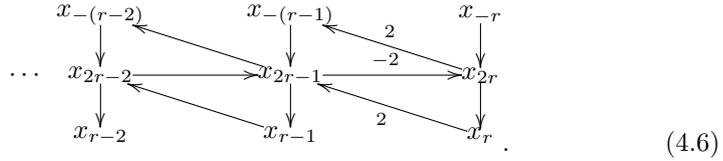
\begin{equation}\label{gammaex-2}
\begin{xy}
(60,98)*{x_{-r}}="r+k",
(60,90) *{x_{2r}}="k",
(60,82)*{x_r}="-j_k",
(35,98)*{x_{-(r-1)}}="8",
(35,90) *{x_{2r-1}}="4",
(35,82)*{x_{r-1}}="-2",
(10,98)*{x_{-(r-2)}}="6",
(10,90) *{x_{2r-2}}="2",
(10,82)*{x_{r-2}}="-1",
(0,90)*{\cdots}="emp",
\ar@{->} "r+k";"k"
\ar@{->} "8";"4"
\ar@{->} "6";"2"
\ar@{->} "k";"-j_k"
\ar@{->} "k";"8"_2
\ar@{->} "4";"-2"
\ar@{->} "2";"-1"
\ar@{->} "4";"k"^{-2}
\ar@{->} "-j_k";"4"^2
\ar@{->} "2";"4"
\ar@{->} "4";"6"
\ar@{->} "-2";"2"
\end{xy}. \end{equation}
\end{ex}

\section{Monomial realizations and Demazure crystals}

In Sect.\ref{gmc}, we shall describe the cluster variables in a cluster algebra of finite type in terms of the {\it monomial realizations} of Demazure crystals. Let us recall the notion of crystal base and its monomial realization in this section. Let $\ge$ be a complex simple Lie algebra and $I=\{1,2,\cdots,r\}$ the index set.

\subsection{Monomial realizations of crystals}\label{monoreacry}

In this subsection, we review the monomial realizations of crystals~\cite{K, K2, Nj}. First, let us recall the crystals.

\begin{defn}\label{defcry}$\cite{K0}$
A~{\it crystal} associated with a Cartan matrix~$A$ is a~set~$B$ together with the maps $\text{wt}: B \rightarrow P$,
$\tilde{e_{i}}$, $\tilde{f_{i}}: B \cup \{0\} \rightarrow B \cup \{0\}$ and $\varepsilon_i$,
$\varphi_i: B \rightarrow {\mathbb Z} \cup \{-\infty \}$, $i \in I$, satisfying some properties.
\end{defn}
We call $\{\tilde{e}_i,\ \tilde{f}_i\}$ the {\it Kashiwara operators}. Let $U_q(\mathfrak g)$ be the quantum enveloping algebra \cite{K0} associated with a Cartan matrix~$A$, that is,  $U_q(\mathfrak g)$ has generators $\{e_i,\ f_i,\ h_i |\ i\in I\}$ over $\mathbb{C}(q)$ satisfying some relations, where $q$ is an indeterminate. Let $V(\lm)$ ($\lm\in P^+=\oplus_{i\in I}\mathbb{Z}_{\geq0}\Lambda_i$) be the finite dimensional irreducible representation of $U_q(\mathfrak g)$ which has the highest weight vector $v_{\lm}$, and $B(\lm)$ be the crystal base of $V(\lm)$. The crystal base $B(\lm)$ has a crystal structure.

Let us introduce monomial realizations which realize each element of $B(\lambda)$ as a~certain Laurent monomial. First, fix a cyclic sequence of the indices $\cdots(i_1,i_2,\cdots,i_r)(i_1,i_2,\cdots,i_r)\cdots$ such that $\{i_1,i_2,\cdots,i_r\}=I$. And we can associate this sequence with a~set of integers $p=(p_{j,i})_{j,i \in I,\; j \neq i}$ such that
\begin{gather*}
p_{i_a,i_b}=
\begin{cases}
1 & \text{if} \quad  a<b,
\\
0 & \text{if} \quad  a>b.
\end{cases}
\end{gather*}

Second, for doubly-indexed variables $\{Y_{s,i} \,|\, i \in I$, $s\in \mathbb{Z}\}$, we def\/ine the set of monomials
\begin{gather*}
{\mathcal Y}:=\left\{Y=\prod\limits_{s \in \mathbb{Z},\ i \in I}
Y_{s,i}^{\zeta_{s,i}}\, \Bigg| \,\zeta_{s,i} \in \mathbb{Z},\
\zeta_{s,i} =0~\text{except for f\/initely many}~(s,i) \right\}.
\end{gather*}
Finally, we def\/ine maps $\text{wt}: {\mathcal Y} \rightarrow P$, $\varepsilon_i$, $\varphi_i: {\mathcal Y} \rightarrow
\mathbb{Z}$, $i \in I$. %, as follows.
For $Y=\prod\limits_{s \in \mathbb{Z},\; i \in I} Y_{s,i}^{\zeta_{s,i}}\in {\mathcal Y}$,
\begin{gather}\label{wtph}
\text{wt}(Y):= \sum\limits_{i,s}\zeta_{s,i}\Lambda_i,\!
\quad
\varphi_i(Y):=\max\left\{\! \sum\limits_{k\leq s}\zeta_{k,i}  \,|\, s\in \mathbb{Z} \!\right\},\!
\quad
\varepsilon_i(Y):=\varphi_i(Y)-\text{wt}(Y)(h_i).
\end{gather}

We set
\begin{gather}
\label{asidef}
A_{s,i}:=Y_{s,i}Y_{s+1,i}\prod\limits_{j\neq i}Y_{s+p_{j,i},j}^{a_{j,i}}
\end{gather}
and def\/ine the Kashiwara operators as follows
\begin{gather*}
\tilde{f}_iY=
\begin{cases}
A_{n_{f_i},i}^{-1}Y & \text{if} \quad  \varphi_i(Y)>0,
\\
0 & \text{if} \quad  \varphi_i(Y)=0,
\end{cases}
\qquad
\tilde{e}_iY=
\begin{cases}
A_{n_{e_i},i}Y & \text{if} \quad  \varepsilon_i(Y)>0,
\\
0 & \text{if} \quad  \varepsilon_i(Y)=0,
\end{cases}
\end{gather*}
where
\begin{gather*}
n_{f_i}:=\min \left\{n \,\Bigg|\, \varphi_i(Y)= \sum\limits_{k\leq n}\zeta_{k,i}\right\},
\qquad
n_{e_i}:=\max \left\{n \,\Bigg|\, \varphi_i(Y)= \sum\limits_{k\leq n}\zeta_{k,i}\right\}.
\end{gather*}

Then the following theorem holds:

\begin{thm}$\cite{K, Nj}$\label{monorealmain}\quad
\begin{enumerate}\itemsep=0pt
\item[(i)] For the set $p=(p_{j,i})$ as above, $({\mathcal Y}, {\rm wt}, \varphi_i, \varepsilon_i,\tilde{f}_i,
\tilde{e}_i)_{i\in I}$ is a~crystal.
When we emphasize~$p$, we write ${\mathcal Y}$ as ${\mathcal Y}(p)$.
\item[(ii)] If a~monomial $Y \in {\mathcal Y}(p)$ satisfies $\varepsilon_i(Y)=0$ for all $i \in I$,
then the connected component containing~$Y$ is isomorphic to $B({\rm wt}(Y))$.
\end{enumerate}
\end{thm}

\subsection{Demazure crystals}\label{Demcrysub}

For $w\in W$ and $\lambda\in P^+$, a {\it Demazure crystal}
$B(\lambda)_w\subset B(\lambda)$ is inductively defined as follows. 

\begin{defn}
Let $u_{\lambda}$ be the highest weight vector of $B(\lambda)$. For
 the identity element $e$ of $W$, we set
 $B(\lambda)_e:=\{u_{\lambda}\}$. 
For $w\in W$, if $s_iw<w$, 
\[ B(\lambda)_w:=\{\tilde{f}_i^kb\ |\ k\geq0,\ b\in
 B(\lambda)_{s_iw},\ \tilde{e}_ib=0\}\setminus \{0\}. 
\]
\end{defn}

\begin{thm}\label{kashidem}${\cite{K3}}$
For $w\in W$, let $w=s_{i_1}\cdots s_{i_n}$ be an arbitrary reduced
 expression. Let $u_{\lambda}$ be the highest
weight vector of $B(\lambda)$. Then 
\[B(\lambda)_w=\{\tilde{f}_{i_1}^{a(1)}\cdots
\tilde{f}_{i_n}^{a(n)}u_{\lambda}|a(1),\cdots,a(n)\in\mathbb{Z}_{\geq0}
\}
\setminus\{0\}.
\]
\end{thm}

\begin{ex}\label{monorealex0}
Let us consider the case of type ${\rm C}_{2}$ and cyclic sequence is $(2,1)$. In the notation of $(\ref{asidef})$, we can write
\[
 A_{s,i}=
\begin{cases}
\frac{Y_{s,1}Y_{s+1,1}}{Y_{s+1,2}} & {\rm if}\ i=1, \\
\frac{Y_{s,2}Y_{s+1,2}}{Y^2_{s,1}} & {\rm if}\ i=2.
\end{cases}
\]
In general, if each factor of a monomial $Y\in\mathcal{Y}$ has non-negative degree, then $\varepsilon_i(Y)=0$ for all $i\in I=\{1,2\}$.
Therefore, we have $\varepsilon_i(Y_{1,1})=0$. Hence, we can consider the monomial realization of crystal base $B(\Lambda_1)$ such that the highest weight vector in $B(\Lambda_1)$ is realized by $Y_{1,1}:$
\begin{equation}\label{mono0-1}
Y_{1,1}\overset{\tilde{f}_1}{\longrightarrow} \frac{Y_{2,2}}{Y_{2,1}}
\overset{\tilde{f}_2}{\longrightarrow} \frac{Y_{2,1}}{Y_{3,2}}
\overset{\tilde{f}_1}{\longrightarrow} \frac{1}{Y_{3,1}}.
\end{equation}
Similarly, we get the monomial realization of crystal base $B(\Lambda_2)$ such that the highest weight vector in $B(\Lambda_2)$ is realized by $Y_{1,2}:$
\begin{equation}\label{mono0-2}
Y_{1,2}\overset{\tilde{f}_2}{\longrightarrow} \frac{Y^2_{1,1}}{Y_{2,2}}
\overset{\tilde{f}_1}{\longrightarrow} \frac{Y_{1,1}}{Y_{2,1}}
\overset{\tilde{f}_1}{\longrightarrow} \frac{Y_{2,2}}{Y^2_{2,1}}
\overset{\tilde{f}_2}{\longrightarrow} \frac{1}{Y_{3,2}}.
\end{equation}

\end{ex}

\begin{ex}\label{monorealex1}
Let us consider the case of type ${\rm B}_{3}$ and cyclic sequence is $(3,2,1)$.
In the notation of $(\ref{asidef})$, we can write
\[
 A_{s,i}=
\begin{cases}
\frac{Y_{s,1}Y_{s+1,1}}{Y_{s+1,2}} & {\rm if}\ i=1, \\
\frac{Y_{s,2}Y_{s+1,2}}{Y_{s,1}Y^2_{s+1,3}} & {\rm if}\ i=2, \\
\frac{Y_{s,3}Y_{s+1,3}}{Y_{s,2}} & {\rm if}\ i=3.
\end{cases}
\]
We can consider the monomial realization of crystal base $B(\Lambda_1)$ such that the highest weight vector in $B(\Lambda_1)$ is realized by $Y_{1,1}:$
\[
Y_{1,1}\overset{\tilde{f}_1}{\longrightarrow} \frac{Y_{2,2}}{Y_{2,1}}
\overset{\tilde{f}_2}{\longrightarrow} \frac{Y^2_{3,3}}{Y_{3,2}}
\overset{\tilde{f}_3}{\longrightarrow} \frac{Y_{3,3}}{Y_{4,3}}
\overset{\tilde{f}_3}{\longrightarrow} \frac{Y_{3,2}}{Y^2_{4,3}}
\overset{\tilde{f}_2}{\longrightarrow} \frac{Y_{3,1}}{Y_{4,2}}
\overset{\tilde{f}_1}{\longrightarrow} \frac{1}{Y_{4,1}}.
\]

\end{ex}

\section{Cluster variables and crystals}\label{gmc}

Let $G$ be a classical algebraic group of type ${\rm B}_r$, ${\rm C}_r$ or ${\rm D}_r$. In this section, we describe the cluster variables on a double Bruhat cell by the total sum of monomial realizations of Demazure crystals. In the rest of the paper, we only treat the Coxeter element $c\in W$ such that a reduced word ${\rm \bf{i}}$ of $c^2$ can be written as (\ref{redwords2}). Let $j_k$ be the $k$-th index of ${\rm \bf{i}}$ from the right, which implies $j_k=j_{r+k}=k$ $(1\leq k\leq r)$. We shall consider the monomial realization associated with the sequence $(r,\cdots,2,1)$ (Sect.\ref{monoreacry}). 

Let $\textbf{V}:=((\varphi_{\textbf{V}})_{2r},\cdots,(\varphi_{\textbf{V}})_{r+1},(\varphi_{\textbf{V}})_{r},\cdots,(\varphi_{\textbf{V}})_{1},(\varphi_{\textbf{V}})_{-r},\cdots,(\varphi_{\textbf{V}})_{-1})$, where $(\varphi_{\textbf{V}})_k\in \mathbb{C}[G^{e,c^2}]$ are defined as follows:
\begin{equation}\label{inicluvar}
(\varphi_{\textbf{V}})_k=
\begin{cases}
\Delta_{\Lambda_{j_k},c^{2}_{>2r-k}\Lambda_{j_k}} & {\rm if}\ 1\leq k\leq 2r, \\
\Delta_{\Lambda_{|k|},\Lambda_{|k|}} & {\rm if}\ -r\leq k\leq -1.
\end{cases}
\end{equation}
By Theorem \ref{clmainthm}, Example \ref{clmainthmex} and Theorem \ref{finthm2}, we can regard $\mathbb{C}[G^{e,c^2}]$ as a cluster algebra of finite type and $\textbf{V}$ as its initial cluster. Moreover, $(\varphi_{\textbf{V}})_{2r},\cdots,(\varphi_{\textbf{V}})_{r+1}$ and $(\varphi_{\textbf{V}})_{-r},\cdots,(\varphi_{\textbf{V}})_{-1}$ are frozen cluster variables.

Using these notation, we can rewrite the graph (\ref{gammaex-1}) of type ${\rm B}_r$ as
\begin{equation}\label{inidiaB-1}
\begin{xy}
(95,90)*{\cdots}="emp1",
(85,98)*{(\varphi_{\textbf{V}})_{r+k+1}}="7",
(85,90) *{(\varphi_{\textbf{V}})_{k+1}}="3",
(85,82)*{(\varphi_{\textbf{V}})_{-(k+1)}}="-4",
(60,98)*{(\varphi_{\textbf{V}})_{r+k}}="r+k",
(60,90) *{(\varphi_{\textbf{V}})_{k}}="k",
(60,82)*{(\varphi_{\textbf{V}})_{-k}}="-j_k",
(35,98)*{(\varphi_{\textbf{V}})_{r+k-1}}="8",
(35,90) *{(\varphi_{\textbf{V}})_{k-1}}="4",
(35,82)*{(\varphi_{\textbf{V}})_{-(k-1)}}="-2",
(10,98)*{(\varphi_{\textbf{V}})_{r+k-2}}="6",
(10,90) *{(\varphi_{\textbf{V}})_{k-2}}="2",
(10,82)*{(\varphi_{\textbf{V}})_{-(k-2)}}="-1",
(0,90)*{\cdots}="emp",
\ar@{->} "7";"3"
\ar@{->} "r+k";"k"
\ar@{->} "8";"4"
\ar@{->} "6";"2"
\ar@{->} "3";"-4"
\ar@{->} "k";"-j_k"
\ar@{->} "k";"8"
\ar@{->} "4";"-2"
\ar@{->} "2";"-1"
\ar@{->} "k";"3"
\ar@{->} "4";"k"
\ar@{->} "-j_k";"4"
\ar@{->} "2";"4"
\ar@{->} "3";"r+k"
\ar@{->} "4";"6"
\ar@{->} "-4";"k"
\ar@{->} "-2";"2"
\end{xy}
\end{equation}
Similarly, we can also rewrite (\ref{gammaex-2}) by using these notation:
\begin{equation}\label{inidiaB-2}
\begin{xy}
(60,98)*{(\varphi_{\textbf{V}})_{2r}}="r+k",
(60,90) *{(\varphi_{\textbf{V}})_{r}}="k",
(60,82)*{(\varphi_{\textbf{V}})_{-r}}="-j_k",
(35,98)*{(\varphi_{\textbf{V}})_{2r-1}}="8",
(35,90) *{(\varphi_{\textbf{V}})_{r-1}}="4",
(35,82)*{(\varphi_{\textbf{V}})_{-(r-1)}}="-2",
(10,98)*{(\varphi_{\textbf{V}})_{2r-2}}="6",
(10,90) *{(\varphi_{\textbf{V}})_{r-2}}="2",
(10,82)*{(\varphi_{\textbf{V}})_{-(r-2)}}="-1",
(0,90)*{\cdots}="emp",
\ar@{->} "r+k";"k"
\ar@{->} "8";"4"
\ar@{->} "6";"2"
\ar@{->} "k";"-j_k"
\ar@{->} "k";"8"_2
\ar@{->} "4";"-2"
\ar@{->} "2";"-1"
\ar@{->} "4";"k"^{-2}
\ar@{->} "-j_k";"4"^2
\ar@{->} "2";"4"
\ar@{->} "4";"6"
\ar@{->} "-2";"2"
\end{xy}\end{equation}
 In the rest of the paper, when we write a cluster in $\mathbb{C}[G^{e,c^2}]$, we drop frozen variables. For example, $\textbf{V}=((\varphi_{\textbf{V}})_{r},\cdots,(\varphi_{\textbf{V}})_{1})$.
We will order the cluster variables $(\varphi_{\textbf{V}})_{1},\cdots,(\varphi_{\textbf{V}})_{r}$ from the right in $\textbf{V}$ as above, and let $\mu_k$ denote the mutation of the $k$-th cluster variable from the right. For a cluster $\textbf{T}$ in $\mathbb{C}[G^{e,c^2}]$, let $(\varphi_{\textbf{T}})_{k}$ denote the $k$-th (non-frozen) cluster variable from the right:
\[ \textbf{T}:=((\varphi_{\textbf{T}})_{r},\cdots,(\varphi_{\textbf{T}})_{1}).\] 
Each cluster variable is a regular function on $G^{e,c^2}$. On the other hand, by Proposition \ref{gprime}% (resp. Theorem \ref{fp2})
, it can be seen as a function on
$H\times
(\mathbb{C}^{\times})^{2r}$ %(resp. $(\mathbb{C}^{\times})^{2r}$)
. Then, let us consider the following change of variables:
\begin{defn}\label{gendef}
Along (\ref{redwords2}), we set the variables $\textbf{Y}\in (\mathbb{C}^{\times})^{2r}$ as
\begin{equation}\label{yset}
\textbf{Y}:=
(Y_{1,r},Y_{1,r-1},\cdots,Y_{1,2},Y_{1,1},Y_{2,r},Y_{2,r-1},\cdots,Y_{2,2},Y_{2,1}).
\end{equation}
Then for $a\in H$ and cluster $\textbf{T}$ in $\mathbb{C}[G^{e,c^2}]$, we define
\[ 
(\varphi^G_{\textbf{T}})_k(a;\textbf{Y}):=(\varphi_{\textbf{T}})_k\circ
 \ovl{x}^G_{{\rm \bf{i}}}(a;\textbf{Y}),\quad (1\leq k\leq r),
\] 
where $\ovl{x}^G_{{\rm \bf{i}}}$ is as in Proposition \ref{gprime}.
\end{defn}

\begin{ex}\label{initialex4} 
Let us consider the case $G={\rm Sp}_{4}(\mathbb{C})$ $(type\ {\rm C}_2)$ and ${\rm \bf{i}}=(2,1,2,1)$.
In the above setting, for $k\in[1,r]$,
\[ (\varphi^G_{\textbf{V}})_k(a;\textbf{Y})=(\varphi_{\textbf{V}})_k\circ \ovl{x}^G_{{\rm \bf{i}}}(a;\textbf{Y}). \]
The definition $(\ref{inicluvar})$ says $(\varphi_{\textbf{V}})_1=\Delta_{\Lambda_1,c^2_{>3}\Lambda_1}=\Delta_{\Lambda_1,s_1\Lambda_1}$. Using the bilinear form $(\ref{minor-bilin})$, it follows from the actions $(\ref{B-wtv})$, $(\ref{B-f1})$, $(\ref{B-f2})$ on the fundamental representation that 
\begin{eqnarray*}
& &\hspace{-20pt}(\varphi^G_{\textbf{V}})_1(a;\textbf{Y})=(\varphi_{\textbf{V}})_1\circ \ovl{x}^G_{{\rm \bf{i}}}(a;\textbf{Y})\\
&=& \lan a\cdot x_2(Y_{1,2})\al^{\vee}_2(Y_{1,2})x_1(Y_{1,1})\al^{\vee}_1(Y_{1,1})
x_2(Y_{2,2})\al^{\vee}_2(Y_{2,2})x_1(Y_{2,1})\al^{\vee}_1(Y_{2,1}) \ovl{s_1}v_1, v_1\ran\\
&=&a^{\Lm_1}\lan v_2, \al^{\vee}_1(Y_{2,1})y_1(Y_{2,1})\al^{\vee}_2(Y_{2,2})y_2(Y_{2,2}) \al^{\vee}_1(Y_{1,1})y_1(Y_{1,1})\al^{\vee}_2(Y_{1,2})y_2(Y_{1,2})v_1 \ran\\
&=&a^{\Lm_1}\lan v_2,  \al^{\vee}_1(Y_{2,1})y_1(Y_{2,1})\al^{\vee}_2(Y_{2,2})y_2(Y_{2,2}) (Y_{1,1}v_1+v_2) \ran\\
&=&a^{\Lm_1}\lan v_2, \left(Y_{1,1}Y_{2,1}v_1+(Y_{1,1}+\frac{Y_{2,2}}{Y_{2,1}})v_2
+Y_{2,1}v_{\ovl{2}}+v_{\ovl{1}}\right) \ran\\
&=&a^{\Lm_1}(Y_{1,1}+\frac{Y_{2,2}}{Y_{2,1}}).
\end{eqnarray*}
Comparing with (\ref{mono0-1}), the set of monomials $\{Y_{1,1},\frac{Y_{2,2}}{Y_{2,1}}\}$ coincides with the monomial realization of the Demazure crystal $B(\Lm_1)_{s_1}$ in Example \ref{monorealex0} (see Theorem \ref{kashidem}).

Note that the value of $(\varphi^G_{\textbf{V}})_1(a;\textbf{Y})$ was reduced to the calculation of the coefficient of $v_2$ in $\al^{\vee}_1(Y_{2,1})y_1(Y_{2,1})\al^{\vee}_2(Y_{2,2})y_2(Y_{2,2}) \al^{\vee}_1(Y_{1,1})y_1(Y_{1,1})\al^{\vee}_2(Y_{1,2})y_2(Y_{1,2})v_1$. Similarly, the value of $(\varphi^G_{\textbf{V}})_2(a;\textbf{Y})$ was reduced to the calculation of the coefficient of $v_2\wedge v_{\ovl{1}}$ in  $\al^{\vee}_1(Y_{2,1})y_1(Y_{2,1})\cdots\al^{\vee}_2(Y_{1,2})y_2(Y_{1,2})v_1\wedge v_2$. Hence, to calculate $(\varphi^G_{\textbf{V}})_2(a;\textbf{Y})$, we need
\begin{eqnarray*}
& &\hspace{-20pt}\al^{\vee}_1(Y_{2,1})y_1(Y_{2,1})\al^{\vee}_2(Y_{2,2})y_2(Y_{2,2}) \al^{\vee}_1(Y_{1,1})y_1(Y_{1,1})\al^{\vee}_2(Y_{1,2})y_2(Y_{1,2})v_2\\
&=&\al^{\vee}_1(Y_{2,1})y_1(Y_{2,1})\al^{\vee}_2(Y_{2,2})y_2(Y_{2,2})
(\frac{Y_{1,2}}{Y_{1,1}}v_2+Y_{1,1}v_{\ovl{2}}+v_{\ovl{1}})\\
&=&\frac{Y_{1,2}Y_{2,2}}{Y_{1,1}Y_{2,1}}v_2+(\frac{Y_{1,1}Y_{2,1}}{Y_{2,2}}+\frac{Y_{1,2}Y_{2,1}}{Y_{1,1}})v_{\ovl{2}}+(\frac{1}{Y_{2,1}}+\frac{Y_{1,1}}{Y_{2,2}}+\frac{Y_{1,2}}{Y_{1,1}})v_{\ovl{1}}.
\end{eqnarray*}

Just as in the case of $(\varphi^G_{\textbf{V}})_1(a;\textbf{Y})$, one get
\begin{eqnarray*}
& &\hspace{-20pt}(\varphi^G_{\textbf{V}})_2(a;\textbf{Y})=(\varphi_{\textbf{V}})_2\circ x^G_{{\rm \bf{i}}}\circ \phi(a;\textbf{Y})\\
&=&a^{\Lm_2}\lan v_2\wedge v_{\ovl{1}}, \al^{\vee}_1(Y_{2,1})y_1(Y_{2,1})\cdots\al^{\vee}_2(Y_{1,2})y_2(Y_{1,2})v_1\wedge v_2 \ran\\
&=&
a^{\Lm_2}\lan v_2\wedge v_{\ovl{1}}, \left( (Y_{1,1}+\frac{Y_{2,2}}{Y_{2,1}})v_2\wedge (\frac{1}{Y_{2,1}}+\frac{Y_{1,1}}{Y_{2,2}}+\frac{Y_{1,2}}{Y_{1,1}})v_{\ovl{1}}\right)
+v_{\ovl{1}} \wedge \frac{Y_{1,2}Y_{2,2}}{Y_{1,1}Y_{2,1}}v_2\ran
\\
&=&a^{\Lm_2}(Y_{1,2}+\frac{Y^2_{1,1}}{Y_{2,2}}+2\frac{Y_{1,1}}{Y_{2,1}}+\frac{Y_{2,2}}{Y^2_{2,1}}).
\end{eqnarray*}
Comparing with (\ref{mono0-2}), the set of monomials $\{Y_{1,2},\frac{Y^2_{1,1}}{Y_{2,2}},\frac{Y_{1,1}}{Y_{2,1}},\frac{Y_{2,2}}{Y^2_{2,1}}\}$ is equal to the monomial realization of the Demazure crystal $B(\Lm_2)_{s_1s_2}$ in Example \ref{monorealex1}.

\end{ex}

The following theorems are our main results, which mean relations between the cluster variables in $\mathbb{C}[G^{e,c^2}]$ and Demazure crystals. 
\begin{thm}\label{thmB}
Let $G={\rm SO}_{2r+1}(\mathbb{C})$ be the classical algebraic group of type ${\rm B}_r$. The cluster variables in $\mathbb{C}[G^{e,c^2}]$ are the total sums of monomial realizations of certain Demazure crystals. More precisely, each cluster variable is described as follows:
\begin{enumerate}
\item For $k\in[1,r]$, we obtain
\begin{eqnarray}
(\varphi^G_{\textbf{V}})_k(a;\textbf{Y})&=&
a^{\Lambda_k} Y_{1,k}(1+A^{-1}_{1,k}+A^{-1}_{1,k}A^{-1}_{1,k-1}+\cdots+A^{-1}_{1,k}A^{-1}_{1,k-1}\cdots A^{-1}_{1,1})\quad \quad \quad \label{iniBexp2}\\
&=&a^{\Lambda_k} \sum_{b\in B(\Lambda_k)_{c^2_{>2r-k}}} \mu(b)\nonumber,
\end{eqnarray}
where $\mu:B(\Lambda_k)\rightarrow \mathcal{Y}$ is the monomial realization of $B(\Lambda_k)$ such that the highest weight vector in $B(\Lambda_k)$ is realized by $Y_{1,k}\in\mathcal{Y}$. We also obtain
\[
(\varphi^G_{\mu_{k}\cdots\mu_{r-1}\mu_r(\textbf{V})})_k(a;\textbf{Y})=
\begin{cases}
a^{\Lambda_{k-1}+\Lambda_{r-1}} Y_{2,k}=a^{\Lambda_{k-1}+\Lambda_{r-1}}\sum_{b\in B(\Lambda_k)_{e}} \mu'(b) & {\rm if}\ k<r,\\
a^{\Lambda_{r-1}} Y_{2,r}=a^{\Lambda_{r-1}} \sum_{b\in B(\Lambda_k)_{e}} \mu'(b) & {\rm if}\ k=r,
\end{cases}
\]
where $\mu'$ is the monomial realization of $B(\Lambda_k)$ such that the highest weight vector in $B(\Lambda_k)$ is realized by $Y_{2,k}\in\mathcal{Y}$.
\item For $k$ and $l$ with $1\leq k\leq l\leq r-2$,
\begin{eqnarray}
& &\hspace{-20pt}(\varphi^G_{\mu_{l}\cdots\mu_{k+1}\mu_k(\textbf{V})})_l(a;\textbf{Y}) \nonumber\\ 
&=&
a^{\Lm_{k-1}+\Lm_{l+1}} Y_{2,k}Y_{1,l+1} (1+A^{-1}_{1,l+1}+A^{-1}_{1,l+1}A^{-1}_{1,l}+\cdots+A^{-1}_{1,l+1}A^{-1}_{1,l}\cdots A^{-1}_{1,k+1})
\qquad \quad \label{B-exp1}\\
&=&a^{\Lm_{k-1}+\Lm_{l+1}}\sum_{b\in B(\Lambda_k+\Lambda_{l+1})_{s_{k+1}s_{k+2}\cdots s_{l+1}}} \mu(b), \label{B-cry1}
\end{eqnarray}
where $\mu$ is the monomial realization of $B(\Lambda_k+\Lambda_{l+1})$ such that the highest weight vector in $B(\Lambda_k+\Lambda_{l+1})$ is realized by $Y_{2,k}Y_{1,l+1}\in\mathcal{Y}$.
\item For $k\in[1,r-1]$, we obtain
\begin{eqnarray}
& &\hspace{-20pt}(\varphi^G_{\mu_{r-1}\cdots\mu_{k+1}\mu_k(\textbf{V})})_{r-1}(a;\textbf{Y})\nonumber\\
&=&a^{2\Lm_{r}+\Lm_{k-1}} \left(Y^2_{1,r} Y_{2,k}+ 2Y^2_{1,r}Y_{2,k}A^{-1}_{1,r}(1+A^{-1}_{1,r-1}+A^{-1}_{1,r-1}A^{-1}_{1,r-2}+\cdots+A^{-1}_{1,r-1}\cdots A^{-1}_{1,k+1})\right.\nonumber\\
& &\quad +Y^2_{1,r}Y_{2,k}A^{-2}_{1,r}
(1+A^{-1}_{1,r-1}+\cdots+A^{-1}_{1,r-1}\cdots A^{-1}_{1,k+1})^2\nonumber\\
&+&\left.Y_{1,k-1}Y_{1,r-1}
(1+A^{-1}_{1,r-1}+\cdots+A^{-1}_{1,r-1}\cdots A^{-1}_{1,k+1}) (1+A^{-1}_{1,k-1}+\cdots+A^{-1}_{1,k-1}\cdots A^{-1}_{1,1})\right)\nonumber\\
&=&a^{\Lm_{k-1}+2\Lm_{r}} \sum_{b\in B(\Lambda_k+2\Lambda_{r})_{s_{k+1}s_{k+2}\cdots s_{r}}} C(b)\mu(b)\nonumber\\
& &+a^{\Lm_{k-1}+2\Lm_{r}}\sum_{b\in B(\Lambda_{k-1}+ \Lambda_{r-1})_{s_1\cdots s_{k-1}s_{k+1}s_{k+2}\cdots s_{r-1}}} \mu'(b), \label{B-cry2}
\end{eqnarray}
where $\mu$ is the monomial realization of $B(\Lambda_k+2\Lambda_{r})$ such that the highest weight vector is realized by $Y_{2,k}Y^2_{1,r}\in\mathcal{Y}$, $\mu'$ is the monomial realization of $B(\Lambda_{k-1}+\Lambda_{r-1})$ such that the highest weight vector is realized by $Y_{1,k-1}Y_{1,r-1}\in\mathcal{Y}$, and $C(b)$ are some positive integers. We also obtain
\begin{eqnarray}
& &\hspace{-20pt}(\varphi^G_{\mu_r\mu_{r-1}\cdots\mu_{k+1}\mu_k(\textbf{V})})_r(a;\textbf{Y})\nonumber\\
&=&a^{\Lm_{k-1}+\Lm_r}Y_{2,k}Y_{1,r}(1+A^{-1}_{1,r}+A^{-1}_{1,r}A^{-1}_{1,r-1}+\cdots+A^{-1}_{1,r}A^{-1}_{1,r-1}\cdots A^{-1}_{1,k+1}),\qquad \qquad \label{B-lem2pr-7}\\
&=&
a^{\Lm_{k-1}+\Lm_{r}} \sum_{b\in B(\Lambda_k+\Lambda_{r})_{s_{k+1}s_{k+2}\cdots s_{r}}} \mu''(b),\label{B-cry3}
\end{eqnarray}
where $\mu''$ is the monomial realization of $B(\Lambda_k+\Lambda_{r})$ such that the highest weight vector is realized by $Y_{2,k}Y_{1,r}\in\mathcal{Y}$.
\item For $j$ and $k$ with $1\leq j< k\leq r-1$,
\begin{eqnarray}
& &\hspace{-20pt}(\varphi^G_{\mu_j\cdots\mu_{k-2}\mu_{k-1}\mu_{r}\mu_{r-1}\cdots\mu_{k+1}\mu_k(\textbf{V})})_j(a;\textbf{Y})
= a^{\Lm_{j-1}+\Lm_{k-1}+2\Lm_{r}}\left(Y_{2,k}Y_{2,j}Y^2_{1,r}\right.\nonumber\\
&+&\left.2\frac{Y_{2,k}Y_{2,j}Y_{1,r-1}Y_{1,r}}{Y_{2,r}}(1+A^{-1}_{1,r-1}+A^{-1}_{1,r-1}A^{-1}_{1,r-2}+\cdots+A^{-1}_{1,r-1}\cdots A^{-1}_{1,k+1})\right)
\quad\quad \quad \label{B-lem3-exp}\\
&+&a^{\Lm_{j-1}+\Lm_{k-1}+2\Lm_{r}}\frac{Y_{2,k}Y_{2,j}Y^2_{1,r-1}}{Y^2_{2,r}}
(1+A^{-1}_{1,r-1}+A^{-1}_{1,r-1}A^{-1}_{1,r-2}+\cdots+A^{-1}_{1,r-1}\cdots A^{-1}_{1,k+1})\nonumber\\
& &\qquad\qquad\qquad\qquad\qquad\qquad\qquad  \times(1+A^{-1}_{1,r-1}+\cdots+A^{-1}_{1,r-1}\cdots A^{-1}_{1,j+1})\nonumber
\\
&=& a^{\Lm_{j-1}+\Lm_{k-1}+2\Lm_{r}} \sum_{b\in B(\Lm_j+\Lambda_k+2\Lambda_{r})_{s_{k+1}s_{k+2}\cdots s_{r}} } C(b)\mu(b)\label{B-lem3-cry}\\
& &+a^{\Lm_{j-1}+\Lm_{k-1}+2\Lm_{r}}\sum_{b\in B(\Lm_j+\Lambda_{k-1}+ \Lambda_{r-1})_{s_{j+1}\cdots s_{k-1}s_{k+1}s_{k+2}\cdots s_{r-1}}} \mu'(b)\nonumber,
\end{eqnarray}
where $\mu$ is the monomial realization of $B(\Lm_j+\Lambda_k+2\Lambda_{r})$ such that the highest weight vector is realized by $Y_{2,j}Y_{2,k}Y^2_{1,r}\in\mathcal{Y}$, $\mu'$ is the monomial realization of $B(\Lm_j+\Lambda_{k-1}+\Lambda_{r-1})$ such that the highest weight vector is realized by $Y_{2,j}Y_{1,k-1}Y_{1,r-1}\in\mathcal{Y}$, and $C(b)$ are some positive integers.
\end{enumerate}

\end{thm}

\begin{thm}\label{thmC}
Let $G={\rm Sp}_{2r}(\mathbb{C})$ be the classical algebraic group of type ${\rm C}_r$. 
\begin{enumerate}
\item For $k\in[1,r]$, we obtain
\[ 
(\varphi^G_{\textbf{V}})_k(a;\textbf{Y})=a^{\Lambda_k} \sum_{b\in B(\Lambda_k)_{c^2_{>2r-k}}} C(b)\mu(b),
\]
where $\mu:B(\Lambda_k)\rightarrow \mathcal{Y}$ is the monomial realization of $B(\Lambda_k)$ such that the highest weight vector in $B(\Lambda_k)$ is realized by $Y_{1,k}\in\mathcal{Y}$, and $C(b)$ are some positive integers. We also obtain
\[
(\varphi^G_{\mu_{k}\cdots\mu_{r-1}\mu_r(\textbf{V})})_k(a;\textbf{Y})=a^{\Lambda_{k-1}+\Lambda_{r-1}} Y_{2,k}=a^{\Lambda_{k-1}+\Lambda_{r-1}} \sum_{b\in B(\Lambda_k)_{e}} \mu'(b),
\]
where $\mu'$ is the monomial realization of $B(\Lambda_k)$ such that the highest weight vector in $B(\Lambda_k)$ is realized by $Y_{2,k}\in\mathcal{Y}$.
\item  For $k$ and $l$ with $1\leq k\leq l\leq r-2$,
\begin{eqnarray}
& &\hspace{-20pt}(\varphi^G_{\mu_{l}\cdots\mu_{k+1}\mu_k(\textbf{V})})_l(a;\textbf{Y}) \nonumber\\ 
&=&
a^{\Lm_{k-1}+\Lm_{l+1}} Y_{2,k}Y_{1,l+1} (1+A^{-1}_{1,l+1}+A^{-1}_{1,l+1}A^{-1}_{1,l}+\cdots+A^{-1}_{1,l+1}A^{-1}_{1,l}\cdots A^{-1}_{1,k+1})
\qquad \quad \label{C-exp1}\\
&=&a^{\Lm_{k-1}+\Lm_{l+1}}\sum_{b\in B(\Lambda_k+\Lambda_{l+1})_{s_{k+1}s_{k+2}\cdots s_{l+1}}} \mu(b), \label{C-cry1}
\end{eqnarray}
where $\mu$ is the monomial realization of $B(\Lambda_k+\Lambda_{l+1})$ such that the highest weight vector in $B(\Lambda_k+\Lambda_{l+1})$ is realized by $Y_{2,k}Y_{1,l+1}\in\mathcal{Y}$.
\item For $k\in[1,r-1]$, we obtain
\begin{eqnarray}
& &\hspace{-20pt}(\varphi^G_{\mu_{r-1}\cdots\mu_{k+1}\mu_k(\textbf{V})})_{r-1}(a;\textbf{Y})\nonumber\\
&=&a^{\Lm_{k-1}+\Lm_r}\left(Y_{2,k}Y_{1,r}+Y_{2,k}Y_{1,r}A^{-1}_{1,r}(1+A^{-1}_{1,r-1}+\cdots+A^{-1}_{1,r-1}\cdots A^{-1}_{1,k+1})^2\right.\nonumber\\
&+&Y_{1,k-1}Y_{1,r-1}(1+A^{-1}_{1,r-1}+\cdots+A^{-1}_{1,r-1}\cdots A^{-1}_{1,k+1}) \nonumber\\
& &\left.\qquad \qquad \qquad \qquad \times (1+A^{-1}_{1,k-1}+\cdots+A^{-1}_{1,k-1}\cdots A^{-1}_{1,1})\right)\nonumber \\
&=&a^{\Lm_{k-1}+\Lm_{r}} \sum_{b\in B(\Lambda_k+\Lambda_{r})_{s_{k+1}s_{k+2}\cdots s_{r}}} C(b)\mu(b)\nonumber\\
& &+a^{\Lm_{k-1}+\Lm_{r}}\sum_{b\in B(\Lambda_{k-1}+ \Lambda_{r-1})_{s_1\cdots s_{k-1}s_{k+1}s_{k+2}\cdots s_{r-1}}} \mu'(b), \label{C-cry2}
\end{eqnarray}
where $\mu$ is the monomial realization of $B(\Lambda_k+\Lambda_{r})$ such that the highest weight vector is realized by $Y_{2,k}Y_{1,r}\in\mathcal{Y}$, $\mu'$ is the monomial realization of $B(\Lambda_{k-1}+\Lambda_{r-1})$ such that the highest weight vector is realized by $Y_{1,k-1}Y_{1,r-1}\in\mathcal{Y}$, and $C(b)$ are some positive integers. We also obtain
\begin{eqnarray}
& &\hspace{-20pt}(\varphi^G_{\mu_r\mu_{r-1}\cdots\mu_{k+1}\mu_k(\textbf{V})})_{r}(a;\textbf{Y})\nonumber\\
&=&
a^{\Lm_r+2\Lm_{k-1}}Y_{1,r}Y^2_{2,k}(1+A^{-1}_{1,r}(1+A^{-1}_{1,r-1}+\cdots+A^{-1}_{1,r-1}\cdots A^{-1}_{1,k+1})^2)\qquad \quad \label{C-lem2-pr8}\\
&=&a^{2\Lm_{k-1}+\Lm_{r}} \sum_{b\in B(2\Lambda_k+\Lambda_{r})_{s_{k+1}s_{k+2}\cdots s_{r}}} C'(b)\mu''(b),\label{C-cry3}
\end{eqnarray}
where $\mu''$ is the monomial realization of $B(2\Lambda_k+\Lambda_{r})$ such that the highest weight vector is realized by $Y^2_{2,k}Y_{1,r}\in\mathcal{Y}$, and $C'(b)$ are some positive integers.
\item For $j$ and $k$ with $1\leq j< k\leq r-1$,
\begin{eqnarray}
& &\hspace{-20pt}(\varphi^G_{\mu_j\cdots\mu_{k-2}\mu_{k-1}\mu_{r}\mu_{r-1}\cdots\mu_{k+1}\mu_k(\textbf{V})})_j(a;\textbf{Y})= a^{\Lm_{j-1}+\Lm_{k-1}+\Lm_{r}}\left(Y_{2,k}Y_{2,j}Y_{1,r}\nonumber\right.\\
&+&\frac{Y_{2,k}Y_{2,j}Y^2_{1,r-1}}{Y_{2,r}}(1+A^{-1}_{1,r-1}+A^{-1}_{1,r-1}A^{-1}_{1,r-2}+\cdots+A^{-1}_{1,r-1}\cdots A^{-1}_{1,k+1})\ \ \label{C-lem3-exp}\\
& &\left. \qquad\qquad\quad \times
(1+A^{-1}_{1,r-1}+A^{-1}_{1,r-1}A^{-1}_{1,r-2}+\cdots+A^{-1}_{1,r-1}\cdots A^{-1}_{1,j+1})\right)
\quad\quad \quad \nonumber\\
&=& a^{\Lm_{j-1}+\Lm_{k-1}+\Lm_{r}} \sum_{b\in B(\Lm_j+\Lambda_k+\Lambda_{r})_{s_{k+1}s_{k+2}\cdots s_{r}} } C(b)\mu(b)\label{C-lem3-cry}\\
& &+a^{\Lm_{j-1}+\Lm_{k-1}+\Lm_{r}}\sum_{b\in B(\Lm_j+\Lambda_{k-1}+ \Lambda_{r-1})_{s_{j+1}\cdots s_{k-1}s_{k+1}s_{k+2}\cdots s_{r-1}}} \mu'(b)\nonumber,
\end{eqnarray}
where $\mu$ is the monomial realization of $B(\Lm_j+\Lambda_k+\Lambda_{r})$ such that the highest weight vector is realized by $Y_{2,j}Y_{2,k}Y_{1,r}\in\mathcal{Y}$, $\mu'$ is the monomial realization of $B(\Lm_j+\Lambda_{k-1}+\Lambda_{r-1})$ such that the highest weight vector is realized by $Y_{2,j}Y_{1,k-1}Y_{1,r-1}\in\mathcal{Y}$, and $C(b)$ are some positive integers.
\end{enumerate}

\end{thm}

\begin{thm}\label{thmD}
Let $G={\rm SO}_{2r}(\mathbb{C})$ be the classical algebraic group of type ${\rm D}_r$. 
\begin{enumerate}
\item For $k\in[1,r]$, we obtain
\[ 
(\varphi^G_{\textbf{V}})_k(a;\textbf{Y})=a^{\Lambda_k} \sum_{b\in B(\Lambda_k)_{c^2_{>2r-k}}} C(b)\mu(b),
\]
where $\mu:B(\Lambda_k)\rightarrow \mathcal{Y}$ is the monomial realization of $B(\Lambda_k)$ such that the highest weight vector in $B(\Lambda_k)$ is realized by $Y_{1,k}\in\mathcal{Y}$, and $C(b)$ are some positive integers. Furthermore, we have
\[ 
(\varphi^G_{\mu_{r-1}(\textbf{V})})_{r-1}(a;\textbf{Y})=
a^{\Lambda_{r-2}} Y_{2,r-1}=a^{\Lambda_{r-2}}\sum_{b\in B(\Lambda_{r-1})_{e}} \mu'_{r-1}(b),\]
\[
(\varphi^G_{\mu_{r}\mu_{r-1}(\textbf{V})})_{r}(a;\textbf{Y})=
a^{\Lambda_{r-2}} Y_{2,r}=a^{\Lambda_{r-2}}\sum_{b\in B(\Lambda_{r})_{e}} \mu'_{r}(b).
\]
For $k\in[1,r-2]$, 
\[ 
(\varphi^G_{\mu_k\mu_{k+1}\cdots\mu_{r-3}\mu_{r-2}\mu_r\mu_{r-1}(\textbf{V})})_k(a;\textbf{Y})=
a^{\Lambda_{k-1}+\Lambda_{r-2}} Y_{2,k}=a^{\Lambda_{k-1}+\Lambda_{r-2}} 
\sum_{b\in B(\Lambda_{k})_{e}} \mu'_{k}(b),
\]
where $\mu'_s$ $(1\leq s\leq r)$ is the monomial realization of $B(\Lm_s)$ such that the highest weight vector is realized by $Y_{2,s}\in\cY$.
\item  For $k$ and $l$ with $1\leq k\leq l\leq r-3$,
\begin{eqnarray}
& &\hspace{-20pt}(\varphi^G_{\mu_{l}\cdots\mu_{k+1}\mu_k(\textbf{V})})_l(a;\textbf{Y}) \nonumber\\ 
&=&
a^{\Lm_{k-1}+\Lm_{l+1}} Y_{2,k}Y_{1,l+1} (1+A^{-1}_{1,l+1}+A^{-1}_{1,l+1}A^{-1}_{1,l}+\cdots+A^{-1}_{1,l+1}A^{-1}_{1,l}\cdots A^{-1}_{1,k+1})
\qquad \quad \label{D-exp1}\\
&=&a^{\Lm_{k-1}+\Lm_{l+1}}\sum_{b\in B(\Lambda_k+\Lambda_{l+1})_{s_{k+1}s_{k+2}\cdots s_{l+1}}} \mu(b), \label{D-cry1}
\end{eqnarray}
where $\mu$ is the monomial realization of $B(\Lambda_k+\Lambda_{l+1})$ such that the highest weight vector in $B(\Lambda_k+\Lambda_{l+1})$ is realized by $Y_{2,k}Y_{1,l+1}\in\mathcal{Y}$.
\item For $k\in[1,r-2]$, we obtain
\begin{eqnarray}
& &\hspace{-20pt}(\varphi^G_{\mu_{r-2}\cdots\mu_{k+1}\mu_k(\textbf{V})})_{r-2}(a;\textbf{Y})\nonumber\\
&=&a^{\Lm_{k-1}+\Lm_{r-1}+\Lm_r}Y_{2,k}Y_{1,r-1}Y_{1,r}
(1+A^{-1}_{1,r-1}+\cdots+A^{-1}_{1,r-1}\cdots A^{-1}_{1,k+1})\nonumber\\
& &\times (1+A^{-1}_{1,r}+A^{-1}_{1,r}A^{-1}_{1,r-2}+ \cdots+A^{-1}_{1,r}A^{-1}_{1,r-2}A^{-1}_{1,r-3}\cdots A^{-1}_{1,k+1})\nonumber\\
& &+a^{\Lm_{k-1}+\Lm_{r-1}+\Lm_r}Y_{1,k-1}Y_{1,r-2}
(1+A^{-1}_{1,r-2}+\cdots+A^{-1}_{1,r-2}\cdots A^{-1}_{1,k+1})\qquad \quad \label{D-lem2-pr2-1}\\
& &\times (1+A^{-1}_{1,k-1}+A^{-1}_{1,k-1}A^{-1}_{1,k-2}+\cdots+A^{-1}_{1,k-1}A^{-1}_{1,k-2}\cdots A^{-1}_{1,1})\nonumber\\
&=&a^{\Lm_{k-1}+\Lm_{r-1}+\Lm_r} \sum_{b\in B(\Lm_k+\Lm_{r-1}+\Lm_r)_{s_{k+1}\cdots s_{r-1}s_r}} C(b)\mu(b)\label{D-lem2-cry}\\
& & + a^{\Lm_{k-1}+\Lm_{r-1}+\Lm_r} \sum_{b\in B(\Lm_{k-1}+\Lm_{r-2})_{s_1\cdots s_{k-1}s_{k+1}\cdots s_{r-3}s_{r-2}}} \mu'(b), \nonumber
\end{eqnarray}
where $\mu$ is the monomial realization of $B(\Lambda_k+\Lambda_{r-1}+\Lm_r)$ such that the highest weight vector is realized by $Y_{2,k}Y_{1,r-1}Y_{1,r}\in\mathcal{Y}$ and $C(b)$ are some integers, $\mu'$ is the monomial realization of $B(\Lambda_{k-1}+\Lambda_{r-2})$ such that the highest weight vector is realized by $Y_{1,k-1}Y_{1,r-2}\in\mathcal{Y}$.
\item For $k\in[1,r-2]$,
\begin{eqnarray}
& &\hspace{-20pt}(\varphi^G_{\mu_{r-1}\mu_{r-2}\cdots\mu_{k+1}\mu_k(\textbf{V})})_{r-1}(a;\textbf{Y})\nonumber\\
&=&a^{\Lm_{k-1}+\Lm_r}Y_{2,k}Y_{1,r}
(1+A^{-1}_{1,r}+A^{-1}_{1,r}A^{-1}_{1,r-2}+ \cdots+A^{-1}_{1,r}A^{-1}_{1,r-2}A^{-1}_{1,r-3}\cdots A^{-1}_{1,k+1}),\qquad \quad \label{D-lem2-pr3-1}\\
&=&a^{\Lm_{k-1}+\Lm_r} \sum_{b\in B(\Lm_k+\Lm_r)_{s_{k+1}\cdots s_{r-1}s_r}} \mu(b),\label{D-lem3-cry}
\end{eqnarray}
where $\mu$ is the monomial realization of $B(\Lambda_k+\Lm_r)$ such that the highest weight vector is realized by $Y_{2,k}Y_{1,r}\in\mathcal{Y}$.
\begin{eqnarray}
& &\hspace{-20pt}(\varphi^G_{\mu_{r}\mu_{r-1}\cdots\mu_{k+1}\mu_k(\textbf{V})})_{r}(a;\textbf{Y})\nonumber\\
&=&
a^{\Lm_{k-1}+\Lm_{r-1}}Y_{2,k}Y_{1,r-1}
(1+A^{-1}_{1,r-1}+A^{-1}_{1,r-1}A^{-1}_{1,r-2}+ \cdots+A^{-1}_{1,r-1}A^{-1}_{1,r-2}\cdots A^{-1}_{1,k+1})\qquad \quad \label{D-lem2-pr3-2}\\
&=&a^{\Lm_{k-1}+\Lm_{r-1}} \sum_{b\in B(\Lm_k+\Lm_{r-1})_{s_{k+1}s_{k+2}\cdots s_{r-1}}} \mu'(b),\label{D-lem4-cry}
\end{eqnarray}
where $\mu'$ is the monomial realization of $B(\Lambda_k+\Lm_{r-1})$ such that the highest weight vector is realized by $Y_{2,k}Y_{1,r-1}\in\mathcal{Y}$.
\item For $j$ and $k$ with $1\leq j< k\leq r-2$,
\begin{eqnarray}
& &\hspace{-20pt}(\varphi^G_{\mu_j\cdots\mu_{k-2}\mu_{k-1}\mu_{r}\mu_{r-1}\cdots\mu_{k+1}\mu_k(\textbf{V})})_j(a;\textbf{Y})\nonumber\\
&=&a^{\Lm_{j-1}+\Lm_{k-1}+\Lm_{r-1}+\Lm_r} Y_{2,j}Y_{2,k}Y_{1,r-1}Y_{1,r}
(1+A^{-1}_{1,r-1}+\cdots+A^{-1}_{1,r-1}A^{-1}_{1,r-2}\cdots A^{-1}_{1,k+1})\nonumber\\
& &\times 
(1+A^{-1}_{1,r}+A^{-1}_{1,r}A^{-1}_{1,r-2}+ \cdots+A^{-1}_{1,r}A^{-1}_{1,r-2}A^{-1}_{1,r-3}\cdots A^{-1}_{1,k+1})\nonumber\\
&+&a^{\Lm_{j-1}+\Lm_{k-1}+\Lm_{r-1}+\Lm_r} Y_{2,j}Y_{1,k-1}Y_{1,r-2}
(1+A^{-1}_{1,r-2}+\cdots+A^{-1}_{1,r-2}A^{-1}_{1,r-3}\cdots A^{-1}_{1,k+1})\qquad \quad \label{D-lem3-exp1}\\
& &\times (1+A^{-1}_{1,k-1}+\cdots+A^{-1}_{1,k-1}A^{-1}_{1,k-2}\cdots A^{-1}_{1,j+1})\nonumber\\
&=&a^{\Lm_{j-1}+\Lm_{k-1}+\Lm_{r-1}+\Lm_r} \sum_{b\in B(\Lm_j+\Lambda_k+\Lm_{r-1}+\Lambda_{r})_{s_{k+1}\cdots s_{r-1}s_r}}C(b)\mu(b)\nonumber\\
&+&a^{\Lm_{j-1}+\Lm_{k-1}+\Lm_{r-1}+\Lm_r} \sum_{B(\Lm_j+\Lambda_{k-1}+\Lm_{r-2})_{s_{j+1}\cdots s_{k-2}s_{k-1} s_{k+1}\cdots s_{r-3}s_{r-2}}}\mu'(b), \label{D-lem3-cry1}
\end{eqnarray}
where $\mu$ is the monomial realization of $B(\Lm_j+\Lambda_k+\Lm_{r-1}+\Lambda_{r})$ such that the highest weight vector is realized by $Y_{2,j}Y_{2,k}Y_{1,r-1}Y_{1,r}\in\cY$, and $\mu'$ is the monomial realization of $B(\Lm_j+\Lambda_{k-1}+\Lm_{r-2})$ such that the highest weight vector is mapped to $Y_{2,j}Y_{1,k-1}Y_{1,r-2}\in\cY$, and $C(b)$ are some positive integers.
\end{enumerate}

\end{thm}

\begin{rem}
In Theorem \ref{thmB} (resp. \ref{thmC}, \ref{thmD}),
we saw explicit formulas of $r^2+r$ (resp. $r^2+r$, $r^2$)
cluster variables 
in $\mathbb{C}[G^{e,c^2}]$.
The number $r^2+r$ (resp. $r^2+r$, $r^2$) coincides 
with 
$|\Phi_{\geq -1}|$ of type ${\rm B}_r$ 
(resp. ${\rm C}_r$,${\rm D}_r$).
Thus, by Theorem \ref{finthm2}, all cluster variables 
in $\mathbb{C}[G^{e,c^2}]$ appear
in 
Theorem \ref{thmB}-\ref{thmD}.
\end{rem}

\section{The proof of main theorem}\label{mainthmpr}

In this section, we prove Theorem \ref{thmB}, \ref{thmC} and \ref{thmD}. Let $\Sigma_0:=(\textbf{V},\tilde{B}(\textbf{i}))$ be the initial seed of $\mathbb{C}[G^{e,c^2}]$.

\subsection{The proof of Theorem \ref{thmB}}

First, we shall prove the case $G={\rm SO}_{2r+1}(\mathbb{C})$. We start by setting the Laurent monomials as follows:
\begin{equation}\label{bbbar}
B(l,k):=
\begin{cases}
\frac{Y_{l,k}}{Y_{l,k-1}} & {\rm if}\ 1\leq k\leq r-1, \\
\frac{Y^2_{l,r}}{Y_{l,r-1}} & {\rm if}\ k= r, \\
\frac{Y_{l,r}}{Y_{l+1,r}} & {\rm if}\ k= 0, \\
\frac{Y_{l,r-1}}{Y^2_{l+1,r}} & {\rm if}\ k= \ovl{r}, \\
\frac{Y_{l,|k|-1}}{Y_{l+1,|k|}} & {\rm if}\ \ovl{r-1}\leq k\leq \ovl{1},
\end{cases}
\end{equation}
where for $1\leq l\leq r$, we set $|l|=|\ovl{l}|=l$. 

\begin{prop}\label{iniB}
\begin{enumerate}
\item For $k\in[1,r]$,
the initial cluster variables $(\varphi^G_{\textbf{V}})_k(a;\textbf{Y})$ in $\mathbb{C}[G^{e,c^2}]$ are
described as
\[ 
(\varphi^G_{\textbf{V}})_k(a;\textbf{Y})=a^{\Lambda_k} \sum_{b\in B(\Lambda_k)_{c^2_{>2r-k}}} \mu(b),
\]
where $\mu:B(\Lambda_k)\rightarrow \mathcal{Y}$ is the monomial realization of $B(\Lambda_k)$ in Theorem \ref{thmB} $(i)$.
\item
For $k\in[1,r]$, the frozen cluster variables $(\varphi^G_{\textbf{V}})_{-k}(a;\textbf{Y})$, $(\varphi^G_{\textbf{V}})_{r+k}(a;\textbf{Y})$ in $\mathbb{C}[G^{e,c^2}]$ are described as
\[
(\varphi^G_{\textbf{V}})_{-k}(a;\textbf{Y})=a^{\Lm_k}Y_{1,k}Y_{2,k},
\quad
(\varphi^G_{\textbf{V}})_{r+k}(a;\textbf{Y})=a^{\Lm_k}.
\]
\end{enumerate}
\end{prop}
\nd
[{\sl Proof.}]

First, let $k\in[1,r-1]$, and recall the fundamental representation $V(\Lambda_k)$ of type ${\rm B}_r$ in \ref{SectFundB}. By (\ref{B-f1}), (\ref{B-f2}) and (\ref{smpl}), for $i\in[1,r-1]$, $j\in[1,r]$, we get
\[ \ovl{s_i} v_j=
\begin{cases}
v_{i+1} & {\rm if}\ j=i, \\
-v_{i} & {\rm if}\ j=i+1, \\
v_j & {\rm if}\ {\rm otherwise},
\end{cases}\q 
\ovl{s_r} v_j=
\begin{cases}
v_{\ovl{r}} & {\rm if}\ j=r, \\
v_j & {\rm if}\ 1\leq j< r.
\end{cases}
\]
Taking into account these formulas, we obtain
\begin{equation}\label{B-ext}
\ovl{c^2_{>2r-k}}v_1\wedge\cdots\wedge v_k=\ovl{s_1}\cdots \ovl{s_k}(v_1\wedge\cdots\wedge v_k)
=v_2\wedge\cdots\wedge v_{k+1}. 
\end{equation}
Just as in Example \ref{initialex4}, the value of $(\varphi^G_{\textbf{V}})_k(a;\textbf{Y})$ coincides with
\begin{multline*}
a^{\Lm_k}\lan v_2\wedge\cdots\wedge v_{k+1}, \al^{\vee}_1(Y_{2,1})y_1(Y_{2,1})\cdots \al^{\vee}_r(Y_{2,r})y_r(Y_{2,r})\\ \al^{\vee}_1(Y_{1,1})y_1(Y_{1,1})\cdots\al^{\vee}_k(Y_{1,k})y_k(Y_{1,k})v_1\wedge\cdots\wedge v_k \ran.
\end{multline*}

Using $(\ref{B-wtv})$, $(\ref{B-f1})$ and $(\ref{B-f2})$ repeatedly, for $i\in[1,k]$, one obtain
\[
 \al^{\vee}_1(Y_{1,1})y_1(Y_{1,1})\cdots\al^{\vee}_k(Y_{1,k})y_k(Y_{1,k})v_i=
\frac{Y_{1,i}}{Y_{1,i-1}}v_i+v_{i+1},
\]
which means that $\al^{\vee}_1(Y_{1,1})y_1(Y_{1,1})\cdots\al^{\vee}_k(Y_{1,k})y_k(Y_{1,k})v_1\wedge\cdots\wedge v_k$ is a linear combination of $v_1\wedge\cdots \wedge v_s\wedge v_{s+2}\wedge\cdots \wedge v_{k+1}$ $(0\leq s\leq k)$ with the coefficient $B(1,1)B(1,2)\cdots B(1,s)$. Similarly, we can also verify that for $s$ $(0\leq s\leq k)$ the coefficient of $v_2\wedge\cdots\wedge v_{k+1}$ in $\al^{\vee}_1(Y_{2,1})y_1(Y_{2,1})\cdots\al^{\vee}_r(Y_{2,r})y_r(Y_{2,r})(v_1\wedge\cdots \wedge v_s\wedge v_{s+2}\wedge\cdots \wedge v_{k+1})$ is $B(2,s+2)B(2,s+3)\cdots B(2,k+1)$. Hence, we get
\begin{equation}\label{iniBexp}
(\varphi^G_{\textbf{V}})_k(a;\textbf{Y})
=a^{\Lambda_k} \sum_{0\leq s\leq k}B(1,1)B(1,2)\cdots B(1,s)B(2,s+2)\cdots B(2,k+1).
\end{equation}
The definition of the monomial realization implies that
\begin{eqnarray*}
& &\tilde{f}_s B(1,1)\cdots B(1,s-1)B(1,s)B(2,s+2)\cdots B(2,k+1)\\
&=&B(1,1)\cdots B(1,s-1)B(1,s)B(2,s+2)\cdots B(2,k+1)A^{-1}_{1,s}\\
&=&B(1,1)\cdots B(1,s-1)B(2,s+1)B(2,s+2)\cdots B(2,k+1),
\end{eqnarray*}
where we use $B(1,s)A^{-1}_{1,s}=\frac{Y_{1,s}}{Y_{1,s-1}}\frac{Y_{1,s-1}Y_{2,s+1}}{Y_{1,s}Y_{2,s}}=B(2,s+1)$.
Therefore, the conclusion
\begin{eqnarray}
(\varphi^G_{\textbf{V}})_k(a;\textbf{Y})&=&
a^{\Lambda_k} Y_{1,k}(1+A^{-1}_{1,k}+A^{-1}_{1,k}A^{-1}_{1,k-1}+\cdots+A^{-1}_{1,k}A^{-1}_{1,k-1}\cdots A^{-1}_{1,1})\quad \quad \quad \nonumber\\
&=&a^{\Lambda_k} \sum_{b\in B(\Lambda_k)_{c^2_{>2r-k}}} \mu(b)\nonumber
\end{eqnarray}
follows from Theorem \ref{kashidem} and the easy fact $B(1,1)B(1,2)\cdots B(1,k)=Y_{1,k}$. By the same argument, the frozen cluster variable $(\varphi^G_{\textbf{V}})_{-k}(a;\textbf{Y})=\Delta_{\Lm_k,\Lm_k}\circ \ovl{x}^G(a;\textbf{Y})$ coincides with the coefficient of $v_1\wedge \cdots \wedge v_k$ in 
\begin{equation}\label{iniB-pr-1}
a \al^{\vee}_1(Y_{2,1})y_1(Y_{2,1})\cdots \al^{\vee}_r(Y_{2,r})y_r(Y_{2,r}) \al^{\vee}_1(Y_{1,1})y_1(Y_{1,1})\cdots\al^{\vee}_k(Y_{1,k})y_k(Y_{1,k})v_1\wedge\cdots\wedge v_k,
\end{equation}
which equals to $a^{\Lm_k}B(1,1)B(1,2)\cdots B(1,k)B(2,1)B(2,2)\cdots B(2,k)=a^{\Lm_k}Y_{1,k}Y_{2,k}$. The frozen cluster variable $(\varphi^G_{\textbf{V}})_{r+k}(a;\textbf{Y})=\Delta_{\Lm_k,c^{-2}\Lm_k}\circ \ovl{x}^G(a;\textbf{Y})$ coincides with the coefficient of 
\[
\begin{cases}
v_3\wedge \cdots \wedge v_{k+2} & {\rm if}\ k<r-1,\\
v_3\wedge\cdots\wedge v_r\wedge v_{\ovl{1}} & {\rm if}\ k=r-1,
\end{cases}
\]
in (\ref{iniB-pr-1}). It is equal to $a^{\Lm_k}$.

Next, we now turn to the case $k=r$. We need to recall the spin representation in \ref{SectFundB}. From (\ref{Bsp-f1}), (\ref{Bsp-f2}) and (\ref{smpl}), we see that
\[
\ovl{c^2_{>r}}(+,+,+,\cdots,+)=\ovl{s_1}\cdots \ovl{s_r}(+,+,+,\cdots,+)
=(-,+,+,\cdots,+). 
\]
Just as in the case $k<r$, the value of $(\varphi^G_{\textbf{V}})_r(a;\textbf{Y})$ coincides with
\begin{multline*}
a^{\Lm_r}\lan (-,+,+,\cdots,+), \al^{\vee}_1(Y_{2,1})y_1(Y_{2,1})\cdots \al^{\vee}_r(Y_{2,r})y_r(Y_{2,r})\\ \al^{\vee}_1(Y_{1,1})y_1(Y_{1,1})\cdots\al^{\vee}_r(Y_{1,r})y_r(Y_{1,r})(+,+,+,\cdots,+) \ran.
\end{multline*}

From (\ref{Bsp-f0}) and (\ref{Bsp-f1}), we see that $\al^{\vee}_1(Y_{1,1})y_1(Y_{1,1})\cdots\al^{\vee}_r(Y_{1,r})y_r(Y_{1,r})(+,+,+,\cdots,+)$ is a linear combination of $(+,+,\cdots,+)$ and $(+,\cdots,+,\overset{i}{-},\overset{i+1}{+},\cdots,+)$ $(1\leq i\leq r)$ whose coefficients are $Y_{1,r}$ and $Y_{1,i-1}$, respectively. Similarly, the coefficient of $(-,+,\cdots,+,+)$ in $\al^{\vee}_1(Y_{2,1})y_1(Y_{2,1})\cdots \al^{\vee}_r(Y_{2,r})y_r(Y_{2,r})(+,+,\cdots,+,+)$ is $1$, the one in $\al^{\vee}_1(Y_{2,1})y_1(Y_{2,1})\cdots \al^{\vee}_r(Y_{2,r})y_r(Y_{2,r})(+,\cdots,+,\overset{i}{-},\overset{i+1}{+},\cdots,+)$ is $\frac{Y_{2,r}}{Y_{2,i}}$ if $i<r$, and is $\frac{1}{Y_{2,r}}$ if $i=r$. Hence, it follows that
\begin{eqnarray}
(\varphi^G_{\textbf{V}})_r(a;\textbf{Y})&=&a^{\Lambda_r} \left(Y_{1,r}+\frac{Y_{1,r-1}}{Y_{2,r}}+\sum^{r-1}_{i=1} \frac{Y_{1,i-1}Y_{2,r}}{Y_{2,i}}\right) \label{iniBexp-r} \\
&=&a^{\Lambda_r}\sum_{b\in B(\Lambda_r)_{c^2_{>r}}} \mu(b). \nonumber
\end{eqnarray}

By the same argument, we can prove that  $(\varphi^G_{\textbf{V}})_{-r}(a;\textbf{Y})=a^{\Lm_r}Y_{1,r}Y_{2,r}$ and $(\varphi^G_{\textbf{V}})_{2r}(a;\textbf{Y})=a^{\Lm_r}$.
\qed

In the proof, we found the explicit form (\ref{iniBexp}) of $(\varphi^G_{\textbf{V}})_k(a;\textbf{Y})$ for $k\in[1,r-1]$, which can be rewritten as
\begin{eqnarray*}
& &\hspace{-20pt}(\varphi^G_{\textbf{V}})_k(a;\textbf{Y})=a^{\Lambda_k} \sum_{0\leq s\leq k}B(1,1)B(1,2)\cdots B(1,s)B(2,s+2)\cdots B(2,k+1)\\
&=&a^{\Lambda_k} (B(1,1)\cdots B(1,k)+\sum_{0\leq s\leq k-1}B(1,1)\cdots B(1,s)B(2,s+2)\cdots B(2,k+1))\\
&=&a^{\Lambda_k}(Y_{1,k}+B(2,k+1)\sum_{0\leq s\leq k-1}B(1,1)\cdots B(1,s)B(2,s+2)\cdots B(2,k))\\
&=&a^{\Lambda_k}(Y_{1,k}+B(2,k+1)(a^{-\Lm_{k-1}}(\varphi^G_{\textbf{V}})_{k-1}(a;\textbf{Y}))).
\end{eqnarray*}
Thus, we obtain the relation between two cluster variables:
\begin{equation}\label{twovarsB}
a^{-\Lambda_k}(\varphi^G_{\textbf{V}})_k(a;\textbf{Y})=Y_{1,k}+B(2,k+1)(a^{-\Lm_{k-1}}(\varphi^G_{\textbf{V}})_{k-1}(a;\textbf{Y})).
\end{equation}

\begin{lem}\label{B-lem1}
For $k\in[1,r]$, the cluster variables $(\varphi^G_{\mu_k\mu_{k+1}\cdots\mu_r(\textbf{V})})_k(a;\textbf{Y})$ in $\mathbb{C}[G^{e,c^2}]$ are
described as
\[ 
(\varphi^G_{\mu_k\mu_{k+1}\cdots\mu_r(\textbf{V})})_k(a;\textbf{Y})=
\begin{cases}
a^{\Lambda_{k-1}+\Lambda_{r-1}} Y_{2,k} & {\rm if}\ k<r,\\
a^{\Lambda_{r-1}} Y_{2,r} & {\rm if}\ k=r.
\end{cases}
\]
\end{lem}
\nd
[{\sl Proof.}]

We prove this statement by induction on $(r-k)$. If $r-k=0$, so that $k=r$, then the mutation diagram (\ref{inidiaB-2}) implies that
\[
(\varphi^G_{\mu_r(\textbf{V})})_r(a;\textbf{Y})
=\frac{(\varphi^G_{\textbf{V}})_{r-1}(a;\textbf{Y})(\varphi^G_{\textbf{V}})_{2r}(a;\textbf{Y})+(\varphi^G_{\textbf{V}})_{2r-1}(a;\textbf{Y})(\varphi^G_{\textbf{V}})_{-r}(a;\textbf{Y})}{(\varphi^G_{\textbf{V}})_r(a;\textbf{Y})}.
\]
Applying Proposition \ref{iniB} (ii), (\ref{iniBexp2}) for $k=r-1$ and (\ref{iniBexp-r}), we obtain
\begin{eqnarray*}
& &\hspace{-20pt}(\varphi^G_{\mu_r(\textbf{V})})_r(a;\textbf{Y})\\
&=& \frac{a^{\Lambda_{r-1}+\Lambda_r}}{a^{\Lambda_r}}
\frac{Y_{1,r-1}(1+A^{-1}_{1,r-1}+A^{-1}_{1,r-1}A^{-1}_{1,r-2}+\cdots+A^{-1}_{1,r-1}A^{-1}_{1,r-2}\cdots A^{-1}_{1,1})+Y_{1,r}Y_{2,r}}{Y_{1,r}+\frac{Y_{1,r-1}}{Y_{2,r}}+\sum^{r-1}_{i=1} \frac{Y_{1,i-1}Y_{2,r}}{Y_{2,i}}}\\
&=& a^{\Lambda_{r-1}} \frac{Y_{1,r}Y_{2,r}+Y_{1,r-1}(1+A^{-1}_{1,r-1}+A^{-1}_{1,r-1}A^{-1}_{1,r-2}+\cdots+A^{-1}_{1,r-1}A^{-1}_{1,r-2}\cdots A^{-1}_{1,1})}{Y_{1,r}+\frac{Y_{1,r-1}}{Y_{2,r}}(1+A^{-1}_{1,r-1}+A^{-1}_{1,r-1}A^{-1}_{1,r-2}+\cdots+A^{-1}_{1,r-1}A^{-1}_{1,r-2}\cdots A^{-1}_{1,1})}\\
&=& a^{\Lambda_{r-1}} Y_{2,r}.
\end{eqnarray*}

Next, we consider the the cluster variable $(\varphi^G_{\mu_{r-1}\mu_r(\textbf{V})})_{r-1}(a;\textbf{Y})$. In the rest of this proof, we abbreviate $(\varphi^G_{\textbf{T}})_{s}(a;\textbf{Y})$ to $(\varphi^G_{\textbf{T}})_{s}$ for $s\in[1,r]$ and clusters $\textbf{T}$.
By Lemma \ref{mutgamlem}, the mutation diagram of $\mu_r(\Sigma_0)$ is as follows:
\[
\begin{xy}
(60,98)*{(\varphi_{\textbf{V}})_{2r}}="r+k",
(60,90) *{\ (\varphi^G_{\mu_r(\textbf{V})})_r}="k",
(60,82)*{(\varphi_{\textbf{V}})_{-r}}="-j_k",
(35,98)*{(\varphi_{\textbf{V}})_{2r-1}}="8",
(35,90) *{(\varphi_{\textbf{V}})_{r-1}}="4",
(35,82)*{(\varphi_{\textbf{V}})_{-(r-1)}}="-2",
(10,98)*{(\varphi_{\textbf{V}})_{2r-2}}="6",
(10,90) *{(\varphi_{\textbf{V}})_{r-2}}="2",
(10,82)*{(\varphi_{\textbf{V}})_{-(r-2)}}="-1",
(0,90)*{\cdots}="emp",
\ar@{->} "k";"r+k"
\ar@{->} "4";"8"
\ar@{->} "6";"2"
\ar@{->} "-j_k";"k"
\ar@{->} "8";"k"^{-2}
\ar@{->} "4";"-2"
\ar@{->} "2";"-1"
\ar@{->} "k";"4"^{2}
\ar@{->} "2";"4"
\ar@{->} "4";"6"
\ar@{->} "-2";"2"
\end{xy}
\]
Therefore, we get
\begin{eqnarray*}
& &\hspace{-20pt}(\varphi^G_{\mu_{r-1}\mu_r(\textbf{V})})_{r-1}
= \frac{(\varphi^G_{\textbf{V}})_{r-2} (\varphi^G_{\mu_r(\textbf{V})})^2_{r}+(\varphi^G_{\textbf{V}})_{2r-2}(\varphi^G_{\textbf{V}})_{2r-1} (\varphi^G_{\textbf{V}})_{-(r-1)}}{(\varphi^G_{\textbf{V}})_{r-1}}
\\  
&=& \frac{(a^{\Lambda_{r-2}-\Lambda_{r-1}}\frac{Y_{2,r-1}}{Y^2_{2,r}}(\varphi^G_{\textbf{V}})_{r-1}-a^{\Lambda_{r-2}}\frac{Y_{1,r-1}Y_{2,r-1}}{Y^2_{2,r}})
a^{2\Lambda_{r-1}}Y^2_{2,r}
+a^{\Lambda_{r-2}+2\Lambda_{r-1}}Y_{1,r-1}Y_{2,r-1}}
{(\varphi^G_{\textbf{V}})_{r-1}}
\\
&=& a^{\Lambda_{r-2}+\Lambda_{r-1}}Y_{2,r-1},
\end{eqnarray*}
where we use Proposition \ref{iniB} (ii) and (\ref{twovarsB}) in the second equality.

Next, we assume that $r-k>1$. It follows from Lemma \ref{mutgamlem} that the arrows incident to the vertex $(\varphi_{\textbf{V}})_{k}$ in the mutation diagram of $\mu_{k+1}\cdots\mu_{r-1}\mu_r(\Sigma_0)$ are as follows:
\[
\begin{xy}
(-60,98)*{\ \ \ (\varphi_{\textbf{V}})_{r+k-1}}="r+k",
(-60,90) *{\ (\varphi_{\textbf{V}})_{k-1}}="k",
(-35,90) *{\ (\varphi_{\textbf{V}})_{k}}="4",
(-35,82)*{(\varphi_{\textbf{V}})_{-k}}="-2",
(20,98)*{(\varphi_{\textbf{V}})_{2r-1}}="6",
(0,90) *{(\varphi_{\mu_{k+1}\cdots\mu_{r-1}\mu_r(\textbf{V})})_{k+1}}="2",
\ar@{->} "4";"r+k"
\ar@{->} "4";"-2"
\ar@{->} "k";"4"
\ar@{->} "2";"4"
\ar@{->} "4";"6"
\end{xy}
\]
Thus, we obtain
\begin{eqnarray*}
& &\hspace{-20pt}(\varphi^G_{\mu_k\mu_{k+1}\cdots\mu_{r-1}\mu_r(\textbf{V})})_{k}= \frac{(\varphi^G_{\textbf{V}})_{k-1} (\varphi^G_{\mu_{k+1}\cdots\mu_{r-1}\mu_r(\textbf{V})})_{k+1}+(\varphi^G_{\textbf{V}})_{r+k-1}(\varphi^G_{\textbf{V}})_{2r-1} (\varphi^G_{\textbf{V}})_{-k}}{(\varphi^G_{\textbf{V}})_{k}} \\  
&=& \frac{(a^{\Lambda_{k-1}-\Lambda_k}\frac{Y_{2,k}}{Y_{2,k+1}}(\varphi^G_{\textbf{V}})_{k}-a^{\Lambda_{k-1}}\frac{Y_{1,k}Y_{2,k}}{Y_{2,k+1}})
a^{\Lambda_k+\Lambda_{r-1}}Y_{2,k+1}
+
a^{\Lambda_k+\Lambda_{k-1}+\Lambda_{r-1}}Y_{1,k}Y_{2,k}}{(\varphi^G_{\textbf{V}})_{k}} \\
&=& a^{\Lambda_{k-1}+\Lambda_{r-1}}Y_{2,k},
\end{eqnarray*}
where we use Proposition \ref{iniB} (ii), (\ref{twovarsB}) and the induction hypothesis in the second equality. \qed

\vspace{3mm}

\nd
[{\sl Proof of Theorem \ref{thmB} (ii)}]

For fixed $k$, we use the induction on $(l-k)$ to prove (\ref{B-exp1}). If $l-k=0$, so that $k=l$, then the mutation diagram (\ref{inidiaB-1}) means that
\begin{eqnarray*}
(\varphi^G_{\mu_k(\textbf{V})})_k
&=&\frac{(\varphi^G_{\textbf{V}})_{k-1}(\varphi^G_{\textbf{V}})_{r+k}(\varphi^G_{\textbf{V}})_{-(k+1)}
+(\varphi^G_{\textbf{V}})_{k+1}(\varphi^G_{\textbf{V}})_{r+k-1}(\varphi^G_{\textbf{V}})_{-k}}{(\varphi^G_{\textbf{V}})_k}\\
&=&\frac{(\varphi^G_{\textbf{V}})_{k-1}a^{\Lm_{k}+\Lm_{k+1}}Y_{1,k+1}Y_{2,k+1}
+(\varphi^G_{\textbf{V}})_{k+1}a^{\Lm_{k-1}+\Lm_{k}}Y_{1,k}Y_{2,k}}{(\varphi^G_{\textbf{V}})_k}\\
&=&\frac{(\varphi^G_{\textbf{V}})_{k}a^{\Lm_{k-1}+\Lm_{k+1}}Y_{2,k}Y_{1,k+1}
+(\varphi^G_{\textbf{V}})_{k}a^{\Lm_{k-1}+\Lm_{k+1}}\frac{Y_{1,k}Y_{2,k}Y_{2,k+2}}{Y_{2,k+1}}}{(\varphi^G_{\textbf{V}})_k}\\
&=&a^{\Lm_{k-1}+\Lm_{k+1}}\left(Y_{2,k}Y_{1,k+1}+\frac{Y_{1,k}Y_{2,k}Y_{2,k+2}}{Y_{2,k+1}}\right),
\end{eqnarray*}
where we use Proposition \ref{iniB} (ii) in the second equality, and (\ref{twovarsB}) for $k$ and $k+1$ in the third equality.

Next, we assume that $l-k>0$. The vertices and arrows around the vertex $(\varphi_{\textbf{V}})_{l}$ in the mutation diagram of $\mu_{l-1}\cdots\mu_k(\Sigma_0)$ $(l\leq r-2)$ are as follows (Lemma \ref{mutgamlem}):
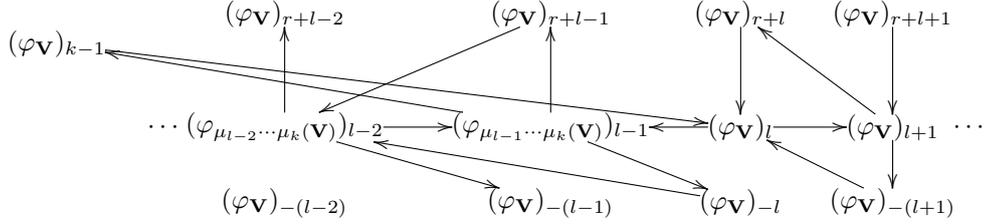
\[
\begin{xy}
(90,90)*{\cdots}="emp1",
(80,105)*{(\varphi_{\textbf{V}})_{r+l+1}}="7",
(80,90) *{(\varphi_{\textbf{V}})_{l+1}}="3",
(80,80)*{(\varphi_{\textbf{V}})_{-(l+1)}}="-4",
(60,105)*{(\varphi_{\textbf{V}})_{r+l}}="r+k",
(60,90) *{(\varphi_{\textbf{V}})_{l}}="k",
(60,80)*{(\varphi_{\textbf{V}})_{-l}}="-j_k",
(35,105)*{(\varphi_{\textbf{V}})_{r+l-1}}="8",
(35,90) *{(\varphi_{\mu_{l-1}\cdots\mu_k(\textbf{V})})_{l-1}}="4",
(35,80)*{(\varphi_{\textbf{V}})_{-(l-1)}}="-2",
(0,105)*{(\varphi_{\textbf{V}})_{r+l-2}}="6",
(0,90) *{(\varphi_{\mu_{l-2}\cdots\mu_k(\textbf{V})})_{l-2}}="2",
(0,80)*{(\varphi_{\textbf{V}})_{-(l-2)}}="-1",
(-16,90)*{\cdots}="emp",
(-30,101)*{(\varphi_{\textbf{V}})_{k-1}}="k-1",
\ar@{->} "7";"3"
\ar@{->} "r+k";"k"
\ar@{->} "4";"8"
\ar@{->} "2";"6"
\ar@{->} "3";"-4"
\ar@{->} "k";"3"
\ar@{->} "k";"4"
\ar@{->} "4";"-j_k"
\ar@{->} "2";"4"
\ar@{->} "3";"r+k"
\ar@{->} "8";"2"
\ar@{->} "-4";"k"
\ar@{->} "k-1";"k"
\ar@{->} "4";"k-1"
\ar@{->} "-j_k";"2"
\ar@{->} "2";"-2"
\end{xy}
\]
It follows from this diagram that
\begin{equation}\label{B-lem2pr-1}
(\varphi^G_{\mu_l\mu_{l-1}\cdots\mu_k(\textbf{V})})_l
=\frac{(\varphi^G_{\textbf{V}})_{l+1}(\varphi^G_{\mu_{l-1}\cdots\mu_k(\textbf{V})})_{l-1}
+(\varphi^G_{\textbf{V}})_{k-1}(\varphi^G_{\textbf{V}})_{r+l}(\varphi^G_{\textbf{V}})_{-(l+1)}}{(\varphi^G_{\textbf{V}})_l}.
\end{equation}
Using the induction hypothesis, the cluster variable $(\varphi^G_{\mu_{l-1}\cdots\mu_k(\textbf{V})})_{l-1}$ is
\[
a^{\Lm_{k-1}+\Lm_{l}} Y_{2,k}Y_{1,l} (1+A^{-1}_{1,l}+A^{-1}_{1,l}A^{-1}_{1,l-1}+\cdots+A^{-1}_{1,l}A^{-1}_{1,l-1}\cdots A^{-1}_{1,k+1}).
\]
In conjunction with (\ref{iniBexp2}), it follows that $(\varphi^G_{\mu_{l-1}\cdots\mu_k(\textbf{V})})_{l-1}$ coincides with
\begin{multline}
a^{\Lm_{k-1}+\Lm_{l}} Y_{2,k}Y_{1,l} (a^{-\Lm_l}Y^{-1}_{1,l}
(\varphi^G_{\textbf{V}})_{l}
-a^{-\Lm_{k-1}}Y^{-1}_{1,k-1} (\varphi^G_{\textbf{V}})_{k-1}
A^{-1}_{1,l}A^{-1}_{1,l-1}\cdots A^{-1}_{1,k}) \\
=a^{\Lm_{k-1}}Y_{2,k}(\varphi^G_{\textbf{V}})_{l}
-a^{\Lm_{l}}Y_{2,l+1} (\varphi^G_{\textbf{V}})_{k-1}. \label{B-lem2pr-11}
\end{multline}

Applying Proposition \ref{iniB} (ii), (\ref{twovarsB}) and this formula to (\ref{B-lem2pr-1}), we see that
\begin{eqnarray*}
& &\hspace{-20pt}(\varphi^G_{\mu_l\mu_{l-1}\cdots\mu_k(\textbf{V})})_l \\
&=& \frac{(a^{\Lm_{l+1}}Y_{1,l+1}+a^{\Lm_{l+1}-\Lm_l}\frac{Y_{2,l+2}}{Y_{2,l+1}}(\varphi^G_{\textbf{V}})_{l})(a^{\Lm_{k-1}}Y_{2,k}(\varphi^G_{\textbf{V}})_{l}
-a^{\Lm_{l}}Y_{2,l+1} (\varphi^G_{\textbf{V}})_{k-1})}{(\varphi^G_{\textbf{V}})_l}\\
& &+\frac{
(\varphi^G_{\textbf{V}})_{k-1}a^{\Lm_{l}+\Lm_{l+1}}Y_{1,l+1}Y_{2,l+1}}{(\varphi^G_{\textbf{V}})_l}\\
&=& a^{\Lm_{l+1}-\Lm_{l}+\Lm_{k-1}}\frac{Y_{2,l+2}Y_{2,k}}{Y_{2,l+1}}(\varphi^G_{\textbf{V}})_l+
a^{\Lm_{l+1}+\Lm_{k-1}}Y_{1,l+1}Y_{2,k}-a^{\Lm_{l+1}}Y_{2,l+2}(\varphi^G_{\textbf{V}})_{k-1}\\
&=& a^{\Lm_{l+1}+\Lm_{k-1}}Y_{1,l+1}Y_{2,k}(1+A^{-1}_{1,l+1}+A^{-1}_{1,l+1}A^{-1}_{1,l}+
\cdots+A^{-1}_{1,l+1}A^{-1}_{1,l}\cdots A^{-1}_{1,1})\\
& &-a^{\Lm_{l+1}+\Lm_{k-1}}Y_{1,k-1}Y_{2,l+2}(1+A^{-1}_{1,k-1}+A^{-1}_{1,k-1}A^{-1}_{1,k-2}+
\cdots+A^{-1}_{1,k-1}A^{-1}_{1,k-2}\cdots A^{-1}_{1,1})\\
&=& a^{\Lm_{l+1}+\Lm_{k-1}}Y_{1,l+1}Y_{2,k}(1+A^{-1}_{1,l+1}+A^{-1}_{1,l+1}A^{-1}_{1,l}+
\cdots+A^{-1}_{1,l+1}A^{-1}_{1,l}\cdots A^{-1}_{1,k+1}),
\end{eqnarray*}
where we use (\ref{iniBexp2}) in the third equality, and $Y_{1,l+1}Y_{2,k}A^{-1}_{1,l+1}A^{-1}_{1,l}\cdots A^{-1}_{1,k}=Y_{1,k-1}Y_{2,l+2}$ in the fourth equality. Hence, we get (\ref{B-exp1}) for all $l$ $(k\leq l\leq r-2)$.

By the same argument as in the proof of Proposition \ref{iniB}, the equation (\ref{B-cry1}) follows from (\ref{B-exp1}). \qed

\vspace{3mm}

\nd
[{\sl Proof of Theorem \ref{thmB} (iii)}]

Using Lemma \ref{mutgamlem} repeatedly, we see that the mutation diagram of $\mu_{r-2}\cdots\mu_{k+1}\mu_k(\Sigma_0)$ is as follows:
\begin{equation}\label{B-dia-1}
\begin{xy}
(80,108)*{(\varphi_{\textbf{V}})_{2r}}="7",
(80,90) *{(\varphi_{\textbf{V}})_{r}}="3",
(80,77)*{(\varphi_{\textbf{V}})_{-r}}="-4",
(60,108)*{(\varphi_{\textbf{V}})_{2r-1}}="r+k",
(60,90) *{(\varphi_{\textbf{V}})_{r-1}}="k",
(60,77)*{(\varphi_{\textbf{V}})_{-(r-1)}}="-j_k",
(35,108)*{(\varphi_{\textbf{V}})_{2r-2}}="8",
(35,90) *{(\varphi_{\mu_{r-2}\cdots\mu_k(\textbf{V})})_{r-2}}="4",
(35,77)*{(\varphi_{\textbf{V}})_{-(r-2)}}="-2",
(0,108)*{(\varphi_{\textbf{V}})_{2r-3}}="6",
(0,90) *{(\varphi_{\mu_{r-3}\cdots\mu_k(\textbf{V})})_{r-3}}="2",
(0,77)*{(\varphi_{\textbf{V}})_{-(r-3)}}="-1",
(-16,90)*{\cdots}="emp",
(-30,102)*{(\varphi_{\textbf{V}})_{k-1}}="k-1",
\ar@{->} "7";"3"
\ar@{->} "r+k";"k"
\ar@{->} "4";"8"
\ar@{->} "2";"6"
\ar@{->} "3";"-4"
\ar@{->} "k";"3"^{-2}
\ar@{->} "k";"4"
\ar@{->} "4";"-j_k"
\ar@{->} "2";"4"
\ar@{->} "3";"r+k"_2
\ar@{->} "8";"2"
\ar@{->} "-4";"k"^{2}
\ar@{->} "k-1";"k"
\ar@{->} "4";"k-1"
\ar@{->} "-j_k";"2"
\ar@{->} "2";"-2"
\end{xy},
\end{equation}
which implies that
\begin{equation}\label{B-lem2pr-2}
(\varphi^G_{\mu_{r-1}\mu_{r-2}\cdots\mu_k(\textbf{V})})_{r-1}
=\frac{(\varphi^G_{\textbf{V}})^2_{r}(\varphi^G_{\mu_{r-2}\cdots\mu_k(\textbf{V})})_{r-2}
+(\varphi^G_{\textbf{V}})_{k-1}(\varphi^G_{\textbf{V}})_{2r-1}(\varphi^G_{\textbf{V}})^2_{-r}}{(\varphi^G_{\textbf{V}})_{r-1}}.
\end{equation}
On the other hand, from (\ref{iniBexp2}) and (\ref{iniBexp-r}), we have
\begin{equation}\label{B-lem2pr-3}
(\varphi^G_{\textbf{V}})_{r}=a^{\Lm_r}\left(Y_{1,r}+\frac{a^{-\Lm_{r-1}}}{Y_{2,r}}(\varphi^G_{\textbf{V}})_{r-1}\right).
\end{equation}
Just as in (\ref{B-lem2pr-11}), we see that
\begin{equation}\label{B-lem2pr-4}
(\varphi^G_{\mu_{r-2}\cdots\mu_k(\textbf{V})})_{r-2}=
a^{\Lm_{k-1}}Y_{2,k}(\varphi^G_{\textbf{V}})_{r-1}-
a^{\Lm_{r-1}}Y^2_{2,r}(\varphi^G_{\textbf{V}})_{k-1}.
\end{equation}
Substituting (\ref{B-lem2pr-3}) and (\ref{B-lem2pr-4}) for (\ref{B-lem2pr-2}), it follows
\begin{eqnarray}
& &\hspace{-20pt}(\varphi^G_{\mu_{r-1}\mu_{r-2}\cdots\mu_k(\textbf{V})})_{r-1} \nonumber \\
&=& \frac{a^{2\Lm_{r}}(Y^2_{1,r}+2\frac{a^{-\Lm_{r-1}}Y_{1,r}}{Y_{2,r}}(\varphi^G_{\textbf{V}})_{r-1}+\frac{a^{-2\Lm_{r-1}}}{Y^2_{2,r}}(\varphi^G_{\textbf{V}})^2_{r-1})
(a^{\Lm_{k-1}}Y_{2,k}(\varphi^G_{\textbf{V}})_{r-1}-
a^{\Lm_{r-1}}Y^2_{2,r}(\varphi^G_{\textbf{V}})_{k-1})}{(\varphi^G_{\textbf{V}})_{r-1}} \nonumber\\
& &+\frac{a^{2\Lm_{r}+\Lm_{r-1}}Y^2_{1,r}Y^2_{2,r}(\varphi^G_{\textbf{V}})_{k-1}
}{(\varphi^G_{\textbf{V}})_{r-1}}\nonumber\\
&=& a^{2\Lm_{r}+\Lm_{k-1}} Y^2_{1,r} Y_{2,k}+2a^{2\Lm_{r}-\Lm_{r-1}+\Lm_{k-1}}\frac{Y_{2,k}Y_{1,r}}{Y_{2,r}}(\varphi^G_{\textbf{V}})_{r-1}
\label{B-lem2pr-5}
\\
& &+a^{2\Lm_{r}-2\Lm_{r-1}+\Lm_{k-1}}\frac{Y_{2,k}}{Y^2_{2,r}}(\varphi^G_{\textbf{V}})^2_{r-1}-2a^{2\Lm_{r}}Y_{1,r}Y_{2,r}(\varphi^G_{\textbf{V}})_{k-1}-a^{2\Lm_{r}-\Lm_{r-1}}(\varphi^G_{\textbf{V}})_{r-1}(\varphi^G_{\textbf{V}})_{k-1}\nonumber \\
&=&
a^{2\Lm_{r}+\Lm_{k-1}} \left(Y^2_{1,r} Y_{2,k}+ 2\frac{Y_{2,k}Y_{1,r}Y_{1,r-1}}{Y_{2,r}}(1+A^{-1}_{1,r-1}+\cdots+A^{-1}_{1,r-1}\cdots A^{-1}_{1,k+1})\right. \label{B-lem2pr-5-1} \\
& &\left.+\frac{Y_{2,k}Y^2_{1,r-1}}{Y^2_{2,r}}
(1+A^{-1}_{1,r-1}+\cdots+A^{-1}_{1,r-1}\cdots A^{-1}_{1,1})
(1+A^{-1}_{1,r-1}+\cdots+A^{-1}_{1,r-1}\cdots A^{-1}_{1,k+1})\right)\nonumber \\
&=&
a^{2\Lm_{r}+\Lm_{k-1}} \left(Y^2_{1,r} Y_{2,k}+ 2Y^2_{1,r}Y_{2,k}A^{-1}_{1,r}(1+A^{-1}_{1,r-1}+\cdots+A^{-1}_{1,r-1}\cdots A^{-1}_{1,k+1})\right.\nonumber\\
& &\quad +Y^2_{1,r}Y_{2,k}A^{-2}_{1,r}
(1+A^{-1}_{1,r-1}+\cdots+A^{-1}_{1,r-1}\cdots A^{-1}_{1,k+1})^2\nonumber\\
& &+\frac{Y_{2,k}Y^2_{1,r-1}}{Y^2_{2,r}}
A^{-1}_{1,r-1}\cdots A^{-1}_{1,k+1}A^{-1}_{1,k}
(1+A^{-1}_{1,r-1}+\cdots+A^{-1}_{1,r-1}\cdots A^{-1}_{1,k+1})\nonumber\\
& &\qquad \qquad \qquad \qquad \qquad \qquad\left.\times (1+A^{-1}_{1,k-1}+A^{-1}_{1,k-1}A^{-1}_{1,k-2}+\cdots+A^{-1}_{1,k-1}\cdots A^{-1}_{1,1})\right),\nonumber
\end{eqnarray}
where we use (\ref{iniBexp2}) in the third equality. Note that by Theorem \ref{kashidem}, the monomial realization $\mu$ of the Demazure crystal $B(2\Lm_r+\Lm_k)_{s_{k+1}\cdots s_{r-1}s_{r}}$ such that the highest weight vector is mapped to $Y:=Y^2_{1,r}Y_{2,k}$ is as follows:
\[
\begin{xy}
(-75,15)*{Y}="h",
(-60,15)*{A^{-1}_{1,r}Y}="hr",
(-35,15)*{A^{-1}_{1,r-1}A^{-1}_{1,r}Y}="hrr-1",
(-10,15)*{\ \cdots\ }="dot1",
(20,15)*{\ A^{-1}_{1,k+1}\cdots A^{-1}_{1,r-1}A^{-1}_{1,r}Y}="hr1",
(-60,0)*{A^{-2}_{1,r}Y}="hrr",
(-35,0)*{A^{-1}_{1,r-1}A^{-2}_{1,r}Y}="hrrr-1",
(-10,0)*{\ \cdots\ }="dot2",
(20,0)*{\ A^{-1}_{1,k+1}\cdots A^{-1}_{1,r-1}A^{-2}_{1,r}Y}="hrr1",
(-35,-15)*{A^{-2}_{1,r-1}A^{-2}_{1,r}Y}="hrrr-1r-1",
(-10,-15)*{\ \ \ \cdots\ }="dot3",
(30,-15)*{\ A^{-1}_{1,k+1}\cdots A^{-1}_{1,r-2}A^{-2}_{1,r-1}A^{-2}_{1,r}Y}="hrrr-11",
(-10,-25)*{\ \ \ \ \ \cdots\ }="dot4",
(30,-25)*{\cdots\ \ \ }="dot5",
(30,-35)*{\ A^{-2}_{1,k+1}\cdots A^{-2}_{1,r-2}A^{-2}_{1,r-1}A^{-2}_{1,r}Y}="l",
\ar@{->} "h";"hr"^{r\ \ }
\ar@{->} "hr";"hrr"^{r}
\ar@{->} "hr";"hrr-1"^{r-1\ \ }
\ar@{->} "hrr-1";"dot1"^{\ \ r-2}
\ar@{->} "dot1";"hr1"^{k+1\ \ \qquad \ \ \ \ }
\ar@{->} "hrr";"hrrr-1"^{r-1\ \ }
\ar@{->} "hrrr-1";"dot2"^{\ \ r-2}
\ar@{->} "dot2";"hrr1"^{k+1\ \ \qquad \ \ \ \ }
\ar@{->} "hrrr-1";"hrrr-1r-1"^{r-1\ \ }
\ar@{->} "hrrr-1r-1";"dot3"^{\ \ r-2}
\ar@{->} "dot3";"hrrr-11"^{k+1\ \ \qquad \ \ \ \ \ \ \ }
\ar@{->} "dot3";"dot4"^{r-2}
\ar@{->} "dot5";"l"^{k+1}
\end{xy}
\]
Thus, we see that
\begin{eqnarray}
& &\hspace{-20pt}a^{2\Lm_{r}+\Lm_{k-1}} (Y^2_{1,r} Y_{2,k}+ 2Y^2_{1,r}Y_{2,k}A^{-1}_{1,r}(1+A^{-1}_{1,r-1}+A^{-1}_{1,r-1}A^{-1}_{1,r-2}+\cdots+A^{-1}_{1,r-1}\cdots A^{-1}_{1,k+1}))\nonumber\\
& &+a^{2\Lm_{r}+\Lm_{k-1}}Y^2_{1,r}Y_{2,k}A^{-2}_{1,r}
(1+A^{-1}_{1,r-1}+\cdots+A^{-1}_{1,r-1}\cdots A^{-1}_{1,k+1})^2\nonumber\\
&=&a^{2\Lm_{r}+\Lm_{k-1}} \sum_{b\in B(\Lambda_k+2\Lambda_{r})_{s_{k+1}s_{k+2}\cdots s_{r}}} C(b)\mu(b),\label{B-lem2pr-c}
\end{eqnarray}
where each coefficient $C(b)$ is either $1$ or $2$. Similarly, we also see that
\begin{eqnarray*}
& &\hspace{-20pt}a^{2\Lm_{r}+\Lm_{k-1}}\frac{Y_{2,k}Y^2_{1,r-1}}{Y^2_{2,r}}
A^{-1}_{1,r-1}\cdots A^{-1}_{1,k+1}A^{-1}_{1,k}
(1+A^{-1}_{1,r-1}+\cdots+A^{-1}_{1,r-1}\cdots A^{-1}_{1,k+1})\\
& &\qquad \qquad \qquad \qquad \qquad \times (1+A^{-1}_{1,k-1}+A^{-1}_{1,k-1}A^{-1}_{1,k-2}+\cdots+A^{-1}_{1,k-1}\cdots A^{-1}_{1,1}),
\\
&=&a^{2\Lm_{r}+\Lm_{k-1}}Y_{1,r-1}Y_{1,k-1}
(1+A^{-1}_{1,r-1}+\cdots+A^{-1}_{1,r-1}\cdots A^{-1}_{1,k+1})\\
& &\qquad \qquad \qquad \qquad \qquad \times (1+A^{-1}_{1,k-1}+A^{-1}_{1,k-1}A^{-1}_{1,k-2}+\cdots+A^{-1}_{1,k-1}\cdots A^{-1}_{1,1}),
\\
&=&a^{2\Lm_{r}+\Lm_{k-1}} \sum_{b\in B(\Lambda_{k-1}+\Lambda_{r-1})_{s_1\cdots s_{k-1}s_{k+1}s_{k+2}\cdots s_{r-1}}} \mu'(b),
\end{eqnarray*}
where $\mu'$ is the monomial realization of the crystal base $B(\Lambda_{k-1}+\Lambda_{r-1})$ in our claim. Hence, we obtain (\ref{B-cry2}).

Finally, we prove (\ref{B-cry3}). Applying the mutation $\mu_{r-1}$ to the diagram (\ref{B-dia-1}), we obtain the mutation diagram of $\mu_{r-1}\mu_{r-2}\cdots\mu_{k+1}\mu_k(\Sigma_0)$:
\begin{equation}\label{B-dia-2}
\begin{xy}
(93,108)*{(\varphi_{\textbf{V}})_{2r}}="7",
(93,90) *{(\varphi_{\textbf{V}})_{r}}="3",
(93,77)*{(\varphi_{\textbf{V}})_{-r}}="-4",
(68,108)*{(\varphi_{\textbf{V}})_{2r-1}}="r+k",
(68,90) *{(\varphi_{\mu_{r-1}\cdots\mu_k(\textbf{V})})_{r-1}}="k",
(68,77)*{(\varphi_{\textbf{V}})_{-(r-1)}}="-j_k",
(35,108)*{(\varphi_{\textbf{V}})_{2r-2}}="8",
(35,90) *{(\varphi_{\mu_{r-2}\cdots\mu_k(\textbf{V})})_{r-2}}="4",
(35,77)*{(\varphi_{\textbf{V}})_{-(r-2)}}="-2",
(0,108)*{(\varphi_{\textbf{V}})_{2r-3}}="6",
(0,90) *{(\varphi_{\mu_{r-3}\cdots\mu_k(\textbf{V})})_{r-3}}="2",
(0,77)*{(\varphi_{\textbf{V}})_{-(r-3)}}="-1",
(-16,90)*{\cdots}="emp",
(-30,104)*{(\varphi_{\textbf{V}})_{k-1}}="k-1",
\ar@{->} "7";"3"
\ar@{->} "k";"r+k"
\ar@{->} "4";"8"
\ar@{->} "2";"6"
\ar@{->} "-4";"3"
\ar@{->} "-4";"4"_2
\ar@{->} "3";"k"^{\qquad \quad 2}
\ar@{->} "4";"k"
\ar@{->} "4";"-j_k"
\ar@{->} "2";"4"
\ar@{->} "r+k";"4"
\ar@{->} "8";"2"
\ar@{->} "k";"-4"^{-2}
\ar@{->} "k";"k-1"
\ar@{->} "k-1";"3"^{-2}
\ar@{->} "-j_k";"2"
\ar@{->} "2";"-2"
\end{xy}
\end{equation}
This diagram says that
\begin{equation}\label{B-lem2pr-6}
(\varphi^G_{\mu_{r}\mu_{r-1}\cdots\mu_{k+1}\mu_k(\textbf{V})})_r
=\frac{1}{(\varphi^G_{\textbf{V}})_r}\left((\varphi^G_{\mu_{r-1}\cdots\mu_{k+1}\mu_k(\textbf{V})})_{r-1}+
(\varphi^G_{\textbf{V}})_{k-1}(\varphi^G_{\textbf{V}})_{-r}(\varphi^G_{\textbf{V}})_{2r}\right).
\end{equation}
Using (\ref{B-lem2pr-3}) and (\ref{B-lem2pr-5}), the following holds:
\[
(\varphi^G_{\mu_{r-1}\cdots\mu_{k+1}\mu_k(\textbf{V})})_{r-1}=a^{\Lm_{k-1}}(\varphi^G_{\textbf{V}})^2_r Y_{2,k}-a^{2\Lm_{r}}Y_{1,r}Y_{2,r}(\varphi^G_{\textbf{V}})_{k-1}
-a^{\Lm_{r}}Y_{2,r}(\varphi^G_{\textbf{V}})_{r}(\varphi^G_{\textbf{V}})_{k-1}.
\]
Applying this to (\ref{B-lem2pr-6}), it is easy to see that
\begin{eqnarray}
& &\hspace{-20pt}(\varphi^G_{\mu_{r}\mu_{r-1}\cdots\mu_{k+1}\mu_k(\textbf{V})})_r\nonumber\\
&=&\frac{1}{(\varphi^G_{\textbf{V}})_r}\left(a^{\Lm_{k-1}}Y_{2,k}(\varphi^G_{\textbf{V}})^2_r-a^{\Lm_{r}}Y_{2,r}(\varphi^G_{\textbf{V}})_r(\varphi^G_{\textbf{V}})_{k-1}\right)\nonumber\\
&=&a^{\Lm_{k-1}}Y_{2,k}(\varphi^G_{\textbf{V}})_r-a^{\Lm_{r}}Y_{2,r}(\varphi^G_{\textbf{V}})_{k-1}\nonumber\\
&=&a^{\Lm_{k-1}+\Lm_r}Y_{2,k}Y_{1,r}(1+A^{-1}_{1,r}+A^{-1}_{1,r}A^{-1}_{1,r-1}+\cdots+A^{-1}_{1,r}A^{-1}_{1,r-1}\cdots A^{-1}_{1,k+1}),\qquad \qquad \nonumber
\end{eqnarray}
which means (\ref{B-cry3}). \qed

\vspace{3mm}

\nd
[{\sl Proof of Theorem \ref{thmB} (iv).}]

We put $A:=1+A^{-1}_{1,r-1}+A^{-1}_{1,r-1}A^{-1}_{1,r-2}+\cdots+A^{-1}_{1,r-1}A^{-1}_{1,r-2}\cdots A^{-1}_{1,k+1}$.
It now follows at once from (\ref{B-lem2pr-5-1}) and $(\ref{B-lem2pr-7})$ that
\begin{eqnarray}
& &\hspace{-20pt}
(\varphi^G_{\mu_{r-1}\cdots\mu_{k+1}\mu_k(\textbf{V})})_{r-1}\nonumber\\
&=&
a^{2\Lm_r+\Lm_{k-1}}Y_{2,k}(Y^2_{1,r}+2\frac{Y_{1,r-1}Y_{1,r}}{Y_{2,r}}A
+\frac{Y^2_{1,r-1}}{Y^2_{2,r}}A(1+A^{-1}_{1,r-1}+\cdots+A^{-1}_{1,r-1}\cdots A^{-1}_{1,1}))\ \ \ \qquad \label{B-lem3-pr1}\\
&=&a^{2\Lm_r+\Lm_{k-1}}Y_{2,k}(Y^2_{1,r}+2\frac{Y_{1,r-1}Y_{1,r}}{Y_{2,r}}A)\nonumber\\
&+&a^{2\Lm_r+\Lm_{k-1}}\frac{Y_{2,k}Y^2_{1,r-1}}{Y^2_{2,r}}(A^2+A\cdot A^{-1}_{1,r-1}\cdots A^{-1}_{1,k}(1+A^{-1}_{1,k-1}+\cdots+A^{-1}_{1,k-1}A^{-1}_{1,k-2}\cdots A^{-1}_{1,1}))\nonumber\\
&=&a^{2\Lm_r+\Lm_{k-1}}Y_{2,k}(Y^2_{1,r}+2\frac{Y_{1,r-1}Y_{1,r}}{Y_{2,r}}A+\frac{Y^2_{1,r-1}}{Y^2_{2,r}}A^2+a^{-\Lm_{k-1}}\frac{Y_{1,r-1}}{Y_{2,k}}A(\varphi^G_{\textbf{V}})_{k-1}),\label{B-lem3-pr2}
\end{eqnarray}
\begin{equation}\label{B-lem3-pr2-1}
(\varphi^G_{\mu_r\mu_{r-1}\cdots\mu_{k+1}\mu_k(\textbf{V})})_{r}
=a^{\Lm_{k-1}+\Lm_r}Y_{2,k}Y_{1,r}(1+A\times A^{-1}_{1,r}).
\end{equation}
To prove (\ref{B-lem3-exp}), we use the induction on $(k-j)$. First, we consider the case $k-j=1$. Applying the mutation $\mu_r$ to (\ref{B-dia-2}), one can verify that the vertices and arrows around $(\varphi_{\textbf{V}})_{k-1}$ in the mutation diagram of $\mu_{r}\mu_{r-1}\cdots\mu_{k+1}\mu_k(\Sigma_0)$ are as follows (Lemma \ref{mutgamlem}):
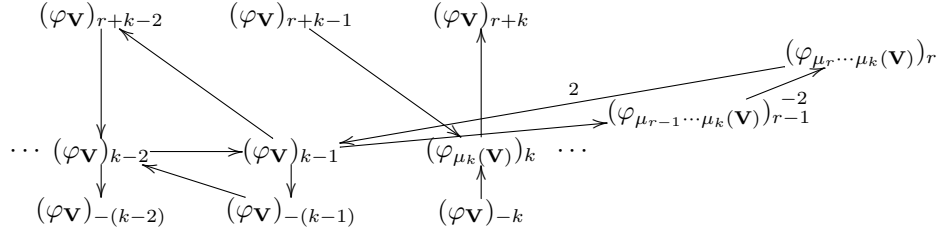
\[
\begin{xy}
(110,103) *{(\varphi_{\mu_{r}\cdots\mu_k(\textbf{V})})_{r}}="r",
(90,95) *{(\varphi_{\mu_{r-1}\cdots\mu_k(\textbf{V})})_{r-1}}="r-1",
(72,90)*{\cdots}="emp1",
(60,108)*{(\varphi_{\textbf{V}})_{r+k}}="r+k",
(60,90) *{(\varphi_{\mu_k(\textbf{V})})_{k}}="k",
(60,82)*{(\varphi_{\textbf{V}})_{-k}}="-j_k",
(35,108)*{(\varphi_{\textbf{V}})_{r+k-1}}="8",
(35,90) *{(\varphi_{\textbf{V}})_{k-1}}="4",
(35,82)*{(\varphi_{\textbf{V}})_{-(k-1)}}="-2",
(10,108)*{(\varphi_{\textbf{V}})_{r+k-2}}="6",
(10,90) *{(\varphi_{\textbf{V}})_{k-2}}="2",
(10,82)*{(\varphi_{\textbf{V}})_{-(k-2)}}="-1",
(0,90)*{\cdots}="emp",
\ar@{->} "r";"4"_2
\ar@{->} "4";"r-1"
\ar@{->} "k";"r+k"
\ar@{->} "6";"2"
\ar@{->} "-j_k";"k"
\ar@{->} "8";"k"
\ar@{->} "4";"-2"
\ar@{->} "2";"-1"
\ar@{->} "2";"4"
\ar@{->} "4";"6"
\ar@{->} "-2";"2"
\ar@{->} "r-1";"r"_{-2}
\end{xy}
\]
Therefore, the cluster variable $(\varphi^G_{\mu_{k-1}\mu_{r}\mu_{r-1}\cdots\mu_{k+1}\mu_k(\textbf{V})})_{k-1}$ is
\begin{equation}\label{B-lem3-pr3}
\frac{(\varphi^G_{\mu_{r}\mu_{r-1}\cdots\mu_{k+1}\mu_k(\textbf{V})})^2_{r}(\varphi^G_{\textbf{V}})_{k-2}+(\varphi^G_{\mu_{r-1}\cdots\mu_{k+1}\mu_k(\textbf{V})})_{r-1}(\varphi^G_{\textbf{V}})_{r+k-2}(\varphi^G_{\textbf{V}})_{-(k-1)}}{(\varphi^G_{\textbf{V}})_{k-1}}.
\end{equation}

Substituting (\ref{B-lem3-pr2}), (\ref{B-lem3-pr2-1}) and $(\varphi^G_{\textbf{V}})_{k-2}=a^{\Lm_{k-2}-\Lm_{k-1}}\frac{Y_{2,k-1}}{Y_{2,k}}(\varphi^G_{\textbf{V}})_{k-1}-a^{\Lm_{k-2}}\frac{Y_{1,k-1}Y_{2,k-1}}{Y_{2,k}}$ (\ref{twovarsB}) for (\ref{B-lem3-pr3}), we get
\begin{eqnarray*}
& &\hspace{-20pt}(\varphi^G_{\mu_{k-1}\mu_{r}\mu_{r-1}\cdots\mu_{k+1}\mu_k(\textbf{V})})_{k-1}\\
&=&\frac{a^{2\Lm_{k-1}+2\Lm_r}Y^2_{2,k}Y^2_{1,r}(1+A\times A^{-1}_{1,r})^2
(a^{\Lm_{k-2}-\Lm_{k-1}}\frac{Y_{2,k-1}}{Y_{2,k}}(\varphi^G_{\textbf{V}})_{k-1}-a^{\Lm_{k-2}}\frac{Y_{1,k-1}Y_{2,k-1}}{Y_{2,k}})}{(\varphi^G_{\textbf{V}})_{k-1}}\\
&+&\frac{a^{2\Lm_r+2\Lm_{k-1}+\Lm_{k-2}}Y_{1,k-1}Y_{2,k-1}Y_{2,k}(Y^2_{1,r}+2\frac{Y_{1,r-1}Y_{1,r}}{Y_{2,r}}A+\frac{Y^2_{1,r-1}}{Y^2_{2,r}}A^2+a^{-\Lm_{k-1}}\frac{Y_{1,r-1}}{Y_{2,k}}A(\varphi^G_{\textbf{V}})_{k-1})
}{(\varphi^G_{\textbf{V}})_{k-1}}\\
&=&\frac{a^{2\Lm_{k-1}+2\Lm_r}Y^2_{2,k}Y^2_{1,r}(1+A\times A^{-1}_{1,r})^2
(a^{\Lm_{k-2}-\Lm_{k-1}}\frac{Y_{2,k-1}}{Y_{2,k}}(\varphi^G_{\textbf{V}})_{k-1})}{(\varphi^G_{\textbf{V}})_{k-1}}\\
&+&\frac{a^{2\Lm_r+2\Lm_{k-1}+\Lm_{k-2}}Y_{1,k-1}Y_{2,k-1}Y_{2,k}(a^{-\Lm_{k-1}}\frac{Y_{1,r-1}}{Y_{2,k}}A(\varphi^G_{\textbf{V}})_{k-1})
}{(\varphi^G_{\textbf{V}})_{k-1}}\\
&=&a^{2\Lm_r+\Lm_{k-1}+\Lm_{k-2}}
(Y_{2,k-1}Y_{2,k}Y^2_{1,r}(1+A\times A^{-1}_{1,r})^2+Y_{1,k-1}Y_{2,k-1}Y_{1,r-1}A)\\
&=&a^{2\Lm_r+\Lm_{k-1}+\Lm_{k-2}}(Y_{2,k-1}Y_{2,k}Y^2_{1,r}+2\frac{Y_{2,k-1}Y_{2,k}Y_{1,r-1}Y_{1,r}}{Y_{2,r}}A)\\
&+& a^{2\Lm_r+\Lm_{k-1}+\Lm_{k-2}} (\frac{Y_{2,k-1}Y_{2,k}Y^2_{1,r-1}}{Y^2_{2,r}}A^2
+Y_{1,k-1}Y_{2,k-1}Y_{1,r-1}A)\\
&=&a^{2\Lm_r+\Lm_{k-1}+\Lm_{k-2}}(Y_{2,k-1}Y_{2,k}Y^2_{1,r}+2\frac{Y_{2,k-1}Y_{2,k}Y_{1,r-1}Y_{1,r}}{Y_{2,r}}A)\\
&+& a^{2\Lm_r+\Lm_{k-1}+\Lm_{k-2}} \frac{Y_{2,k-1}Y_{2,k}Y^2_{1,r-1}}{Y^2_{2,r}}A(A+A^{-1}_{1,r-1}A^{-1}_{1,r-2}\cdots A^{-1}_{1,k}).
\end{eqnarray*}
Thus, we obtain (\ref{B-lem3-exp}) for $j=k-1$.

Next, we consider the case $k-j>1$. The vertices and arrows around $(\varphi_{\textbf{V}})_{j}$ in the mutation diagram of $\mu_{j+1}\cdots\mu_{k-2}\mu_{k-1}\mu_{r}\cdots\mu_{k+1}\mu_k(\Sigma_0)$ are as follows:
\begin{equation}\label{same1}
\begin{xy}
(110,100) *{(\varphi_{\mu_{r-1}\cdots\mu_k(\textbf{V})})_{r-1}}="r-1",
(92,90)*{\cdots}="emp1",
(70,105)*{(\varphi_{\textbf{V}})_{r+j+1}}="r+k",
(70,90)*{ (\varphi_{\mu_{j+1}\cdots\mu_{k-1}\mu_{r}\cdots\mu_k(\textbf{V})})_{j+1}}="k",
(70,82)*{(\varphi_{\textbf{V}})_{-(j+1)}}="-j_k",
(35,105)*{(\varphi_{\textbf{V}})_{r+j}}="8",
(35,90) *{(\varphi_{\textbf{V}})_{j}}="4",
(35,82)*{(\varphi_{\textbf{V}})_{-j}}="-2",
(10,105)*{(\varphi_{\textbf{V}})_{r+j-1}}="6",
(10,90) *{(\varphi_{\textbf{V}})_{j-1}}="2",
(10,82)*{(\varphi_{\textbf{V}})_{-(j-1)}}="-1",
(0,90)*{\cdots}="emp",
\ar@{->} "4";"r-1"
\ar@{->} "k";"4"
\ar@{->} "r-1";"k"
\ar@{->} "-j_k";"k"
\ar@{->} "8";"k"
\ar@{->} "4";"-2"
\ar@{->} "2";"-1"
\ar@{->} "2";"4"
\ar@{->} "4";"6"
\ar@{->} "6";"2"
\ar@{->} "-2";"2"
\end{xy}
\end{equation}
By this diagram and the induction hypothesis,
\begin{eqnarray*}
& &\hspace{-20pt}(\varphi^G_{\mu_j\cdots\mu_{k-2}\mu_{k-1}\mu_{r}\mu_{r-1}\cdots\mu_{k+1}\mu_k(\textbf{V})})_j\\
&=& \frac{(\varphi^G_{\mu_{r-1}\cdots\mu_{k+1}\mu_k(\textbf{V})})_{r-1}
(\varphi^G_{\textbf{V}})_{r+j-1}(\varphi^G_{\textbf{V}})_{-j}
+(\varphi^G_{\textbf{V}})_{j-1}(\varphi^G_{\mu_{j+1}\cdots\mu_{k-2}\mu_{k-1}\mu_{r}\cdots\mu_{k+1}\mu_k(\textbf{V})})_{j+1}}{(\varphi^G_{\textbf{V}})_j}\\
&=&\frac{a^{\Lm_{j-1}+\Lm_j+2\Lm_r+\Lm_{k-1}}Y_{1,j}Y_{2,j}Y_{2,k}}{(\varphi^G_{\textbf{V}})_j}\left(Y^2_{1,r}+2\frac{Y_{1,r-1}Y_{1,r}}{Y_{2,r}}A\right.\\
&+&\left.
\frac{Y^2_{1,r-1}}{Y^2_{2,r}}A(1+A^{-1}_{1,r-1}+\cdots+A^{-1}_{1,r-1}\cdots A^{-1}_{1,1})\right)
+\frac{1}{(\varphi^G_{\textbf{V}})_j}(a^{\Lm_{j-1}-\Lm_j}\frac{Y_{2,j}(\varphi^G_{\textbf{V}})_j}{Y_{2,j+1}}-a^{\Lm_{j-1}}
\frac{Y_{1,j}Y_{2,j}}{Y_{2,j+1}})\\
&\times&a^{\Lm_{j}+\Lm_{k-1}+2\Lm_r}Y_{2,k}Y_{2,j+1}
(Y^2_{1,r}+2\frac{Y_{1,r-1}Y_{1,r}}{Y_{2,r}}A+\frac{Y^2_{1,r-1}}{Y^2_{2,r}}A(1+A^{-1}_{1,r-1}+\cdots+A^{-1}_{1,r-1}\cdots A^{-1}_{1,j+2}))\\
&=&a^{\Lm_{j-1}+2\Lm_r+\Lm_{k-1}}Y_{2,j}Y_{2,k}
(Y^2_{1,r}+2\frac{Y_{1,r-1}Y_{1,r}}{Y_{2,r}}A+\frac{Y^2_{1,r-1}}{Y^2_{2,r}}A(1+A^{-1}_{1,r-1}+\cdots+A^{-1}_{1,r-1}\cdots A^{-1}_{1,j+2}))\\
&+&\frac{a^{\Lm_{j-1}+2\Lm_r+\Lm_{k-1}}Y_{2,j}Y_{2,k}}{(\varphi^G_{\textbf{V}})_j}
a^{\Lm_j}\frac{Y_{1,j}Y^2_{1,r-1}}{Y^2_{2,r}}A\cdot A^{-1}_{1,r-1}\cdots A^{-1}_{1,j+1}
(1+A^{-1}_{1,j}+\cdots+A^{-1}_{1,j}\cdots A^{-1}_{1,1})\\
&=&a^{\Lm_{j-1}+2\Lm_r+\Lm_{k-1}}Y_{2,j}Y_{2,k}(Y^2_{1,r}+2\frac{Y_{1,r-1}Y_{1,r}}{Y_{2,r}}A+\frac{Y^2_{1,r-1}}{Y^2_{2,r}}A(1+A^{-1}_{1,r-1}+\cdots+A^{-1}_{1,r-1}\cdots A^{-1}_{1,j+1})).
\end{eqnarray*}
Hence, we obtain (\ref{B-lem3-exp}). By using (\ref{B-lem3-exp}), it is easy to see that 
\begin{eqnarray*}
& &\hspace{-20pt}a^{-(\Lm_{j-1}+\Lm_{k-1}+2\Lm_{r})}(\varphi^G_{\mu_j\cdots\mu_{k-2}\mu_{k-1}\mu_{r}\mu_{r-1}\cdots\mu_{k+1}\mu_k(\textbf{V})})_j(a;\textbf{Y})\\
&=& Y_{2,k} Y_{2,j}Y^2_{1,r}
+2Y_{2,k} Y_{2,j}Y^2_{1,r}A^{-1}_{1,r}
(1+A^{-1}_{1,r-1}+A^{-1}_{1,r-1}A^{-1}_{1,r-2}+\cdots+A^{-1}_{1,r-1}\cdots A^{-1}_{1,k+1})\\
&+&Y_{2,k} Y_{2,j}Y^2_{1,r}A^{-2}_{1,r}
(1+A^{-1}_{1,r-1}+A^{-1}_{1,r-1}A^{-1}_{1,r-2}+\cdots+A^{-1}_{1,r-1}\cdots A^{-1}_{1,k+1})^2\\
& &Y_{2,k} Y_{2,j}Y^2_{1,r}A^{-2}_{1,r}(1+A^{-1}_{1,r-1}+A^{-1}_{1,r-1}A^{-1}_{1,r-2}+\cdots+A^{-1}_{1,r-1}\cdots A^{-1}_{1,k+1})\\
& &\qquad \qquad \qquad \times A^{-1}_{1,r-1}A^{-1}_{1,r-2}\cdots A^{-1}_{1,k}(1+A^{-1}_{1,k-1}+\cdots+A^{-1}_{1,k-1}\cdots A^{-1}_{1,j+1}).
\end{eqnarray*}

Just as in the proof of (\ref{B-lem2pr-c}), we obtain
\begin{eqnarray*}
& &\hspace{-20pt}Y_{2,k} Y_{2,j}Y^2_{1,r}
+2Y_{2,k} Y_{2,j}Y^2_{1,r}A^{-1}_{1,r}
(1+A^{-1}_{1,r-1}+A^{-1}_{1,r-1}A^{-1}_{1,r-2}+\cdots+A^{-1}_{1,r-1}\cdots A^{-1}_{1,k+1})\\
&+&Y_{2,k} Y_{2,j}Y^2_{1,r}A^{-2}_{1,r}
(1+A^{-1}_{1,r-1}+A^{-1}_{1,r-1}A^{-1}_{1,r-2}+\cdots+A^{-1}_{1,r-1}\cdots A^{-1}_{1,k+1})^2\\
&=&\sum_{b\in B(\Lm_j+\Lambda_k+2\Lambda_{r})_{s_{k+1}s_{k+2}\cdots s_{r}} } C(b)\mu(b),
\end{eqnarray*}
where $\mu$ is the monomial realization in our claim, and each coefficient $C(b)$ is either $1$ or $2$. Similarly, we also obtain
\begin{eqnarray*}
& &\hspace{-20pt}Y_{2,k} Y_{2,j}Y^2_{1,r}A^{-2}_{1,r}(1+A^{-1}_{1,r-1}+A^{-1}_{1,r-1}A^{-1}_{1,r-2}+\cdots+A^{-1}_{1,r-1}\cdots A^{-1}_{1,k+1})\\
& &\qquad \qquad \qquad \times A^{-1}_{1,r-1}A^{-1}_{1,r-2}\cdots A^{-1}_{1,k}(1+A^{-1}_{1,k-1}+\cdots+A^{-1}_{1,k-1}\cdots A^{-1}_{1,j+1})\\
&=&Y_{1,k-1}Y_{2,j}Y_{1,r-1}(1+A^{-1}_{1,r-1}+A^{-1}_{1,r-1}A^{-1}_{1,r-2}+\cdots+A^{-1}_{1,r-1}\cdots A^{-1}_{1,k+1})\\
& &\qquad \qquad \qquad\times(1+A^{-1}_{1,k-1}+\cdots+A^{-1}_{1,k-1}\cdots A^{-1}_{1,j+1})\\
&=&\sum_{b\in B(\Lm_j+\Lambda_{k-1}+ \Lambda_{r-1})_{s_{j+1}\cdots s_{k-1}s_{k+1}s_{k+2}\cdots s_{r-1}}} \mu'(b).
\end{eqnarray*}
Therefore, we have (\ref{B-lem3-cry}). \qed

\subsection{The proof of Theorem \ref{thmC}}

Next, let us prove Theorem \ref{thmC}. It will be proved in a similar way to Theorem \ref{thmB}. Let $G={\rm Sp}_{2r}$. First, we set the Laurent monomials as follows:
\begin{equation}\label{ccbar}
C(l,k):=
\begin{cases}
\frac{Y_{l,k}}{Y_{l,k-1}} & {\rm if}\ 1\leq k\leq r, \\
\frac{Y_{l,|k|-1}}{Y_{l+1,|k|}} & {\rm if}\ \ovl{r}\leq k\leq \ovl{1},
\end{cases}
\end{equation}
where for $1\leq l\leq r$, set $|l|=|\ovl{l}|=l$. The mutation diagram (Definition \ref{wei-dia}) of the initial seed $\Sigma_0$ is
\begin{equation}\label{inidiaC-2}
\begin{xy}
(90,90)*{\cdots}="emp1",
(80,98)*{(\varphi_{\textbf{V}})_{r+k+1}}="7",
(80,90) *{(\varphi_{\textbf{V}})_{k+1}}="3",
(80,82)*{(\varphi_{\textbf{V}})_{-(k+1)}}="-4",
(55,98)*{(\varphi_{\textbf{V}})_{r+k}}="r+k",
(55,90) *{(\varphi_{\textbf{V}})_{k}}="k",
(55,82)*{(\varphi_{\textbf{V}})_{-k}}="-j_k",
(30,98)*{(\varphi_{\textbf{V}})_{r+k-1}}="8",
(30,90) *{(\varphi_{\textbf{V}})_{k-1}}="4",
(30,82)*{(\varphi_{\textbf{V}})_{-(k-1)}}="-2",
(0,98)*{(\varphi_{\textbf{V}})_{r+1}}="6",
(0,90) *{(\varphi_{\textbf{V}})_{1}}="2",
(0,82)*{(\varphi_{\textbf{V}})_{-1}}="-1",
(20,90)*{\ \cdots}="emp",
(20,98)*{}="emp2",
(20,82)*{}="emp3",
(125,98)*{(\varphi_{\textbf{V}})_{2r}}="2r",
(125,90) *{(\varphi_{\textbf{V}})_{r}}="r",
(125,82)*{(\varphi_{\textbf{V}})_{-r}}="-j_r",
(100,98)*{(\varphi_{\textbf{V}})_{2r-1}}="2r-1",
(100,90) *{(\varphi_{\textbf{V}})_{r-1}}="r-1",
(100,82)*{(\varphi_{\textbf{V}})_{-(r-1)}}="-r+1",
\ar@{->} "7";"3"
\ar@{->} "r+k";"k"
\ar@{->} "8";"4"
\ar@{->} "6";"2"
\ar@{->} "3";"-4"
\ar@{->} "k";"-j_k"
\ar@{->} "k";"8"
\ar@{->} "4";"-2"
\ar@{->} "2";"-1"
\ar@{->} "k";"3"
\ar@{->} "4";"k"
\ar@{->} "-j_k";"4"
\ar@{->} "2";"emp"
\ar@{->} "3";"r+k"
\ar@{->} "emp";"6"
\ar@{->} "-4";"k"
\ar@{->} "emp3";"2"
\ar@{->} "2r";"r"
\ar@{->} "2r-1";"r-1"
\ar@{->} "r";"-j_r"
\ar@{->} "r";"2r-1"_{-2}
\ar@{->} "r-1";"-r+1"
\ar@{->} "r-1";"r"^{2}
\ar@{->} "-j_r";"r-1"^{-2}
\end{xy}
\end{equation}

\begin{prop}\label{iniC}
\begin{enumerate}
\item For $k\in[1,r]$,
the initial cluster variables $(\varphi^G_{\textbf{V}})_k(a;\textbf{Y})$ in $\mathbb{C}[G^{e,c^2}]$ are
described as
\[ 
(\varphi^G_{\textbf{V}})_k(a;\textbf{Y})=a^{\Lambda_k} \sum_{b\in B(\Lambda_k)_{c^2_{>2r-k}}} C(b)\mu(b),
\]
where $\mu:B(\Lambda_k)\rightarrow \mathcal{Y}$ is the monomial realization of $B(\Lambda_k)$ in Theorem \ref{thmC} $(i)$, $C(b)$ are some positive integers.
\item
For $k\in[1,r]$, the frozen cluster variables $(\varphi^G_{\textbf{V}})_{-k}(a;\textbf{Y})$, $(\varphi^G_{\textbf{V}})_{r+k}(a;\textbf{Y})$ in $\mathbb{C}[G^{e,c^2}]$ are described as
\[
(\varphi^G_{\textbf{V}})_{-k}(a;\textbf{Y})=a^{\Lm_k}Y_{1,k}Y_{2,k},
\quad
(\varphi^G_{\textbf{V}})_{r+k}(a;\textbf{Y})=a^{\Lm_k}.
\]
\end{enumerate}
\end{prop}
\nd
[{\sl Proof.}]

For $1\leq k\leq r-1$, we can prove 
\begin{eqnarray}
(\varphi^G_{\textbf{V}})_k(a;\textbf{Y})&=&
a^{\Lambda_k} Y_{1,k}(1+A^{-1}_{1,k}+A^{-1}_{1,k}A^{-1}_{1,k-1}+\cdots+A^{-1}_{1,k}\cdots A^{-1}_{1,1})\qquad \label{iniC-exp1} \\
&=&a^{\Lambda_k} \sum_{b\in B(\Lambda_k)_{c^2_{>2r-k}}}\mu(b), \nonumber
\end{eqnarray}
by the same argument as in Proposition \ref{iniB}.

Next, let us consider the case $k=r$. Just as in (\ref{B-ext}), it follows
\[
\ovl{c^2_{>r}}v_1\wedge\cdots\wedge v_r=\ovl{s_1}\cdots \ovl{s_r}(v_1\wedge\cdots\wedge v_r)
=v_2\wedge\cdots\wedge v_{r}\wedge v_{\ovl{1}}. 
\]
Using the bilinear form in \ref{bilingen}, the cluster variable $(\varphi^G_{\textbf{V}})_r(a;\textbf{Y})$ is described as
\begin{multline*}
a^{\Lm_r}\lan v_2\wedge\cdots\wedge v_{r}\wedge v_{\ovl{1}}, \al^{\vee}_1(Y_{2,1})y_1(Y_{2,1})\cdots \al^{\vee}_r(Y_{2,r})y_r(Y_{2,r})\\ \al^{\vee}_1(Y_{1,1})y_1(Y_{1,1})\cdots\al^{\vee}_r(Y_{1,r})y_r(Y_{1,r})v_1\wedge\cdots\wedge v_r \ran.
\end{multline*}
Using $(\ref{C-wtv})$, $(\ref{C-f1})$ and $(\ref{C-f2})$ repeatedly, one obtain
\[
 \al^{\vee}_1(Y_{1,1})y_1(Y_{1,1})\cdots\al^{\vee}_r(Y_{1,r})y_r(Y_{1,r})v_i=
\begin{cases}
\frac{Y_{1,i}}{Y_{1,i-1}}v_i+v_{i+1} & {\rm if}\ 1\leq i\leq r-1, \\
\frac{Y_{1,r}}{Y_{1,r-1}}v_r+\sum^r_{j=1} Y_{1,j-1}v_{\ovl{j}} & {\rm if}\ i=r,
\end{cases}
\]
which means that $\al^{\vee}_1(Y_{1,1})y_1(Y_{1,1})\cdots\al^{\vee}_r(Y_{1,r})y_r(Y_{1,r})v_1\wedge\cdots\wedge v_r$ is a linear combination of $v_1\wedge\cdots \wedge v_r$ and $v_1\wedge\cdots \wedge v_s\wedge v_{s+2}\wedge\cdots \wedge v_{r}\wedge v_{\ovl{t+1}}$ $(0\leq s,t\leq r-1)$ with the coefficients $C(1,1)C(1,2)\cdots C(1,r)$ and $C(1,1)C(1,2)\cdots C(1,s)Y_{1,t}$, respectively. Similarly, we see that the coefficient of $v_2\wedge\cdots\wedge v_r\wedge v_{\ovl{1}}$ in $\al^{\vee}_1(Y_{2,1})$ $y_1(Y_{2,1})\cdots\al^{\vee}_r(Y_{2,r})y_r(Y_{2,r})v_1\wedge\cdots \wedge v_r$ is $1$, and the coefficient in $\al^{\vee}_1(Y_{2,1})y_1(Y_{2,1})\cdots\al^{\vee}_r(Y_{2,r})y_r(Y_{2,r})v_1\wedge\cdots \wedge v_s\wedge v_{s+2}\wedge\cdots \wedge v_{r}\wedge v_{\ovl{t+1}}$ is $C(2,s+2)\cdots C(2,r)Y^{-1}_{2,t+1}$. By the above argument, we obtain
\begin{eqnarray}
& &\hspace{-20pt}(\varphi^G_{\textbf{V}})_r(a;\textbf{Y})\nonumber\\
&=&a^{\Lm_r}(Y_{1,r}
+\sum_{0\leq s,t \leq r-1} C(1,1)C(1,2)\cdots C(1,s)
C(2,s+2)\cdots C(2,r)C(1,\ovl{t+1}))\nonumber\\
&=&a^{\Lm_r}(Y_{1,r}+Y_{1,r}A^{-1}_{1,r}(1+A^{-1}_{1,r-1}+A^{-1}_{1,r-1}A^{-1}_{1,r-2}+\cdots+A^{-1}_{1,r-1}\cdots A^{-1}_{1,1})^2).\qquad \quad \label{iniC-exp2}
\end{eqnarray}

From Theorem \ref{kashidem}, the conclusion $(\varphi^G_{\textbf{V}})_r(a;\textbf{Y})=a^{\Lambda_r} \sum_{b\in B(\Lambda_k)_{c^2_{>r}}} C(b)\mu(b)$ follows, where each coefficient $C(b)$ is either $1$ or $2$. Thus, we obtain the claim (i). By the same argument, we can also obtain (ii). \qed

By the same way as in (\ref{twovarsB}), we see that for $k\in[1,r-2]$, 
\begin{equation}\label{twovarsC}
a^{-\Lambda_k}(\varphi^G_{\textbf{V}})_k(a;\textbf{Y})=Y_{1,k}+C(2,k+1)(a^{-\Lm_{k-1}}(\varphi^G_{\textbf{V}})_{k-1}(a;\textbf{Y})).
\end{equation}

\begin{lem}\label{C-lem1}
For $k\in[1,r]$, the cluster variables $(\varphi^G_{\mu_k\mu_{k+1}\cdots\mu_r(\textbf{V})})_k(a;\textbf{Y})$ in $\mathbb{C}[G^{e,c^2}]$ are
described as
\[ 
(\varphi^G_{\mu_k\mu_{k+1}\cdots\mu_r(\textbf{V})})_k(a;\textbf{Y})=
a^{\Lambda_{k-1}+\Lambda_{r-1}} Y_{2,k}.
\]
\end{lem}
\nd
[{\sl Proof.}]

Using the mutation diagram (\ref{inidiaC-2}), we can prove this lemma in the same way as Lemma \ref{B-lem1}. \qed

\vspace{3mm}

\nd
[{\sl Proof of Theorem \ref{thmC} (ii) and (iii).}]

Our claim (ii) is obtained by the same calculation as in Theorem \ref{thmB} (ii). So let us consider the claim (iii). By Lemma \ref{mutgamlem}, the mutation diagram of $\mu_{r-2}\cdots\mu_{k+1}\mu_k(\Sigma_0)$ is as follows:
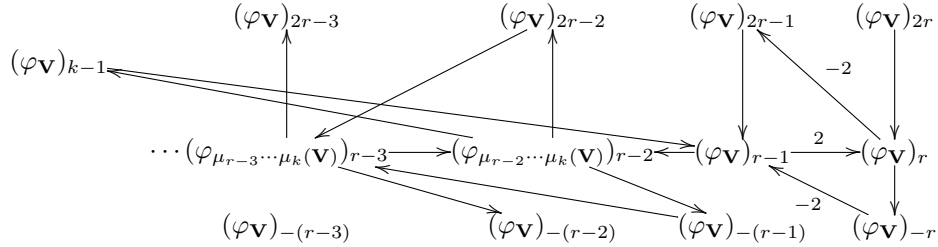
\[
\begin{xy}
(80,108)*{(\varphi_{\textbf{V}})_{2r}}="7",
(80,90) *{(\varphi_{\textbf{V}})_{r}}="3",
(80,80)*{(\varphi_{\textbf{V}})_{-r}}="-4",
(60,108)*{(\varphi_{\textbf{V}})_{2r-1}}="r+k",
(60,90) *{(\varphi_{\textbf{V}})_{r-1}}="k",
(60,80)*{(\varphi_{\textbf{V}})_{-(r-1)}}="-j_k",
(35,108)*{(\varphi_{\textbf{V}})_{2r-2}}="8",
(35,90) *{(\varphi_{\mu_{r-2}\cdots\mu_k(\textbf{V})})_{r-2}}="4",
(35,80)*{(\varphi_{\textbf{V}})_{-(r-2)}}="-2",
(0,108)*{(\varphi_{\textbf{V}})_{2r-3}}="6",
(0,90) *{(\varphi_{\mu_{r-3}\cdots\mu_k(\textbf{V})})_{r-3}}="2",
(0,80)*{(\varphi_{\textbf{V}})_{-(r-3)}}="-1",
(-16,90)*{\cdots}="emp",
(-30,102)*{(\varphi_{\textbf{V}})_{k-1}}="k-1",
\ar@{->} "7";"3"
\ar@{->} "r+k";"k"
\ar@{->} "4";"8"
\ar@{->} "2";"6"
\ar@{->} "3";"-4"
\ar@{->} "k";"3"^{2}
\ar@{->} "k";"4"
\ar@{->} "4";"-j_k"
\ar@{->} "2";"4"
\ar@{->} "3";"r+k"_{-2}
\ar@{->} "8";"2"
\ar@{->} "-4";"k"^{-2}
\ar@{->} "k-1";"k"
\ar@{->} "4";"k-1"
\ar@{->} "-j_k";"2"
\ar@{->} "2";"-2"
\end{xy}
\]
From this diagram, we have
\begin{equation}\label{C-lem2-pr1}
(\varphi^G_{\mu_{r-1}\cdots\mu_{k+1}\mu_k(\textbf{V})})_{r-1}
=\frac{(\varphi^G_{\textbf{V}})_{r}(\varphi^G_{\mu_{r-2}\cdots\mu_{k+1}\mu_k(\textbf{V})})_{r-2}+(\varphi^G_{\textbf{V}})_{k-1}(\varphi^G_{\textbf{V}})_{2r-1}(\varphi^G_{\textbf{V}})_{-r}}{(\varphi^G_{\textbf{V}})_{r-1}}.
\end{equation}
Using the notation  $A:=(1+A^{-1}_{1,k-1}+A^{-1}_{1,k-1}A^{-1}_{1,k-2}+\cdots+A^{-1}_{1,k-1}A^{-1}_{1,k-2}\cdots A^{-1}_{1,1})$,
we can write
\begin{equation}\label{C-lem2-pr2}
(\varphi^G_{\textbf{V}})_{k-1}
=a^{\Lm_{k-1}}Y_{1,k-1}A.
\end{equation}
Note that by (\ref{iniC-exp1}) and (\ref{C-exp1}),
\begin{eqnarray}
& &\hspace{-20pt}(\varphi^G_{\mu_{r-2}\cdots\mu_{k+1}\mu_k(\textbf{V})})_{r-2}(a;\textbf{Y}) \nonumber\\ 
&=&
a^{\Lm_{k-1}+\Lm_{r-1}} Y_{2,k}(a^{-\Lm_{r-1}}(\varphi^G_{\textbf{V}})_{r-1}-Y_{1,r-1}A^{-1}_{1,r-1}A^{-1}_{1,r-2}\cdots A^{-1}_{1,k} \cdot A),
\qquad \quad\nonumber \\
&=&a^{\Lm_{k-1}+\Lm_{r-1}} Y_{2,k}(a^{-\Lm_{r-1}}(\varphi^G_{\textbf{V}})_{r-1}-\frac{Y_{1,k-1}Y_{2,r}}{Y_{2,k}} \cdot A) \label{C-lem2-pr3}
\end{eqnarray}
and by (\ref{iniC-exp1}) and (\ref{iniC-exp2}),
\begin{equation}\label{C-lem2-pr4}
(\varphi^G_{\textbf{V}})_{r}=a^{\Lm_r}(Y_{1,r}+a^{-2\Lm_{r-1}}Y^{-1}_{2,r}(\varphi^G_{\textbf{V}})^2_{r-1}).
\end{equation}

Substituting (\ref{C-lem2-pr2}), (\ref{C-lem2-pr3}) and (\ref{C-lem2-pr4}) for (\ref{C-lem2-pr1}), we get
\begin{eqnarray}
& &\hspace{-20pt}(\varphi^G_{\mu_{r-1}\cdots\mu_{k+1}\mu_k(\textbf{V})})_{r-1}(a;\textbf{Y}) \nonumber\\ 
&=&\frac{a^{\Lm_{k-1}+\Lm_r}Y_{2,k}(\varphi^G_{\textbf{V}})_{r-1}(Y_{1,r}+a^{-2\Lm_{r-1}}Y^{-1}_{2,r}(\varphi^G_{\textbf{V}})^2_{r-1}) 
-Y_{1,k-1}A a^{\Lm_{k-1}-\Lm_{r-1}+\Lm_r}(\varphi^G_{\textbf{V}})^2_{r-1}}{(\varphi^G_{\textbf{V}})_{r-1}}\nonumber\\
&=&
a^{\Lm_{k-1}+\Lm_r}Y_{2,k}(Y_{1,r}+a^{-2\Lm_{r-1}}Y^{-1}_{2,r}(\varphi^G_{\textbf{V}})^2_{r-1})-a^{\Lm_{k-1}+\Lm_r-\Lm_{r-1}}Y_{1,k-1}A(\varphi^G_{\textbf{V}})_{r-1}\label{C-lem2-pr5}\\
&=&a^{\Lm_{k-1}+\Lm_r}Y_{2,k}Y_{1,r}\nonumber\\
&+&a^{\Lm_{k-1}+\Lm_r-\Lm_{r-1}}(\varphi^G_{\textbf{V}})_{r-1}(
\frac{Y_{1,r-1}Y_{2,k}}{Y_{2,r}}(1+A^{-1}_{1,r-1}+\cdots+A^{-1}_{1,r-1}\cdots A^{-1}_{1,1})-Y_{1,k-1}A)\nonumber\\
&=&a^{\Lm_{k-1}+\Lm_r}Y_{2,k}Y_{1,r}+a^{\Lm_{k-1}+\Lm_r}\frac{Y_{2,k}Y^2_{1,r-1}}{Y_{2,r}}(1+A^{-1}_{1,r-1}+\cdots+A^{-1}_{1,r-1}\cdots A^{-1}_{1,k+1})\label{C-lem2-pr5-1}\\
& &\qquad \qquad \qquad \qquad \qquad \qquad \qquad \qquad \times (1+A^{-1}_{1,r-1}+\cdots+A^{-1}_{1,r-1}\cdots A^{-1}_{1,1})\nonumber\\
&=&a^{\Lm_{k-1}+\Lm_r}\left(Y_{2,k}Y_{1,r}+Y_{2,k}Y_{1,r}A^{-1}_{1,r}(1+A^{-1}_{1,r-1}+\cdots+A^{-1}_{1,r-1}\cdots A^{-1}_{1,k+1})^2\right.\nonumber\\
&+&\left.Y_{1,k-1}Y_{1,r-1}(1+A^{-1}_{1,r-1}+\cdots+A^{-1}_{1,r-1}\cdots A^{-1}_{1,k+1})(1+A^{-1}_{1,k-1}+\cdots+A^{-1}_{1,k-1}\cdots A^{-1}_{1,1})\right)\nonumber,
\end{eqnarray}
which implies (\ref{C-cry2}).

Next, let us consider the cluster variable $(\varphi^G_{\mu_r\mu_{r-1}\cdots\mu_{k+1}\mu_k(\textbf{V})})_{r}(a;\textbf{Y})$.
The mutation diagram of $\mu_{r-1}\cdots\mu_{k+1}\mu_k(\Sigma_0)$ is as follows:
\[
\begin{xy}
(93,108)*{(\varphi_{\textbf{V}})_{2r}}="7",
(93,90) *{(\varphi_{\textbf{V}})_{r}}="3",
(93,77)*{(\varphi_{\textbf{V}})_{-r}}="-4",
(68,108)*{(\varphi_{\textbf{V}})_{2r-1}}="r+k",
(68,90) *{(\varphi_{\mu_{r-1}\cdots\mu_k(\textbf{V})})_{r-1}}="k",
(68,77)*{(\varphi_{\textbf{V}})_{-(r-1)}}="-j_k",
(35,108)*{(\varphi_{\textbf{V}})_{2r-2}}="8",
(35,90) *{(\varphi_{\mu_{r-2}\cdots\mu_k(\textbf{V})})_{r-2}}="4",
(35,77)*{(\varphi_{\textbf{V}})_{-(r-2)}}="-2",
(0,108)*{(\varphi_{\textbf{V}})_{2r-3}}="6",
(0,90) *{(\varphi_{\mu_{r-3}\cdots\mu_k(\textbf{V})})_{r-3}}="2",
(0,77)*{(\varphi_{\textbf{V}})_{-(r-3)}}="-1",
(-16,90)*{\cdots}="emp",
(-30,104)*{(\varphi_{\textbf{V}})_{k-1}}="k-1",
\ar@{->} "7";"3"
\ar@{->} "k";"r+k"
\ar@{->} "4";"8"
\ar@{->} "2";"6"
\ar@{->} "-4";"3"
\ar@{->} "-4";"4"_{-2}
\ar@{->} "3";"k"^{\qquad \quad -2}
\ar@{->} "4";"k"
\ar@{->} "4";"-j_k"
\ar@{->} "2";"4"
\ar@{->} "r+k";"4"
\ar@{->} "8";"2"
\ar@{->} "k";"-4"^{2}
\ar@{->} "k";"k-1"
\ar@{->} "k-1";"3"^2
\ar@{->} "-j_k";"2"
\ar@{->} "2";"-2"
\end{xy}
\]
Thus,
\begin{equation}\label{C-lem2-pr6}
(\varphi^G_{\mu_r\mu_{r-1}\cdots\mu_{k+1}\mu_k(\textbf{V})})_{r}
=\frac{(\varphi^G_{\mu_{r-1}\cdots\mu_{k+1}\mu_k(\textbf{V})})^2_{r-1}+
(\varphi^G_{\textbf{V}})^2_{k-1}(\varphi^G_{\textbf{V}})_{2r}(\varphi^G_{\textbf{V}})_{-r}}{(\varphi^G_{\textbf{V}})_{r}}.
\end{equation}
From (\ref{C-lem2-pr4}) and (\ref{C-lem2-pr5}), we see that
\begin{eqnarray*}
& &\hspace{-20pt}(\varphi^G_{\mu_{r-1}\cdots\mu_{k+1}\mu_k(\textbf{V})})^2_{r-1} \\ 
&=&a^{2\Lm_{k-1}}(Y_{2,k}(\varphi^G_{\textbf{V}})_{r}-a^{\Lm_r-\Lm_{r-1}}Y_{1,k-1}A(\varphi^G_{\textbf{V}})_{r-1})^2\\
&=&a^{2\Lm_{k-1}}(Y^2_{2,k}(\varphi^G_{\textbf{V}})^2_{r}-2a^{\Lm_r-\Lm_{r-1}}Y_{2,k}Y_{1,k-1}A(\varphi^G_{\textbf{V}})_{r-1}(\varphi^G_{\textbf{V}})_{r}+a^{2\Lm_r-2\Lm_{r-1}}Y^2_{1,k-1}A^2(\varphi^G_{\textbf{V}})^2_{r-1})\\
&=&a^{2\Lm_{k-1}}(\varphi^G_{\textbf{V}})_{r}(Y^2_{2,k}(\varphi^G_{\textbf{V}})_r
-2a^{\Lm_r-\Lm_{r-1}}Y_{2,k}Y_{1,k-1}A(\varphi^G_{\textbf{V}})_{r-1}
+a^{\Lm_r}Y^2_{1,k-1}Y_{2,r}A^2)\\
& &-a^{2\Lm_{k-1}+2\Lm_r}Y^2_{1,k-1}Y_{1,r}Y_{2,r}A^2.
\end{eqnarray*}
Since we know that $(\varphi^G_{\textbf{V}})_{k-1}=a^{\Lm_{k-1}}Y_{1,k-1}A$ (\ref{iniC-exp1}), from (\ref{C-lem2-pr6}), one obtain
\begin{eqnarray}
& &\hspace{-20pt}(\varphi^G_{\mu_r\mu_{r-1}\cdots\mu_{k+1}\mu_k(\textbf{V})})_{r}\nonumber\\
&=&a^{2\Lm_{k-1}}(Y^2_{2,k}(\varphi^G_{\textbf{V}})_r
-2a^{\Lm_r-\Lm_{r-1}}Y_{2,k}Y_{1,k-1}A(\varphi^G_{\textbf{V}})_{r-1}
+a^{\Lm_r}Y^2_{1,k-1}Y_{2,r}A^2). \qquad \quad \label{C-lem2-pr7}
\end{eqnarray}
By (\ref{iniC-exp1}) and (\ref{iniC-exp2}), we get
\begin{eqnarray*}
& &\hspace{-20pt}Y^2_{2,k}(\varphi^G_{\textbf{V}})_r
-a^{\Lm_r-\Lm_{r-1}}Y_{2,k}Y_{1,k-1}A(\varphi^G_{\textbf{V}})_{r-1}\\
&=&a^{\Lm_r}Y_{1,r}Y^2_{2,k}+a^{\Lm_r}Y_{1,r}Y^2_{2,k}A^{-1}_{1,r}
(1+A^{-1}_{1,r-1}+\cdots+A^{-1}_{1,r-1}\cdots A^{-1}_{1,1})\\
& &\qquad\qquad\qquad\qquad\qquad\times (1+A^{-1}_{1,r-1}+\cdots+A^{-1}_{1,r-1}\cdots A^{-1}_{1,k+1}),
\end{eqnarray*}
and
\begin{eqnarray*}
& &\hspace{-20pt}-a^{\Lm_r-\Lm_{r-1}}Y_{2,k}Y_{1,k-1}A(\varphi^G_{\textbf{V}})_{r-1}
+a^{\Lm_r}Y^2_{1,k-1}Y_{2,r}A^2
\\
&=&-a^{\Lm_r}Y_{2,k}Y_{1,k-1}Y_{1,r-1}A (1+A^{-1}_{1,r-1}+\cdots+A^{-1}_{1,r-1}\cdots A^{-1}_{1,k+1}).
\end{eqnarray*}
Thus, it follows
\begin{eqnarray}
& &\hspace{-20pt}(\varphi^G_{\mu_r\mu_{r-1}\cdots\mu_{k+1}\mu_k(\textbf{V})})_{r}\nonumber\\
&=&
a^{\Lm_r+2\Lm_{k-1}}Y_{1,r}Y^2_{2,k}(1+A^{-1}_{1,r}(1+A^{-1}_{1,r-1}+\cdots+A^{-1}_{1,r-1}\cdots A^{-1}_{1,k+1})^2).\qquad \quad \nonumber
\end{eqnarray}
From this explicit formula, the conclusion (\ref{C-cry3}) follows. \qed

\vspace{3mm}

\nd
[{\sl Proof of Theorem \ref{thmC} (iv).}]

To prove our claim, we use the induction on $(k-j)$. First, let $k-j=1$. It follows by Lemma \ref{mutgamlem} that the mutation diagram of $\mu_{r}\mu_{r-1}\cdots\mu_{k+1}\mu_k(\Sigma_0)$ is
\[
\begin{xy}
(110,103) *{(\varphi_{\mu_{r}\cdots\mu_k(\textbf{V})})_{r}}="r",
(90,95) *{(\varphi_{\mu_{r-1}\cdots\mu_k(\textbf{V})})_{r-1}}="r-1",
(72,90)*{\cdots}="emp1",
(60,106)*{(\varphi_{\textbf{V}})_{r+k}}="r+k",
(60,90) *{(\varphi_{\mu_k(\textbf{V})})_{k}}="k",
(60,82)*{(\varphi_{\textbf{V}})_{-k}}="-j_k",
(35,106)*{(\varphi_{\textbf{V}})_{r+k-1}}="8",
(35,90) *{(\varphi_{\textbf{V}})_{k-1}}="4",
(35,82)*{(\varphi_{\textbf{V}})_{-(k-1)}}="-2",
(10,106)*{(\varphi_{\textbf{V}})_{r+k-2}}="6",
(10,90) *{(\varphi_{\textbf{V}})_{k-2}}="2",
(10,82)*{(\varphi_{\textbf{V}})_{-(k-2)}}="-1",
(0,90)*{\cdots}="emp",
\ar@{->} "r";"4"_{-2}
\ar@{->} "r-1";"r"_2
\ar@{->} "4";"r-1"
\ar@{->} "k";"r+k"
\ar@{->} "6";"2"
\ar@{->} "-j_k";"k"
\ar@{->} "8";"k"
\ar@{->} "4";"-2"
\ar@{->} "2";"-1"
\ar@{->} "2";"4"
\ar@{->} "4";"6"
\ar@{->} "-2";"2"
\end{xy},
\]
which means that
\begin{eqnarray*}
& &\hspace{-20pt}(\varphi^G_{\mu_{k-1}\mu_{r}\mu_{r-1}\cdots\mu_{k+1}\mu_k(\textbf{V})})_{k-1}\\
&=&
\frac{(\varphi^G_{\mu_{r}\mu_{r-1}\cdots\mu_{k+1}\mu_k(\textbf{V})})_{r}(\varphi^G_{\textbf{V}})_{k-2}+(\varphi^G_{\mu_{r-1}\cdots\mu_{k+1}\mu_k(\textbf{V})})_{r-1}(\varphi^G_{\textbf{V}})_{r+k-2}(\varphi^G_{\textbf{V}})_{-(k-1)}
}{(\varphi^G_{\textbf{V}})_{k-1}}.
\end{eqnarray*}
By (\ref{C-lem2-pr4}) and (\ref{C-lem2-pr5}), we get
\[
(\varphi^G_{\mu_{r-1}\cdots\mu_{k+1}\mu_k(\textbf{V})})_{r-1}
=a^{\Lm_{k-1}}Y_{2,k}(\varphi^G_{\textbf{V}})_{r}-a^{\Lm_{r}-\Lm_{r-1}}
(\varphi^G_{\textbf{V}})_{k-1}(\varphi^G_{\textbf{V}})_{r-1},
\]
and (\ref{C-lem2-pr7}) implies that
\begin{eqnarray}
& &\hspace{-20pt}(\varphi^G_{\mu_r\mu_{r-1}\cdots\mu_{k+1}\mu_k(\textbf{V})})_{r}
=a^{2\Lm_{k-1}}Y^2_{2,k}(\varphi^G_{\textbf{V}})_{r}\nonumber\\
&-&2a^{\Lm_{r}-\Lm_{r-1}+\Lm_{k-1}}Y_{2,k}
(\varphi^G_{\textbf{V}})_{k-1}(\varphi^G_{\textbf{V}})_{r-1}+a^{\Lm_{r}}Y_{2,r}
(\varphi^G_{\textbf{V}})^2_{k-1}=a^{\Lm_{k-1}}Y_{2,k}(\varphi^G_{\mu_{r-1}\cdots\mu_{k+1}\mu_k(\textbf{V})})_{r-1}\nonumber\\
&-&
a^{\Lm_{r}-\Lm_{r-1}+\Lm_{k-1}}Y_{2,k}
(\varphi^G_{\textbf{V}})_{k-1}(\varphi^G_{\textbf{V}})_{r-1}+a^{\Lm_{r}}Y_{2,r}
(\varphi^G_{\textbf{V}})^2_{k-1}. \label{C-lem3-pr1}
\end{eqnarray}
As have seen in (\ref{twovarsC}), it follows that
\begin{equation}\label{C-lem3-pr2}
(\varphi^G_{\textbf{V}})_{k-2}=a^{\Lm_{k-2}-\Lm_{k-1}}
\frac{Y_{2,k-1}}{Y_{2,k}}(\varphi^G_{\textbf{V}})_{k-1}-a^{\Lm_{k-2}}\frac{Y_{1,k-1}Y_{2,k-1}}{Y_{2,k}}.
\end{equation}
By using the formulas $(\ref{C-lem3-pr1})$ and $(\ref{C-lem3-pr2})$, we get
\begin{eqnarray*}
& &\hspace{-20pt}(\varphi^G_{\mu_{k-1}\mu_{r}\mu_{r-1}\cdots\mu_{k+1}\mu_k(\textbf{V})})_{k-1}\\
&=&
a^{\Lm_{k-2}-\Lm_{k-1}}\frac{Y_{2,k-1}}{Y_{2,k}}
(\varphi^G_{\mu_r\mu_{r-1}\cdots\mu_{k+1}\mu_k(\textbf{V})})_{r}
+a^{\Lm_{r}-\Lm_{r-1}+\Lm_{k-1}+\Lm_{k-2}}Y_{1,k-1}Y_{2,k-1}(\varphi^G_{\textbf{V}})_{r-1}\\
&-&a^{\Lm_{r}+\Lm_{k-2}}\frac{Y_{1,k-1}Y_{2,k-1}Y_{2,r}}{Y_{2,k}}(\varphi^G_{\textbf{V}})_{k-1}\\
&=&a^{\Lm_{r}+\Lm_{k-1}+\Lm_{k-2}}Y_{2,k-1}Y_{2,k}Y_{1,r}(1+A^{-1}_{1,r}(1+A^{-1}_{1,r-1}+\cdots+A^{-1}_{1,r-1}\cdots A^{-1}_{1,k+1})^2)\\
&+&a^{\Lm_{r}+\Lm_{k-1}+\Lm_{k-2}}Y_{1,k-1}Y_{2,k-1}Y_{1,r-1}(1+A^{-1}_{1,r-1}+\cdots+A^{-1}_{1,r-1}\cdots A^{-1}_{1,k+1})\\
&=&a^{\Lm_{r}+\Lm_{k-1}+\Lm_{k-2}}Y_{2,k-1}Y_{2,k}Y_{1,r}
+a^{\Lm_{r}+\Lm_{k-1}+\Lm_{k-2}}Y_{2,k-1}Y_{2,k}Y_{1,r}A^{-1}_{1,r}\\
& &\times(1+A^{-1}_{1,r-1}+\cdots+A^{-1}_{1,r-1}\cdots A^{-1}_{1,k+1}) (1+A^{-1}_{1,r-1}+\cdots+A^{-1}_{1,r-1}\cdots A^{-1}_{1,k}),
\end{eqnarray*}
where we use (\ref{iniC-exp1}) and (\ref{C-lem2-pr8}) in the second equality. Thus, we obtain our claim (\ref{C-lem3-exp}) for $j=k-1$.

Next, we assume $k-j>1$. The mutation diagram of $\mu_{j+1}\cdots\mu_{k-2}\mu_{k-1}\mu_{r}\mu_{r-1}\cdots\mu_{k+1}\mu_k(\Sigma_0)$ is the same form as in (\ref{same1}), which means
\begin{eqnarray*}
& &\hspace{-20pt}(\varphi^G_{\mu_j\mu_{j+1}\cdots\mu_{k-2}\mu_{k-1}\mu_{r}\mu_{r-1}\cdots\mu_{k+1}\mu_k(\textbf{V})})_j(a;\textbf{Y})\\
&=&\frac{(\varphi^G_{\mu_{r-1}\cdots\mu_{k+1}\mu_k(\textbf{V})})_{r-1}(\varphi^G_{\textbf{V}})_{r+j-1}(\varphi^G_{\textbf{V}})_{-j}+
(\varphi^G_{\mu_{j+1}\cdots\mu_{k-2}\mu_{k-1}\mu_{r}\mu_{r-1}\cdots\mu_{k+1}\mu_k(\textbf{V})})_{j+1}
(\varphi^G_{\textbf{V}})_{j-1}}{(\varphi^G_{\textbf{V}})_j}.
\end{eqnarray*}
By the induction hypothesis, we get
\begin{eqnarray}
& &\hspace{-20pt}(\varphi^G_{\mu_{j+1}\cdots\mu_{k-2}\mu_{k-1}\mu_{r}\mu_{r-1}\cdots\mu_{k+1}\mu_k(\textbf{V})})_{j+1}(a;\textbf{Y})\nonumber\\
&=& a^{\Lm_{j}+\Lm_{k-1}+\Lm_{r}}Y_{2,k}Y_{2,j+1}Y_{1,r}\nonumber\\
&+&a^{\Lm_{j}+\Lm_{k-1}+\Lm_{r}}\frac{Y_{2,k}Y_{2,j+1}Y^2_{1,r-1}}{Y_{2,r}}(1+A^{-1}_{1,r-1}+A^{-1}_{1,r-1}A^{-1}_{1,r-2}+\cdots+A^{-1}_{1,r-1}\cdots A^{-1}_{1,k+1})\qquad \quad \label{C-lem3-pr3}\\
& &\qquad\qquad\qquad\qquad\qquad\qquad \times
(1+A^{-1}_{1,r-1}+A^{-1}_{1,r-1}A^{-1}_{1,r-2}+\cdots+A^{-1}_{1,r-1}\cdots A^{-1}_{1,j+2}).\nonumber
\end{eqnarray}
Since the equation $(\varphi^G_{\textbf{V}})_{j-1}=a^{\Lm_{j-1}-\Lm_j}\frac{Y_{2,j}}{Y_{2,j+1}}(\varphi^G_{\textbf{V}})_j-a^{\Lm_{j-1}}\frac{Y_{1,j}Y_{2,j}}{Y_{2,j+1}}$ holds from $(\ref{twovarsC})$, it follows by (\ref{C-lem2-pr5-1}) and (\ref{C-lem3-pr3}) that
\begin{eqnarray*}
& &\hspace{-20pt}(\varphi^G_{\mu_j\mu_{j+1}\cdots\mu_{k-2}\mu_{k-1}\mu_{r}\mu_{r-1}\cdots\mu_{k+1}\mu_k(\textbf{V})})_j(a;\textbf{Y})\\
&=&\frac{a^{\Lm_{j-1}-\Lm_j}\frac{Y_{2,j}}{Y_{2,j+1}}
(\varphi^G_{\mu_{j+1}\cdots\mu_{k-2}\mu_{k-1}\mu_{r}\mu_{r-1}\cdots\mu_{k+1}\mu_k(\textbf{V})})_{j+1}
(\varphi^G_{\textbf{V}})_j}{(\varphi^G_{\textbf{V}})_j}\\
& &+\frac{a^{\Lm_{j-1}+\Lm_{j}+\Lm_{k-1}+\Lm_r}Y_{1,j}Y_{2,j}Y_{2,k}Y^2_{1,r-1}}{(\varphi^G_{\textbf{V}})_j Y_{2,r}}
(1+A^{-1}_{1,r-1}+\cdots+A^{-1}_{1,r-1}\cdots A^{-1}_{1,k+1})\\
& &\qquad \qquad \ \ \times A^{-1}_{1,r-1}\cdots A^{-1}_{1,j+2}A^{-1}_{1,j+1}(1+A^{-1}_{1,j}+\cdots+A^{-1}_{1,j}\cdots A^{-1}_{1,1})\\
&=&a^{\Lm_{j-1}-\Lm_j}\frac{Y_{2,j}}{Y_{2,j+1}}
(\varphi^G_{\mu_{j+1}\cdots\mu_{k-2}\mu_{k-1}\mu_{r}\mu_{r-1}\cdots\mu_{k+1}\mu_k(\textbf{V})})_{j+1}
\\
&+&a^{\Lm_{j-1}+\Lm_{k-1}+\Lm_r}\frac{Y_{2,j}Y_{2,k}Y^2_{1,r-1}}{Y_{2,r}}
(1+A^{-1}_{1,r-1}+\cdots+A^{-1}_{1,r-1}\cdots A^{-1}_{1,k+1})
A^{-1}_{1,r-1}\cdots A^{-1}_{1,j+2}A^{-1}_{1,j+1}\\
&=&a^{\Lm_{j-1}+\Lm_{k-1}+\Lm_r} Y_{2,j}Y_{2,k}Y_{1,r}
+a^{\Lm_{j-1}+\Lm_{k-1}+\Lm_r}\frac{Y_{2,j}Y_{2,k}Y^2_{1,r-1}}{Y_{2,r}}\\
&\times&(1+A^{-1}_{1,r-1}+\cdots+A^{-1}_{1,r-1}\cdots A^{-1}_{1,k+1})
(1+A^{-1}_{1,r-1}+\cdots+A^{-1}_{1,r-1}\cdots A^{-1}_{1,j+1}),
\end{eqnarray*}
where we use (\ref{iniC-exp1}) in the second equality, and use (\ref{C-lem3-pr3}) in the third equality. Thus, we get (\ref{C-lem3-exp}), and one can rewrite it as follows:
\begin{eqnarray*}
& &\hspace{-20pt}(\varphi^G_{\mu_j\mu_{j+1}\cdots\mu_{k-2}\mu_{k-1}\mu_{r}\mu_{r-1}\cdots\mu_{k+1}\mu_k(\textbf{V})})_j(a;\textbf{Y})\\
&=&a^{\Lm_{j-1}+\Lm_{k-1}+\Lm_r} Y_{2,j}Y_{2,k}Y_{1,r}(1+A^{-1}_{1,r}
(1+A^{-1}_{1,r-1}+\cdots+A^{-1}_{1,r-1}\cdots A^{-1}_{1,k+1})^2)
\\
&+&a^{\Lm_{j-1}+\Lm_{k-1}+\Lm_r} Y_{2,j}Y_{1,k-1}Y_{1,r-1}(1+A^{-1}_{1,r-1}+\cdots+A^{-1}_{1,r-1}\cdots A^{-1}_{1,k+1})\\
& & \qquad \qquad \qquad \qquad \qquad \qquad\times(1+A^{-1}_{1,k-1}+\cdots+A^{-1}_{1,k-1}\cdots A^{-1}_{1,j+1})\\
&=& a^{\Lm_{j-1}+\Lm_{k-1}+\Lm_{r}} \sum_{b\in B(\Lm_j+\Lambda_k+2\Lambda_{r})_{s_{k+1}s_{k+2}\cdots s_{r}} } C(b)\mu(b)\\
& &+a^{\Lm_{j-1}+\Lm_{k-1}+\Lm_{r}}\sum_{b\in B(\Lm_j+\Lambda_{k-1}+ \Lambda_{r-1})_{s_{j+1}\cdots s_{k-1}s_{k+1}s_{k+2}\cdots s_{r-1}}} \mu'(b),
\end{eqnarray*}
where each coefficient $C(b)$ is either $1$ or $2$. Thus, we obtain
(\ref{C-lem3-cry}) in the same way as (\ref{B-lem2pr-c}). \qed 

\subsection{The proof of Theorem \ref{thmD}}

In this final subsection, we shall prove Theorem \ref{thmD}. Let $G={\rm SO}_{2r}$. The mutation diagram of the initial seed $\Sigma_0$ is
\begin{equation}\label{inidiaD-2}
\begin{xy}
(90,90)*{\cdots}="emp1",
(80,105)*{(\varphi_{\textbf{V}})_{r+k+1}}="7",
(80,90) *{(\varphi_{\textbf{V}})_{k+1}}="3",
(80,80)*{(\varphi_{\textbf{V}})_{-(k+1)}}="-4",
(55,105)*{(\varphi_{\textbf{V}})_{r+k}}="r+k",
(55,90) *{(\varphi_{\textbf{V}})_{k}}="k",
(55,80)*{(\varphi_{\textbf{V}})_{-k}}="-j_k",
(30,105)*{(\varphi_{\textbf{V}})_{r+1}}="6",
(30,90) *{(\varphi_{\textbf{V}})_{1}}="2",
(30,80)*{(\varphi_{\textbf{V}})_{-1}}="-1",
(45,90)*{\ \ \cdots}="emp",
(45,105)*{}="emp2",
(45,80)*{}="emp3",
(140,105)*{(\varphi_{\textbf{V}})_{2r}}="2rr",
(140,96) *{(\varphi_{\textbf{V}})_{r}}="rr",
(140,80)*{(\varphi_{\textbf{V}})_{-r}}="-j_rr",
(120,105)*{(\varphi_{\textbf{V}})_{2r-1}}="2r",
(120,90) *{(\varphi_{\textbf{V}})_{r-1}}="r",
(120,80)*{(\varphi_{\textbf{V}})_{-(r-1)}}="-j_r",
(100,105)*{(\varphi_{\textbf{V}})_{2r-2}}="2r-1",
(100,90) *{(\varphi_{\textbf{V}})_{r-2}}="r-1",
(100,80)*{(\varphi_{\textbf{V}})_{-(r-2)}}="-r+1",
\ar@{->} "7";"3"
\ar@{->} "r+k";"k"
\ar@{->} "6";"2"
\ar@{->} "3";"-4"
\ar@{->} "k";"-j_k"
\ar@{->} "2";"-1"
\ar@{->} "k";"3"
\ar@{->} "2";"emp"
\ar@{->} "3";"r+k"
\ar@{->} "emp";"6"
\ar@{->} "-4";"k"
\ar@{->} "emp3";"2"
\ar@{->} "2r";"r"
\ar@{->} "2r-1";"r-1"
\ar@{->} "r";"-j_r"
\ar@{->} "r";"2r-1"
\ar@{->} "rr";"2r-1"
\ar@{->} "r-1";"-r+1"
\ar@{->} "r-1";"r"
\ar@{->} "r-1";"rr"
\ar@{->} "2rr";"rr"
\ar@{->} "rr";"-j_rr"
\ar@{->} "-j_r";"r-1"
\ar@{->} "-j_rr";"r-1"
\end{xy}
\end{equation}
First, we can prove the following proposition in a similar way to Proposition \ref{iniB} and \ref{iniC}.
\begin{prop}\label{iniD}
\begin{enumerate}
\item For $k\in[1,r]$,
the initial cluster variables $(\varphi^G_{\textbf{V}})_k(a;\textbf{Y})$ in $\mathbb{C}[G^{e,c^2}]$ are
described as
\[ 
(\varphi^G_{\textbf{V}})_k(a;\textbf{Y})=a^{\Lambda_k} \sum_{b\in B(\Lambda_k)_{c^2_{>2r-k}}} \mu(b),
\]
where $\mu:B(\Lambda_k)\rightarrow \mathcal{Y}$ is the monomial realization of $B(\Lambda_k)$ in Theorem \ref{thmD} $(i)$.
\item
For $k\in[1,r]$, the frozen cluster variables $(\varphi^G_{\textbf{V}})_{-k}(a;\textbf{Y})$, $(\varphi^G_{\textbf{V}})_{r+k}(a;\textbf{Y})$ in $\mathbb{C}[G^{e,c^2}]$ are described as
\[
(\varphi^G_{\textbf{V}})_{-k}(a;\textbf{Y})=a^{\Lm_k}Y_{1,k}Y_{2,k},
\quad
(\varphi^G_{\textbf{V}})_{r+k}(a;\textbf{Y})=a^{\Lm_k}.
\]
\end{enumerate}
\end{prop}

For $k\in[1,r-1]$, just as in (\ref{iniBexp2}), (\ref{iniBexp-r}) and (\ref{iniC-exp1}), we see that
\begin{equation}\label{iniD-exp1}
(\varphi^G_{\textbf{V}})_{k}(a;\textbf{Y})=a^{\Lm_k}Y_{1,k}(1+A^{-1}_{1,k}+A^{-1}_{1,k}A^{-1}_{1,k-1}+\cdots+A^{-1}_{1,k}\cdots A^{-1}_{1,1}),
\end{equation}
and
\begin{equation}\label{iniD-exp2}
(\varphi^G_{\textbf{V}})_{r}(a;\textbf{Y})=a^{\Lm_r}Y_{1,r}(1+A^{-1}_{1,r}+A^{-1}_{1,r}A^{-1}_{1,r-2}+\cdots+A^{-1}_{1,r}A^{-1}_{1,r-2}A^{-1}_{1,r-3}\cdots A^{-1}_{1,1}).
\end{equation}
Considering (\ref{iniD-exp1}), we see that for $k\in[1,r-2]$,
\begin{equation}\label{twovarsD}
a^{-\Lambda_k}(\varphi^G_{\textbf{V}})_k(a;\textbf{Y})=Y_{1,k}+\frac{Y_{2,k+1}}{Y_{2,k}}(a^{-\Lm_{k-1}}(\varphi^G_{\textbf{V}})_{k-1}(a;\textbf{Y})).
\end{equation}

\begin{lem}\label{D-lem1}
\begin{enumerate}
\item
\[ 
(\varphi^G_{\mu_{r-1}(\textbf{V})})_{r-1}(a;\textbf{Y})=
a^{\Lambda_{r-2}} Y_{2,r-1},\quad
(\varphi^G_{\mu_{r}\mu_{r-1}(\textbf{V})})_{r}(a;\textbf{Y})=
a^{\Lambda_{r-2}} Y_{2,r}.
\]
\item For $k\in[1,r-2]$, the cluster variables $(\varphi^G_{\mu_k\mu_{k+1}\cdots\mu_{r-3}\mu_{r-2}\mu_r\mu_{r-1}(\textbf{V})})_k(a;\textbf{Y})$ in $\mathbb{C}[G^{e,c^2}]$ are
described as
\[ 
(\varphi^G_{\mu_k\mu_{k+1}\cdots\mu_{r-3}\mu_{r-2}\mu_r\mu_{r-1}(\textbf{V})})_k(a;\textbf{Y})=
a^{\Lambda_{k-1}+\Lambda_{r-2}} Y_{2,k}.
\]
\end{enumerate}
\end{lem}
\nd
[{\sl Proof.}]

(i) The mutation diagram (\ref{inidiaD-2}) and the equation (\ref{iniD-exp1}) say that
\begin{eqnarray*}
(\varphi^G_{\mu_{r-1}(\textbf{V})})_{r-1}&=&\frac{(\varphi^G_{\textbf{V}})_{2r-2}(\varphi^G_{\textbf{V}})_{-(r-1)}+(\varphi^G_{\textbf{V}})_{r-2}(\varphi^G_{\textbf{V}})_{2r-1}}{(\varphi^G_{\textbf{V}})_{r-1}}\\
&=&a^{\Lm_{r-1}+\Lm_{r-2}}\frac{Y_{1,r-1}Y_{2,r-1}+Y_{1,r-2}
(1+A^{-1}_{1,r-2}+\cdots+A^{-1}_{1,r-2}\cdots A^{-1}_{1,1})}{a^{\Lm_{r-1}}Y_{1,r-1}(1+A^{-1}_{1,r-1}+\cdots+A^{-1}_{1,r-1}\cdots A^{-1}_{1,1})}\\
&=&a^{\Lambda_{r-2}} Y_{2,r-1}.
\end{eqnarray*}

Next, the mutation diagram of $\mu_{r-1}(\Sigma_0)$ is as follows by Lemma \ref{mutgamlem}:
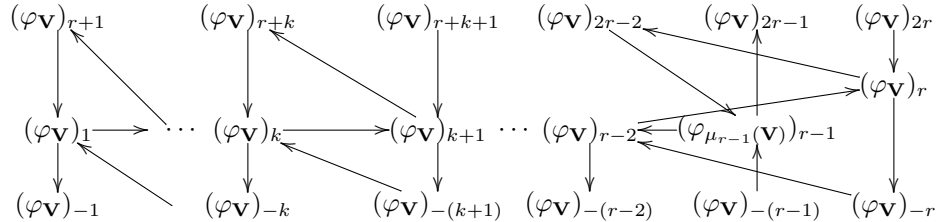
\[
\begin{xy}
(90,90)*{\cdots}="emp1",
(80,105)*{(\varphi_{\textbf{V}})_{r+k+1}}="7",
(80,90) *{(\varphi_{\textbf{V}})_{k+1}}="3",
(80,80)*{(\varphi_{\textbf{V}})_{-(k+1)}}="-4",
(55,105)*{(\varphi_{\textbf{V}})_{r+k}}="r+k",
(55,90) *{(\varphi_{\textbf{V}})_{k}}="k",
(55,80)*{(\varphi_{\textbf{V}})_{-k}}="-j_k",
(30,105)*{(\varphi_{\textbf{V}})_{r+1}}="6",
(30,90) *{(\varphi_{\textbf{V}})_{1}}="2",
(30,80)*{(\varphi_{\textbf{V}})_{-1}}="-1",
(45,90)*{\ \ \cdots}="emp",
(45,105)*{}="emp2",
(45,80)*{}="emp3",
(140,105)*{(\varphi_{\textbf{V}})_{2r}}="2rr",
(140,96) *{(\varphi_{\textbf{V}})_{r}}="rr",
(140,80)*{(\varphi_{\textbf{V}})_{-r}}="-j_rr",
(122,105)*{(\varphi_{\textbf{V}})_{2r-1}}="2r",
(122,90) *{(\varphi_{\mu_{r-1}(\textbf{V})})_{r-1}}="r",
(122,80)*{(\varphi_{\textbf{V}})_{-(r-1)}}="-j_r",
(100,105)*{(\varphi_{\textbf{V}})_{2r-2}}="2r-1",
(100,90) *{(\varphi_{\textbf{V}})_{r-2}}="r-1",
(100,80)*{(\varphi_{\textbf{V}})_{-(r-2)}}="-r+1",
\ar@{->} "7";"3"
\ar@{->} "r+k";"k"
\ar@{->} "6";"2"
\ar@{->} "3";"-4"
\ar@{->} "k";"-j_k"
\ar@{->} "2";"-1"
\ar@{->} "k";"3"
\ar@{->} "2";"emp"
\ar@{->} "3";"r+k"
\ar@{->} "emp";"6"
\ar@{->} "-4";"k"
\ar@{->} "emp3";"2"
\ar@{->} "r";"2r"
\ar@{->} "-j_r";"r"
\ar@{->} "2r-1";"r"
\ar@{->} "rr";"2r-1"
\ar@{->} "r-1";"-r+1"
\ar@{->} "r";"r-1"
\ar@{->} "r-1";"rr"
\ar@{->} "2rr";"rr"
\ar@{->} "rr";"-j_rr"
\ar@{->} "-j_rr";"r-1"
\end{xy}
\]
Hence, by a calculation similar to the one of $(\varphi^G_{\mu_{r-1}(\textbf{V})})_{r-1}$, we obtain $(\varphi^G_{\mu_{r}\mu_{r-1}(\textbf{V})})_{r}(a;\textbf{Y})=
a^{\Lambda_{r-2}} Y_{2,r}$.

(ii) We use the induction on $(r-2-k)$. First, let $r-2-k=0$, so that $k=r-2$. The mutation diagram of $\mu_r\mu_{r-1}(\Sigma_0)$ is
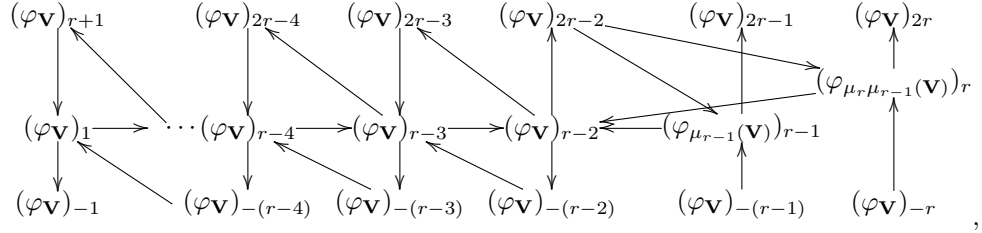
\[
\begin{xy}
(75,105)*{(\varphi_{\textbf{V}})_{2r-3}}="7",
(75,90) *{(\varphi_{\textbf{V}})_{r-3}}="3",
(75,80)*{(\varphi_{\textbf{V}})_{-(r-3)}}="-4",
(55,105)*{(\varphi_{\textbf{V}})_{2r-4}}="r+k",
(55,90) *{(\varphi_{\textbf{V}})_{r-4}}="k",
(55,80)*{(\varphi_{\textbf{V}})_{-(r-4)}}="-j_k",
(30,105)*{(\varphi_{\textbf{V}})_{r+1}}="6",
(30,90) *{(\varphi_{\textbf{V}})_{1}}="2",
(30,80)*{(\varphi_{\textbf{V}})_{-1}}="-1",
(45,90)*{\ \ \cdots}="emp",
(45,105)*{}="emp2",
(45,80)*{}="emp3",
(140,105)*{(\varphi_{\textbf{V}})_{2r}}="2rr",
(140,96) *{(\varphi_{\mu_r\mu_{r-1}(\textbf{V})})_{r}}="rr",
(140,80)*{(\varphi_{\textbf{V}})_{-r}}="-j_rr",
(120,105)*{(\varphi_{\textbf{V}})_{2r-1}}="2r",
(120,90) *{(\varphi_{\mu_{r-1}(\textbf{V})})_{r-1}}="r",
(120,80)*{(\varphi_{\textbf{V}})_{-(r-1)}}="-j_r",
(95,105)*{(\varphi_{\textbf{V}})_{2r-2}}="2r-1",
(95,90) *{(\varphi_{\textbf{V}})_{r-2}}="r-1",
(95,80)*{(\varphi_{\textbf{V}})_{-(r-2)}}="-r+1",
\ar@{->} "7";"3"
\ar@{->} "r+k";"k"
\ar@{->} "6";"2"
\ar@{->} "3";"-4"
\ar@{->} "3";"r-1"
\ar@{->} "k";"-j_k"
\ar@{->} "2";"-1"
\ar@{->} "k";"3"
\ar@{->} "2";"emp"
\ar@{->} "3";"r+k"
\ar@{->} "emp";"6"
\ar@{->} "-4";"k"
\ar@{->} "emp3";"2"
\ar@{->} "r";"2r"
\ar@{->} "-j_r";"r"
\ar@{->} "2r-1";"r"
\ar@{->} "2r-1";"rr"
\ar@{->} "r-1";"-r+1"
\ar@{->} "r-1";"2r-1"
\ar@{->} "r-1";"7"
\ar@{->} "-r+1";"3"
\ar@{->} "r";"r-1"
\ar@{->} "rr";"r-1"
\ar@{->} "rr";"2rr"
\ar@{->} "-j_rr";"rr"
\end{xy},
\]
which yields
\begin{eqnarray*}
& &\hspace{-20pt}(\varphi^G_{\mu_{r-2}\mu_r\mu_{r-1}(\textbf{V})})_{r-2}\\
&=&\frac{(\varphi^G_{\mu_{r}\mu_{r-1}(\textbf{V})})_{r}(\varphi^G_{\mu_{r-1}(\textbf{V})})_{r-1}(\varphi^G_{\textbf{V}})_{r-3}+
(\varphi^G_{\textbf{V}})_{2r-3}(\varphi^G_{\textbf{V}})_{2r-2}
(\varphi^G_{\textbf{V}})_{-(r-2)}}{(\varphi^G_{\textbf{V}})_{r-2}}\\
&=&\frac{a^{2\Lm_{r-2}+\Lm_{r-3}}Y_{2,r-1}Y_{2,r}Y_{1,r-3}(1+A^{-1}_{1,r-3}+\cdots+A^{-1}_{1,r-3}\cdots A^{-1}_{1,1})
+a^{2\Lm_{r-2}+\Lm_{r-3}}Y_{1,r-2}Y_{2,r-2}}{a^{\Lm_{r-2}}Y_{1,r-2}(1+A^{-1}_{1,r-2}+A^{-1}_{1,r-2}A^{-1}_{1,r-3}+\cdots+A^{-1}_{1,r-2}\cdots A^{-1}_{1,1})}\\
&=&a^{\Lm_{r-2}+\Lm_{r-3}}Y_{2,r-2},
\end{eqnarray*}
where we use (\ref{iniD-exp1}) in the second equality. 

Next, let $r-2-k>0$. The vertices and arrows around the vertex $(\varphi^G_{\textbf{V}})_{k}$ in the mutation diagram of $\mu_{k+1}\mu_{k+2}\cdots\mu_{r-3}\mu_{r-2}\mu_r\mu_{r-1}(\Sigma_0)$ are
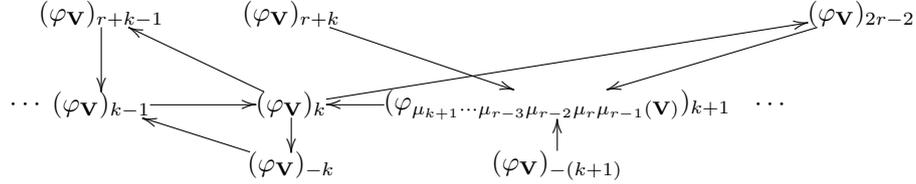
\[
\begin{xy}
(110,102) *{(\varphi_{\textbf{V}})_{2r-2}}="r-1",
(98,90)*{\cdots}="emp1",
(70,90)*{ (\varphi_{\mu_{k+1}\cdots\mu_{r-3}\mu_{r-2}\mu_r\mu_{r-1}(\textbf{V})})_{k+1}}="k",
(70,82)*{(\varphi_{\textbf{V}})_{-(k+1)}}="-j_k",
(35,102)*{(\varphi_{\textbf{V}})_{r+k}}="8",
(35,90) *{(\varphi_{\textbf{V}})_{k}}="4",
(35,82)*{(\varphi_{\textbf{V}})_{-k}}="-2",
(10,102)*{(\varphi_{\textbf{V}})_{r+k-1}}="6",
(10,90) *{(\varphi_{\textbf{V}})_{k-1}}="2",
(0,90)*{\cdots}="emp",
\ar@{->} "4";"r-1"
\ar@{->} "k";"4"
\ar@{->} "r-1";"k"
\ar@{->} "-j_k";"k"
\ar@{->} "8";"k"
\ar@{->} "4";"-2"
\ar@{->} "2";"4"
\ar@{->} "4";"6"
\ar@{->} "6";"2"
\ar@{->} "-2";"2"
\end{xy}
\]
Hence,
\begin{eqnarray*}
& &\hspace{-20pt}(\varphi^G_{\mu_k\mu_{k+1}\mu_{k+2}\cdots\mu_{r-3}\mu_{r-2}\mu_r\mu_{r-1}(\textbf{V})})_{k}\\
&=&\frac{(\varphi^G_{\mu_{k+1}\mu_{k+2}\cdots\mu_{r-3}\mu_{r-2}\mu_r\mu_{r-1}(\textbf{V})})_{k+1}(\varphi^G_{\textbf{V}})_{k-1}+
(\varphi^G_{\textbf{V}})_{r+k-1}(\varphi^G_{\textbf{V}})_{2r-2}
(\varphi^G_{\textbf{V}})_{-k}}{(\varphi^G_{\textbf{V}})_{k}}\\
&=&\frac{a^{\Lm_{k-1}+\Lambda_{k}+\Lambda_{r-2}}(Y_{2,k+1}Y_{1,k-1}
(1+A^{-1}_{1,k-1}+\cdots+A^{-1}_{1,k-1}\cdots A^{-1}_{1,1})+Y_{1,k}Y_{2,k})
}{a^{\Lm_k}Y_{1,k}
(1+A^{-1}_{1,k}+\cdots+A^{-1}_{1,k}\cdots A^{-1}_{1,1})}\\
&=&a^{\Lm_{k-1}+\Lambda_{r-2}} Y_{2,k},
\end{eqnarray*}
where we use the induction hypothesis and (\ref{iniD-exp1}) in the second equality. \qed

\vspace{3mm}

\nd
[{\sl Proof of Theorem \ref{thmD} (ii) and (iii).}]

The claim (ii) is obtained by the same calculation as in Theorem \ref{thmB} (ii). So let us consider the claim (iii).

First, by Lemma \ref{mutgamlem}, the mutation diagram of $\mu_{r-3}\cdots\mu_{k+1}\mu_{k}(\Sigma_0)$ is
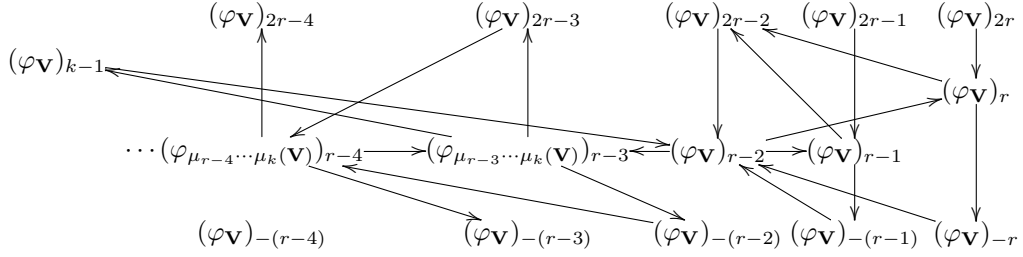
\[
\begin{xy}
(94,108)*{(\varphi_{\textbf{V}})_{2r}}="2rr",
(94,98) *{(\varphi_{\textbf{V}})_{r}}="rr",
(94,79)*{(\varphi_{\textbf{V}})_{-r}}="-rr",
(78,108)*{(\varphi_{\textbf{V}})_{2r-1}}="7",
(78,90) *{(\varphi_{\textbf{V}})_{r-1}}="3",
(78,79)*{(\varphi_{\textbf{V}})_{-(r-1)}}="-4",
(60,108)*{(\varphi_{\textbf{V}})_{2r-2}}="r+k",
(60,90) *{(\varphi_{\textbf{V}})_{r-2}}="k",
(60,79)*{(\varphi_{\textbf{V}})_{-(r-2)}}="-j_k",
(35,108)*{(\varphi_{\textbf{V}})_{2r-3}}="8",
(35,90) *{(\varphi_{\mu_{r-3}\cdots\mu_k(\textbf{V})})_{r-3}}="4",
(35,79)*{(\varphi_{\textbf{V}})_{-(r-3)}}="-2",
(0,108)*{(\varphi_{\textbf{V}})_{2r-4}}="6",
(0,90) *{(\varphi_{\mu_{r-4}\cdots\mu_k(\textbf{V})})_{r-4}}="2",
(0,79)*{(\varphi_{\textbf{V}})_{-(r-4)}}="-1",
(-16,90)*{\cdots}="emp",
(-27,102)*{(\varphi_{\textbf{V}})_{k-1}}="k-1",
\ar@{->} "7";"3"
\ar@{->} "r+k";"k"
\ar@{->} "4";"8"
\ar@{->} "2";"6"
\ar@{->} "3";"-4"
\ar@{->} "k";"3"
\ar@{->} "k";"4"
\ar@{->} "4";"-j_k"
\ar@{->} "2";"4"
\ar@{->} "3";"r+k"
\ar@{->} "8";"2"
\ar@{->} "-4";"k"
\ar@{->} "k-1";"k"
\ar@{->} "4";"k-1"
\ar@{->} "-j_k";"2"
\ar@{->} "2";"-2"
\ar@{->} "rr";"-rr"
\ar@{->} "2rr";"rr"
\ar@{->} "k";"rr"
\ar@{->} "rr";"r+k"
\ar@{->} "-rr";"k"
\end{xy}
\]
From this diagram, 
\begin{eqnarray}
& &\hspace{-20pt}(\varphi^G_{\mu_{r-2}\cdots\mu_{k+1}\mu_k(\textbf{V})})_{r-2}(a;\textbf{Y})\nonumber\\
&=&\frac{(\varphi^G_{\textbf{V}})_{k-1}(\varphi^G_{\textbf{V}})_{2r-2}(\varphi^G_{\textbf{V}})_{-(r-1)}(\varphi^G_{\textbf{V}})_{-r}+(\varphi^G_{\mu_{r-3}\cdots\mu_{k+1}\mu_{k}(\textbf{V})})_{r-3}(\varphi^G_{\textbf{V}})_{r-1}(\varphi^G_{\textbf{V}})_{r}}{(\varphi^G_{\textbf{V}})_{r-2}}.\qquad \quad\label{D-lem2-pr-1}
\end{eqnarray}

Using (\ref{D-exp1}) for $l=r-3$, we obtain
\begin{eqnarray*}
& &\hspace{-20pt}(\varphi^G_{\mu_{r-3}\cdots\mu_{k+1}\mu_{k}(\textbf{V})})_{r-3}\\
&=&a^{\Lm_{k-1}+\Lm_{r-2}} Y_{2,k}Y_{1,r-2} (1+A^{-1}_{1,r-2}+A^{-1}_{1,r-2}A^{-1}_{1,r-3}+\cdots+A^{-1}_{1,r-2}A^{-1}_{1,r-3}\cdots A^{-1}_{1,k+1})\\
&=&a^{\Lm_{k-1}}Y_{2,k}(\varphi^G_{\textbf{V}})_{r-2}-a^{\Lm_{r-2}}Y_{2,r-1}Y_{2,r}(\varphi^G_{\textbf{V}})_{k-1},
\end{eqnarray*}
where we use (\ref{iniD-exp1}) in the second equality. Similarly, using (\ref{iniD-exp1}) and (\ref{iniD-exp2}), we also get
\[
(\varphi^G_{\textbf{V}})_{r-1}=a^{\Lm_{r-1}}Y_{1,r-1}+a^{\Lm_{r-1}-\Lm_{r-2}}\frac{(\varphi^G_{\textbf{V}})_{r-2}}{Y_{2,r-1}},\ \ (\varphi^G_{\textbf{V}})_{r}=a^{\Lm_{r}}Y_{1,r}+a^{\Lm_{r}-\Lm_{r-2}}\frac{(\varphi^G_{\textbf{V}})_{r-2}}{Y_{2,r}}.
\]
Applying these formulas to (\ref{D-lem2-pr-1}), it follows
\begin{eqnarray}
& &\hspace{-20pt}(\varphi^G_{\mu_{r-2}\cdots\mu_{k+1}\mu_k(\textbf{V})})_{r-2}(a;\textbf{Y})\nonumber\\
&=&a^{\Lm_{k-1}}Y_{2,k}(\varphi^G_{\textbf{V}})_{r-1}(\varphi^G_{\textbf{V}})_{r}\nonumber\\
& &-a^{\Lm_{r-1}+\Lm_r}Y_{2,r-1}Y_{2,r}(\varphi^G_{\textbf{V}})_{k-1}(\frac{Y_{1,r-1}}{Y_{2,r}}+\frac{Y_{1,r}}{Y_{2,r-1}}+a^{-\Lm_{r-2}}\frac{(\varphi^G_{\textbf{V}})_{r-2}}{Y_{2,r-1}Y_{2,r}})\label{D-lem2-pr1-1}\\
&=&a^{\Lm_{k-1}+\Lm_{r-1}+\Lm_r}Y_{2,k}Y_{1,r-1}Y_{1,r}
(1+A^{-1}_{1,r-1}+\cdots+A^{-1}_{1,r-1}\cdots A^{-1}_{1,1})\nonumber\\
& & \times (1+A^{-1}_{1,r}+A^{-1}_{1,r}A^{-1}_{1,r-2}+ \cdots+A^{-1}_{1,r}A^{-1}_{1,r-2}A^{-1}_{1,r-3}\cdots A^{-1}_{1,1})\nonumber\\
& &-a^{\Lm_{k-1}+\Lm_{r-1}+\Lm_r}Y_{1,k-1}Y_{2,r-1}Y_{2,r}
(1+A^{-1}_{1,k-1}+\cdots+A^{-1}_{1,k-1}\cdots A^{-1}_{1,1})\nonumber\\
& & \times (\frac{Y_{1,r}}{Y_{2,r-1}}+\frac{Y_{1,r-1}}{Y_{2,r}}(1+A^{-1}_{1,r-1}+A^{-1}_{1,r-1}A^{-1}_{1,r-2}+ \cdots+A^{-1}_{1,r-1}A^{-1}_{1,r-2}\cdots A^{-1}_{1,1}))\nonumber\\
&=&a^{\Lm_{k-1}+\Lm_{r-1}+\Lm_r}Y_{2,k}Y_{1,r-1}Y_{1,r}
(1+A^{-1}_{1,r-1}+\cdots+A^{-1}_{1,r-1}\cdots A^{-1}_{1,1})\nonumber\\
& & \times (1+A^{-1}_{1,r}+A^{-1}_{1,r}A^{-1}_{1,r-2}+ \cdots+A^{-1}_{1,r}A^{-1}_{1,r-2}A^{-1}_{1,r-3}\cdots A^{-1}_{1,k+1})\label{D-lem2-pr2}\\
& &-a^{\Lm_{k-1}+\Lm_{r-1}+\Lm_r}Y_{1,k-1}Y_{1,r}Y_{2,r}
(1+A^{-1}_{1,k-1}+\cdots+A^{-1}_{1,k-1}\cdots A^{-1}_{1,1})\nonumber\\
&=&a^{\Lm_{k-1}+\Lm_{r-1}+\Lm_r}Y_{2,k}Y_{1,r-1}Y_{1,r}
(1+A^{-1}_{1,r-1}+\cdots+A^{-1}_{1,r-1}\cdots A^{-1}_{1,k+1})\nonumber\\
& &\times (1+A^{-1}_{1,r}+A^{-1}_{1,r}A^{-1}_{1,r-2}+ \cdots+A^{-1}_{1,r}A^{-1}_{1,r-2}A^{-1}_{1,r-3}\cdots A^{-1}_{1,k+1})\nonumber\\
& &+a^{\Lm_{k-1}+\Lm_{r-1}+\Lm_r}Y_{1,k-1}Y_{1,r-2}
(1+A^{-1}_{1,r-2}+\cdots+A^{-1}_{1,r-2}\cdots A^{-1}_{1,k+1})\nonumber\\
& &\times (1+A^{-1}_{1,k-1}+A^{-1}_{1,k-1}A^{-1}_{1,k-2}+\cdots+A^{-1}_{1,k-1}A^{-1}_{1,k-2}\cdots A^{-1}_{1,1}).\nonumber
\end{eqnarray}
By this explicit formula, the conclusion (\ref{D-lem2-cry}) follows.

\vspace{3mm}

\nd
[{\sl Proof of Theorem \ref{thmD} (iv)}]

The vertices and arrows around $(\varphi_{\textbf{V}})_{r-1}$ in the mutation diagram of $\mu_{r-2}\mu_{r-3}\cdots\mu_{k+1}\mu_{k}(\Sigma_0)$ are
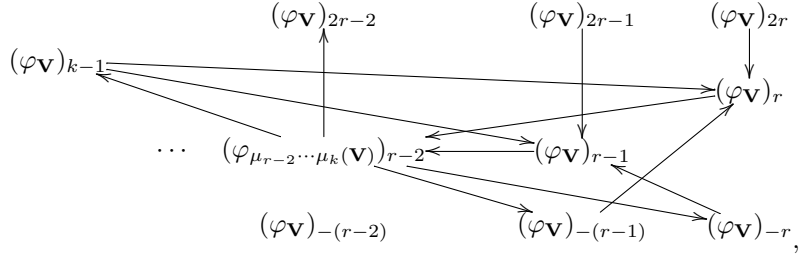
\[
\begin{xy}
(106,108)*{(\varphi_{\textbf{V}})_{2r}}="2rr",
(106,98) *{(\varphi_{\textbf{V}})_{r}}="rr",
(106,80)*{(\varphi_{\textbf{V}})_{-r}}="-rr",
(84,108)*{(\varphi_{\textbf{V}})_{2r-1}}="7",
(84,90) *{(\varphi_{\textbf{V}})_{r-1}}="3",
(84,80)*{(\varphi_{\textbf{V}})_{-(r-1)}}="-4",
(50,108)*{(\varphi_{\textbf{V}})_{2r-2}}="r+k",
(50,90) *{(\varphi_{\mu_{r-2}\cdots\mu_k(\textbf{V})})_{r-2}}="k",
(50,80)*{(\varphi_{\textbf{V}})_{-(r-2)}}="-j_k",
(30,90)*{\cdots}="emp",
(15,102)*{(\varphi_{\textbf{V}})_{k-1}}="k-1",
\ar@{->} "7";"3"
\ar@{->} "k";"r+k"
\ar@{->} "3";"k"
\ar@{->} "k";"-4"
\ar@{->} "k";"k-1"
\ar@{->} "2rr";"rr"
\ar@{->} "rr";"k"
\ar@{->} "k";"-rr"
\ar@{->} "k-1";"rr"
\ar@{->} "k-1";"3"
\ar@{->} "-rr";"3"
\ar@{->} "-4";"rr"
\end{xy},
\]
which yields
\begin{equation}\label{D-lem2-pr3}
(\varphi^G_{\mu_{r-1}\mu_{r-2}\cdots\mu_{k+1}\mu_k(\textbf{V})})_{r-1}(a;\textbf{Y})
=\frac{(\varphi^G_{\mu_{r-2}\cdots\mu_{k+1}\mu_k(\textbf{V})})_{r-2}+(\varphi^G_{\textbf{V}})_{k-1}(\varphi^G_{\textbf{V}})_{2r-1}(\varphi^G_{\textbf{V}})_{-r}}{(\varphi^G_{\textbf{V}})_{r-1}}.
\end{equation}

It follows from (\ref{iniD-exp1}) and (\ref{D-lem2-pr2}) that
\begin{eqnarray}
& &\hspace{-20pt}(\varphi^G_{\mu_{r-2}\cdots\mu_{k+1}\mu_k(\textbf{V})})_{r-2}(a;\textbf{Y})\nonumber\\
&=&a^{\Lm_{k-1}+\Lm_r}Y_{2,k}Y_{1,r}(\varphi^G_{\textbf{V}})_{r-1}
(1+A^{-1}_{1,r}+A^{-1}_{1,r}A^{-1}_{1,r-2}+ \cdots+A^{-1}_{1,r}A^{-1}_{1,r-2}A^{-1}_{1,r-3}\cdots A^{-1}_{1,k+1})\qquad \quad\label{D-lem2-pr4}\\
& &-a^{\Lm_{r-1}+\Lm_r}Y_{1,r}Y_{2,r}
(\varphi^G_{\textbf{V}})_{k-1}\nonumber.
\end{eqnarray}
Substituting this for (\ref{D-lem2-pr3}), we see that
\begin{eqnarray}
& &\hspace{-20pt}(\varphi^G_{\mu_{r-1}\mu_{r-2}\cdots\mu_{k+1}\mu_k(\textbf{V})})_{r-1}(a;\textbf{Y})\nonumber\\
&=&a^{\Lm_{k-1}+\Lm_r}Y_{2,k}Y_{1,r}
(1+A^{-1}_{1,r}+A^{-1}_{1,r}A^{-1}_{1,r-2}+ \cdots+A^{-1}_{1,r}A^{-1}_{1,r-2}A^{-1}_{1,r-3}\cdots A^{-1}_{1,k+1}),\qquad \quad \nonumber
\end{eqnarray}
which implies (\ref{D-lem3-cry}).

Finally, we consider the cluster variable $(\varphi^G_{\mu_{r}\mu_{r-1}\cdots\mu_{k+1}\mu_k(\textbf{V})})_{r}(a;\textbf{Y})$. The mutation diagram of $\mu_{r-1}\mu_{r-2}\mu_{r-3}\cdots\mu_{k+1}\mu_{k}(\Sigma_0)$
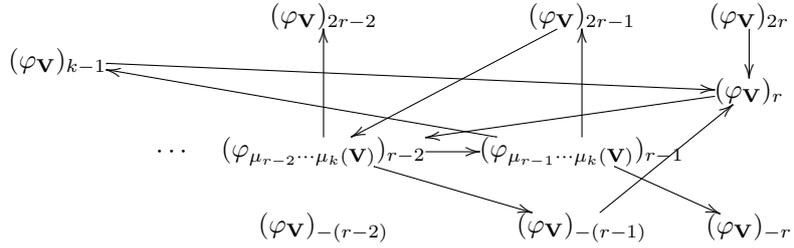
\[
\begin{xy}
(106,108)*{(\varphi_{\textbf{V}})_{2r}}="2rr",
(106,98) *{(\varphi_{\textbf{V}})_{r}}="rr",
(106,80)*{(\varphi_{\textbf{V}})_{-r}}="-rr",
(84,108)*{(\varphi_{\textbf{V}})_{2r-1}}="7",
(84,90) *{(\varphi_{\mu_{r-1}\cdots\mu_k(\textbf{V})})_{r-1}}="3",
(84,80)*{(\varphi_{\textbf{V}})_{-(r-1)}}="-4",
(50,108)*{(\varphi_{\textbf{V}})_{2r-2}}="r+k",
(50,90) *{(\varphi_{\mu_{r-2}\cdots\mu_k(\textbf{V})})_{r-2}}="k",
(50,80)*{(\varphi_{\textbf{V}})_{-(r-2)}}="-j_k",
(30,90)*{\cdots}="emp",
(15,102)*{(\varphi_{\textbf{V}})_{k-1}}="k-1",
\ar@{->} "3";"7"
\ar@{->} "7";"k"
\ar@{->} "k";"r+k"
\ar@{->} "k";"3"
\ar@{->} "k";"-4"
\ar@{->} "2rr";"rr"
\ar@{->} "rr";"k"
\ar@{->} "k-1";"rr"
\ar@{->} "3";"k-1"
\ar@{->} "3";"-rr"
\ar@{->} "-4";"rr"
\end{xy}
\]
implies that
\begin{equation}\label{D-lem2-pr5}
(\varphi^G_{\mu_{r}\mu_{r-1}\mu_{r-2}\cdots\mu_{k+1}\mu_k(\textbf{V})})_{r}(a;\textbf{Y})\\
=\frac{(\varphi^G_{\mu_{r-2}\cdots\mu_{k+1}\mu_k(\textbf{V})})_{r-2}+(\varphi^G_{\textbf{V}})_{k-1}(\varphi^G_{\textbf{V}})_{2r}(\varphi^G_{\textbf{V}})_{-(r-1)}}{(\varphi^G_{\textbf{V}})_{r}}.
\end{equation}
From (\ref{D-lem2-pr1-1}), by the same way as in (\ref{D-lem2-pr4}), we can prove the following:
 \begin{eqnarray*}
& &\hspace{-20pt}(\varphi^G_{\mu_{r-2}\cdots\mu_{k+1}\mu_k(\textbf{V})})_{r-2}(a;\textbf{Y})\\
&=&a^{\Lm_{k-1}+\Lm_{r-1}}Y_{2,k}Y_{1,r-1}(\varphi^G_{\textbf{V}})_{r}
(1+A^{-1}_{1,r-1}+A^{-1}_{1,r-1}A^{-1}_{1,r-2}+ \cdots+A^{-1}_{1,r-1}A^{-1}_{1,r-2}\cdots A^{-1}_{1,k+1})\qquad \\
& &-a^{\Lm_{r-1}+\Lm_r}Y_{1,r-1}Y_{2,r-1}
(\varphi^G_{\textbf{V}})_{k-1}.
\end{eqnarray*}
Substituting this for (\ref{D-lem2-pr5}), we obtain
\begin{eqnarray}
& &\hspace{-20pt}(\varphi^G_{\mu_r\mu_{r-1}\mu_{r-2}\cdots\mu_{k+1}\mu_k(\textbf{V})})_{r}(a;\textbf{Y})\nonumber\\
&=&a^{\Lm_{k-1}+\Lm_{r-1}}Y_{2,k}Y_{1,r-1}
(1+A^{-1}_{1,r-1}+A^{-1}_{1,r-1}A^{-1}_{1,r-2}+ \cdots+A^{-1}_{1,r-1}A^{-1}_{1,r-2}\cdots A^{-1}_{1,k+1}),\qquad \quad \nonumber
\end{eqnarray}
which means (\ref{D-lem4-cry}). \qed

\vspace{3mm}

\nd
[{\sl Proof of Theorem \ref{thmD} (v).}]

Using the induction on $(k-j)$, we shall prove (\ref{D-lem3-exp1}). First, let $k-j=1$. The vertices and arrows around $(\varphi_{\textbf{V}})_{k-1}$ in the mutation diagram of $\mu_{r}\mu_{r-1}\mu_{r-2}\cdots\mu_{k+1}\mu_{k}(\Sigma_0)$ are
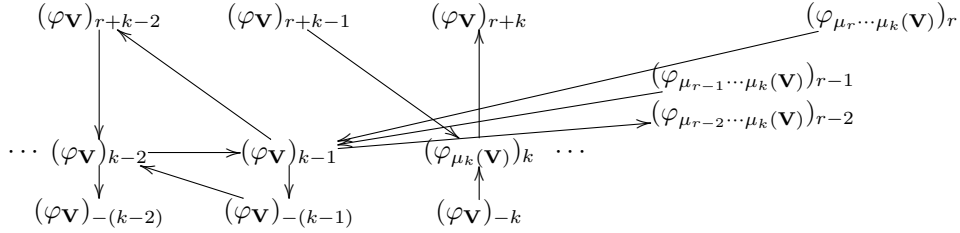
\[
\begin{xy}
(113,108) *{(\varphi_{\mu_{r}\cdots\mu_k(\textbf{V})})_{r}}="r",
(96,100) *{(\varphi_{\mu_{r-1}\cdots\mu_k(\textbf{V})})_{r-1}}="r-1",
(96,95) *{(\varphi_{\mu_{r-2}\cdots\mu_k(\textbf{V})})_{r-2}}="r-2",
(72,90)*{\cdots}="emp1",
(60,108)*{(\varphi_{\textbf{V}})_{r+k}}="r+k",
(60,90) *{(\varphi_{\mu_k(\textbf{V})})_{k}}="k",
(60,82)*{(\varphi_{\textbf{V}})_{-k}}="-j_k",
(35,108)*{(\varphi_{\textbf{V}})_{r+k-1}}="8",
(35,90) *{(\varphi_{\textbf{V}})_{k-1}}="4",
(35,82)*{(\varphi_{\textbf{V}})_{-(k-1)}}="-2",
(10,108)*{(\varphi_{\textbf{V}})_{r+k-2}}="6",
(10,90) *{(\varphi_{\textbf{V}})_{k-2}}="2",
(10,82)*{(\varphi_{\textbf{V}})_{-(k-2)}}="-1",
(0,90)*{\cdots}="emp",
\ar@{->} "r";"4"
\ar@{->} "r-1";"4"
\ar@{->} "4";"r-2"
\ar@{->} "k";"r+k"
\ar@{->} "6";"2"
\ar@{->} "-j_k";"k"
\ar@{->} "8";"k"
\ar@{->} "4";"-2"
\ar@{->} "2";"-1"
\ar@{->} "2";"4"
\ar@{->} "4";"6"
\ar@{->} "-2";"2"
\end{xy}
\]
means that
\begin{eqnarray}
& &\hspace{-20pt}
(\varphi^G_{\mu_{k-1}\mu_{r}\mu_{r-1}\cdots\mu_{k+1}\mu_k(\textbf{V})})_{k-1}(a;\textbf{Y})=  \frac{1}{(\varphi^G_{\textbf{V}})_{k-1}}\left(
(\varphi^G_{\mu_{r-2}\cdots\mu_{k+1}\mu_k(\textbf{V})})_{r-2}(\varphi^G_{\textbf{V}})_{r+k-2}(\varphi^G_{\textbf{V}})_{-(k-1)}\right.\nonumber\\
&+&\left.
(\varphi^G_{\textbf{V}})_{k-2}
(\varphi^G_{\mu_{r-1}\mu_{r-2}\cdots\mu_k(\textbf{V})})_{r-1}
(\varphi^G_{\mu_r\mu_{r-1}\cdots\mu_k(\textbf{V})})_{r}\right).\label{D-lem3-pr1}
\end{eqnarray}

The explicit form (\ref{D-lem2-pr2-1}) of $(\varphi^G_{\mu_{r-2}\cdots\mu_{k+1}\mu_k(\textbf{V})})_{r-2}$ can be rewritten as
\begin{eqnarray}
& &\hspace{-20pt}
(\varphi^G_{\mu_{r-2}\cdots\mu_{k+1}\mu_k(\textbf{V})})_{r-2}\nonumber\\
&=&a^{-\Lm_{k-1}}(\varphi^G_{\mu_{r-1}\mu_{r-2}\cdots\mu_k(\textbf{V})})_{r-1}
(\varphi^G_{\mu_r\mu_{r-1}\cdots\mu_k(\textbf{V})})_{r}\frac{1}{Y_{2,k}}\label{D-lem3-pr2}\\
&+&a^{\Lm_{r-1}+\Lm_r}Y_{1,r-2}(\varphi^G_{\textbf{V}})_{k-1}
(1+A^{-1}_{1,r-2}+A^{-1}_{1,r-2}A^{-1}_{1,r-3}+ \cdots+A^{-1}_{1,r-2}\cdots A^{-1}_{1,k+1})\nonumber
\end{eqnarray}
by using (\ref{iniD-exp1}), (\ref{D-lem2-pr3-1}) and (\ref{D-lem2-pr3-2}). Applying this and $(\varphi^G_{\textbf{V}})_{k-2}=a^{\Lm_{k-2}-\Lm_{k-1}}\frac{Y_{2,k-1}}{Y_{2,k}}(\varphi^G_{\textbf{V}})_{k-1}-a^{\Lm_{k-2}}\frac{Y_{1,k-1}Y_{2,k-1}}{Y_{2,k}}$ (\ref{twovarsD}) to (\ref{D-lem3-pr1}), we obtain
\begin{eqnarray*}
& &\hspace{-20pt}
(\varphi^G_{\mu_{k-1}\mu_{r}\mu_{r-1}\cdots\mu_{k+1}\mu_k(\textbf{V})})_{k-1}(a;\textbf{Y})\\
&=&a^{\Lm_{k-2}-\Lm_{k-1}}\frac{Y_{2,k-1}}{Y_{2,k}}(\varphi^G_{\mu_{r-1}\mu_{r-2}\cdots\mu_k(\textbf{V})})_{r-1}
(\varphi^G_{\mu_r\mu_{r-1}\cdots\mu_k(\textbf{V})})_{r}\\
&+&a^{\Lm_{k-2}+\Lm_{k-1}+\Lm_{r-1}+\Lm_r}Y_{1,k-1}Y_{2,k-1}Y_{1,r-2}
(1+A^{-1}_{1,r-2}+\cdots +A^{-1}_{1,r-2}\cdots A^{-1}_{1,k+1})\\
&=&a^{\Lm_{k-2}+\Lm_{k-1}+\Lm_{r-1}+\Lm_r} Y_{2,k-1}Y_{2,k}Y_{1,r-1}Y_{1,r}
(1+A^{-1}_{1,r-1}+\cdots +A^{-1}_{1,r-1}\cdots A^{-1}_{1,k+1})\\
& &\times (1+A^{-1}_{1,r}+A^{-1}_{1,r}A^{-1}_{1,r-2}+ \cdots+A^{-1}_{1,r}A^{-1}_{1,r-2}A^{-1}_{1,r-3}\cdots A^{-1}_{1,k+1})\\
&+&a^{\Lm_{k-2}+\Lm_{k-1}+\Lm_{r-1}+\Lm_r}Y_{1,k-1}Y_{2,k-1}Y_{1,r-2}
(1+A^{-1}_{1,r-2}+\cdots +A^{-1}_{1,r-2}\cdots A^{-1}_{1,k+1}).
\end{eqnarray*}

Next, let us consider the case $k-j>1$. The vertices and arrows around $(\varphi_{\textbf{V}})_{j}$ in the mutation diagram of $\mu_{j+1}\cdots \mu_{k-2}\mu_{k-1}\mu_{r}\mu_{r-1}\cdots\mu_{k+1}\mu_k(\Sigma_0)$ are as follows:
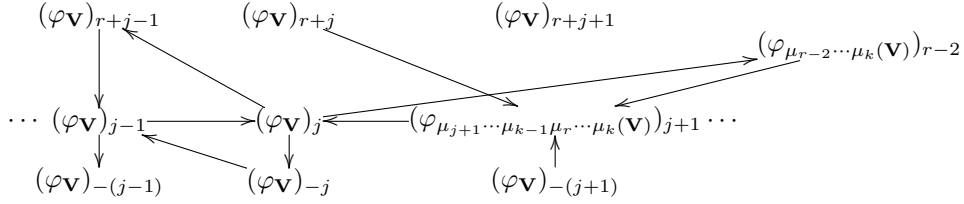
\[
\begin{xy}
(110,100) *{(\varphi_{\mu_{r-2}\cdots\mu_k(\textbf{V})})_{r-2}}="r-1",
(92,90)*{\cdots}="emp1",
(70,104)*{(\varphi_{\textbf{V}})_{r+j+1}}="r+k",
(70,90)*{ (\varphi_{\mu_{j+1}\cdots\mu_{k-1}\mu_{r}\cdots\mu_k(\textbf{V})})_{j+1}}="k",
(70,82)*{(\varphi_{\textbf{V}})_{-(j+1)}}="-j_k",
(35,104)*{(\varphi_{\textbf{V}})_{r+j}}="8",
(35,90) *{(\varphi_{\textbf{V}})_{j}}="4",
(35,82)*{(\varphi_{\textbf{V}})_{-j}}="-2",
(10,104)*{(\varphi_{\textbf{V}})_{r+j-1}}="6",
(10,90) *{(\varphi_{\textbf{V}})_{j-1}}="2",
(10,82)*{(\varphi_{\textbf{V}})_{-(j-1)}}="-1",
(0,90)*{\cdots}="emp",
\ar@{->} "4";"r-1"
\ar@{->} "k";"4"
\ar@{->} "r-1";"k"
\ar@{->} "-j_k";"k"
\ar@{->} "8";"k"
\ar@{->} "4";"-2"
\ar@{->} "2";"-1"
\ar@{->} "2";"4"
\ar@{->} "4";"6"
\ar@{->} "6";"2"
\ar@{->} "-2";"2"
\end{xy}
\]
This diagram says that
\begin{eqnarray}
& &\hspace{-20pt}
(\varphi^G_{\mu_j\mu_{j+1}\cdots \mu_{k-1}\mu_{r}\mu_{r-1}\cdots\mu_k(\textbf{V})})_{j}(a;\textbf{Y})\nonumber\\
&=& \frac{(\varphi^G_{\mu_{r-2}\cdots\mu_{k+1}\mu_k(\textbf{V})})_{r-2}(\varphi^G_{\textbf{V}})_{r+j-1}(\varphi^G_{\textbf{V}})_{-j}
+(\varphi^G_{\mu_{j+1}\cdots \mu_{k-1}\mu_{r}\mu_{r-1}\cdots\mu_k(\textbf{V})})_{j+1}
(\varphi^G_{\textbf{V}})_{j-1}}{(\varphi^G_{\textbf{V}})_{j}}. \qquad \quad \label{D-lem3-pr3}
\end{eqnarray}
It follows by induction hypothesis that
\begin{eqnarray*}
& &\hspace{-20pt}
(\varphi^G_{\mu_{j+1}\cdots \mu_{k-1}\mu_{r}\mu_{r-1}\cdots\mu_k(\textbf{V})})_{j+1}\\
&=&a^{\Lm_{j}-\Lm_{k-1}}(\varphi^G_{\mu_{r-1}\mu_{r-2}\cdots\mu_k(\textbf{V})})_{r-1}
(\varphi^G_{\mu_r\mu_{r-1}\cdots\mu_k(\textbf{V})})_{r}\frac{Y_{2,j+1}}{Y_{2,k}}\\
&+&a^{\Lm_{j}+\Lm_{k-1}+\Lm_{r-1}+\Lm_r} Y_{2,j+1}Y_{1,k-1}Y_{1,r-2}
(1+A^{-1}_{1,r-2}+\cdots+A^{-1}_{1,r-2}A^{-1}_{1,r-3}\cdots A^{-1}_{1,k+1})\\
& &\times (1+A^{-1}_{1,k-1}+\cdots+A^{-1}_{1,k-1}A^{-1}_{1,k-2}\cdots A^{-1}_{1,j+2}).
\end{eqnarray*}
Furthermore, we know that $(\varphi^G_{\textbf{V}})_{j-1}=a^{\Lm_{j-1}-\Lm_j}\frac{Y_{2,j}}{Y_{2,j+1}}(\varphi^G_{\textbf{V}})_{j}-a^{\Lm_{j-1}}\frac{Y_{1,j}Y_{2,j}}{Y_{2,j+1}}$. Thus, taking (\ref{D-lem3-pr2}) and (\ref{D-lem3-pr3}) into account, we obtain
\begin{eqnarray*}
& &\hspace{-20pt}
(\varphi^G_{\mu_j\mu_{j+1}\cdots \mu_{k-1}\mu_{r}\mu_{r-1}\cdots\mu_k(\textbf{V})})_{j}(a;\textbf{Y})=a^{\Lm_{j-1}-\Lm_{j}}\frac{Y_{2,j}}{Y_{2,j+1}}
(\varphi^G_{\mu_{j+1}\cdots \mu_{k-1}\mu_{r}\mu_{r-1}\cdots\mu_k(\textbf{V})})_{j+1}\\
& &
+\frac{a^{\Lm_{j-1}+\Lm_j+\Lm_{k-1}+\Lm_{r-1}+\Lm_r}}{(\varphi^G_{\textbf{V}})_{j}}
Y_{1,j}Y_{2,j}Y_{1,r-2}Y_{1,k-1}\\
&\times& (1+A^{-1}_{1,r-2}+\cdots+A^{-1}_{1,r-2}\cdots A^{-1}_{1,k+1})
A^{-1}_{1,k-1}\cdots A^{-1}_{1,j+1}(1+A^{-1}_{1,j}+\cdots+A^{-1}_{1,j}\cdots A^{-1}_{1,1})\\
&=& a^{\Lm_{j-1}+\Lm_{k-1}+\Lm_{r-1}+\Lm_r} Y_{2,j}Y_{2,k}Y_{1,r-1}Y_{1,r}
(1+A^{-1}_{1,r-1}+\cdots+A^{-1}_{1,r-1}A^{-1}_{1,r-2}\cdots A^{-1}_{1,k+1})\\
& &\times 
(1+A^{-1}_{1,r}+A^{-1}_{1,r}A^{-1}_{1,r-2}+ \cdots+A^{-1}_{1,r}A^{-1}_{1,r-2}A^{-1}_{1,r-3}\cdots A^{-1}_{1,k+1})\\
&+&a^{\Lm_{j-1}+\Lm_{k-1}+\Lm_{r-1}+\Lm_r} Y_{2,j}Y_{1,k-1}Y_{1,r-2}
(1+A^{-1}_{1,r-2}+\cdots+A^{-1}_{1,r-2}A^{-1}_{1,r-3}\cdots A^{-1}_{1,k+1})\qquad \quad \\
& &\times (1+A^{-1}_{1,k-1}+\cdots+A^{-1}_{1,k-1}A^{-1}_{1,k-2}\cdots A^{-1}_{1,j+1}),
\end{eqnarray*}
which is our desired result. The description (\ref{D-lem3-cry1}) immediately follows from (\ref{D-lem3-exp1}). \qed


\begin{thebibliography}{9}
\bibitem{A-F-Z}\label{A-F-Z}A. Berenstein, S. Fomin,  
A. Zelevinsky, Cluster algebras III : 
Upper bounds and double Bruhat cells, 
Duke Math Journal, vol. 126 No.1, 1--52 (2005).

\bibitem{B-Z}\label{B-Z}A. Berenstein, A. Zelevinsky, 
Tensor product multiplicities,
canonical bases and totally positive varieties, 
Invent. Math. 143 No.1, 77--128 (2001). 

\bibitem{F-Z}\label{F-Z}S. Fomin, A. Zelevinsky, 
Double Bruhat cells and total positivity, 
J. Amer. Math. Soc., vol.12, No.2, 335--380 (1999).

\bibitem{FZ2}S. Fomin,
A. Zelevinsky, Cluster algebras I: Foundations, J. Amer. Math. Soc., vol.15,
No.2, 497--529 (2002).

\bibitem{FZ3}S. Fomin,
A. Zelevinsky, Cluster algebras II: Finite type classification, Invent. Math. 154 No.1, 63--121 (2003). 

\bibitem{GLS}C. Geiss, B. Leclerc, J. Schr$\ddot{{\rm o}}$er, Kac-Moody groups and cluster algebras, Adv. Math. 228 No.1, 329--433 (2011).

\bibitem{M-M-A}\label{M-M-A}M. Gekhtman, 
M. Shapiro, A. Vainshtein, 
{\it Cluster Algebras and Poisson Geometry}, AMS (2010). 

\bibitem{GY}\label{GY}K. R. Goodearl, M. T. Yakimov, The Berenstein-Zelevinsky quantum cluster algebra conjecture, arXiv:1602.00498.

\bibitem{GY2}\label{GY2}Goodearl K.R., Yakimov, M.T.
Quantum cluster algebras and quantum nilpotent algebras.
Proc. Natl. Acad. Sci. USA 111 (2014), no. 27, 9696--9703. 

%\bibitem{H-K}\label{H-K} Hong J. and Kang S-J., 
%Quantum Groups and Crystal Bases, 
%Graduate Studies in Mathematics vol.42, AMS, (2002). 

\bibitem{KaN}Y. Kanakubo, T. Nakashima,
Cluster Variables on Certain Double Bruhat Cells of 
Type $(u,e)$ and Monomial Realizations of Crystal Bases 
of Type A, SIGMA, vol.11, 1--32 (2015).

\bibitem{KaN2}Y. Kanakubo, T. Nakashima,
Cluster Variables on Double Bruhat Cells $G^{u,e}$ of Classical
Groups and Monomial Realizations of Demazure Crystals, arXiv:1604.05956.

\bibitem{KaN3}Y. Kanakubo, T. Nakashima, Cluster algebras of finite type via a Coxeter element and Demazure Crystals of type A, arXiv:1703.08323.

\bibitem{K}\label{K}M. Kashiwara, Realizations of crystals, 
Combinatorial and geometric representation theory, Contemporary Mathematics 325, AMS, 133--139 (2003).

\bibitem{K0}M. Kashiwara, Crystallizing the $q$-analogue of universal 
              enveloping algebras, Comm. Math. Phys, vol.133, 249--260 (1990).

\bibitem{K2}\label{K2}M. Kashiwara, On crystal bases of 
the q-analogue of universal enveloping algebras, 
Duke Math Journal vol.63, No.2, 465--516 (1991).

\bibitem{K3}\label{K3}M. Kashiwara, Bases cristallines des groupes
quantiques, edited by Charles Cochet. Cours Specialises, 9, 
Societe Mathematique de France, Paris, (2002).

\bibitem{KN}M. Kashiwara, T. Nakashima, 
        Crystal graphs for representations 
       of the $q$-analogue of classical Lie algebras, 
       {\sl J. Algebra},vol.165, No.2, 295--345 (1994).

\bibitem{Nj}\label{Nj}H. Nakajima, t-analogs of q-characters of quantum affine algebras of 
type $A_n$, $D_n$, Contemp. Math, 325, AMS, Providence, RI, 141--160 (2003).

\bibitem{N1}T. Nakashima, Decorations on Geometric Crystals and 
Monomial Realizations of Crystal Bases for Classical Groups, 
{\sl J. Algebra}, vol.399, 712--769 (2014).

\end{thebibliography}
\end{document}